\documentclass[english,journal]{IEEEtran}

\usepackage[T1]{fontenc}
\usepackage[latin9]{inputenc}
\usepackage{array}
\usepackage{float}
\usepackage{booktabs}
\usepackage{units}
\usepackage{url}
\usepackage{amsmath}
\usepackage{amssymb}
\usepackage{graphicx}
\usepackage{esint}

\makeatletter

\providecommand{\tabularnewline}{\\}
\newcommand{\lyxdot}{.}

\floatstyle{ruled}
\newfloat{algorithm}{tbp}{loa}
\providecommand{\algorithmname}{Algorithm}
\floatname{algorithm}{\protect\algorithmname}

\usepackage{cite}

\usepackage{algorithm}
\usepackage{algpseudocode}

\DeclareMathOperator*{\argmin}{arg\,min}

\DeclareMathOperator*{\st}{s.t.\,}
\DeclareMathOperator{\prox}{prox}
\DeclareMathOperator{\diag}{diag}
\DeclareMathOperator{\sign}{sign}

\usepackage{array}

\usepackage{babel}

\newcommand{\etal}{\emph{et al.\ }}

\@ifundefined{showcaptionsetup}{}{%
 \PassOptionsToPackage{caption=false}{subfig}}
\usepackage{subfig}
\makeatother

\usepackage{babel}
\begin{document}

\title{TRex: A Tomography Reconstruction Proximal Framework for Robust Sparse
View X-Ray Applications}

\author{Mohamed Aly\thanks{M Aly is with the Visual Computing Center, KAUST, KSA and is on leave from Computer Engineering, Cairo University, Egypt},
Guangming Zang\thanks{G Zang is with the Visual Computing Center, KAUST, KSA},
Wolfgang Heidrich\thanks{W Heidrich is with the Visual Computing Center, KAUST, KSA},
and Peter Wonka\thanks{P Wonka is with the Visual Computing Center, KAUST, KSA}}
\maketitle
\begin{abstract}
We present TRex, a flexible and robust Tomographic Reconstruction
framework using proximal algorithms. We provide an overview and perform
an experimental comparison between the famous iterative reconstruction
methods in terms of reconstruction quality in sparse view situations.
We then derive the proximal operators for the four best methods. We
show the flexibility of our framework by deriving solvers for two
noise models: Gaussian and Poisson; and by plugging in three powerful
regularizers. We compare our framework to state of the art methods,
and show superior quality on both synthetic and real datasets. 
\end{abstract}

\begin{IEEEkeywords}
Image reconstruction, X-ray imaging and computed tomography, Simultaneous
Algebraic Reconstruction Technique, SART, Proximal Algorithms, Cone
beam X-ray tomography 
\end{IEEEkeywords}

\newcommand{\Short}[1]{}
\newcommand{\Long}[1]{#1}

\section{Introduction\label{sec:Introduction}}

Reducing the dosage in X-ray tomography is a very important issue
in medical applications, since long term exposure to X-rays can have
adverse health effects. This can be done in at least two ways: (a)
reducing the X-ray beam power, which leads to increased measurement
noise at the detectors; or (b) acquiring fewer projections to reduce
the acquisition time \cite{herman2009fundamentals}. This makes the
reconstruction problem even more ill-posed, since less information
is collected from the volume to be reconstructed; and one has to use
non-linear regularizers (priors) to achieve a reasonable result. This
is typically done using iterative solvers \cite{thibault2007three,zhang2014model}.

Iterative algorithms for X-ray tomography reconstruction have been
around for years. In fact, one of the first implemented tomography reconstruction algorithm
was an iterative one \cite{kak2001principles,gordon1970algebraic,gordon1971reconstruction}. However, non-iterative, transform-based algorithms, such as the filtered
back projection (FBP) \cite{ramachandran1971three,shepp1974fourier,feldkamp1984practical},
have been more popular due to their speed and low computational cost.
Moreover, most commercial X-ray CT scanners employ some variant of
FBP in their reconstruction software \cite{pan2009commercial}. Recently,
interest has been ignited again in iterative algorithms because, although
they are more computationally demanding, they are much more flexible
and yield superior reconstruction quality by employing powerful priors. 

Thus, in this work, we study iterative reconstruction techniques.
We present TRex, a flexible proximal framework for robust X-Ray tomography
reconstruction in sparse view applications. TRex uses iterative algorithms,
especially the SART (Simultaneous ART) \cite{andersen1984simultaneous,andersen1989algebraic},
to solve the tomography proximal operator. We show that they are better
suited for this task and produce better performance than state of
the art, combined with different noise models in the data terms and
with different powerful regularizers. Up to our knowledge, this is
the first time these methods have been used to directly solve the
tomography proximal operator.

We start by conducting a thorough comparison of the famous iterative
algorithms including SART \cite{andersen1984simultaneous}, ART (Algebraic
Reconstruction Technique) \cite{andersen1984simultaneous}, SIRT (Simultaneous
Iterative Reconstruction Technique) \cite{gilbert1972iterative},
BSSART (Block Simplified SART) \cite{censor2002block}, BICAV (Block
Iterative Component Averaging) \cite{censor2001bicav}, Conjugate
Gradient (CG) \cite{bjorck1996numerical}, and OS-SQS (Ordered Subset-Separable
Quadratic Surrogates) \cite{depierro1994modified,hudson1994accelerated,erdogan1999ordered,kim2013accelerating,nien2015fast}.
We establish that SART provides the best performance in the sparse
view measurements situations, followed closely by ART, OS-SQS, and
BICAV.

We then describe our framework, TRex, which is based on using proximal
algorithms \cite{boyd2011distributed,parikh2013proximal} together
with these iterative methods. We derive proximal operators for SART,
ART, BICAV, and OS-SQS. We show how to use these proximal operators
to minimize two data fitting terms: (a) least squares (LS) that assumes
a Gaussian noise model; and (b) weighted least squares (WLS) that
assumes an approximation to a Poisson noise model \cite{clinthorne1993preconditioning}.
We also show how to plug in different powerful regularizers; namely
Isotropic Total Variation (ITV) \cite{rudin1992nonlinear}, Anisotropic
Total Variation (ATV) \cite{sidky2012convex}, and Sum of Absolute
Differences (SAD) \cite{gregson2012stochastic}. We perform thorough
comparisons between the different proximal operators, data terms,
and regularizers using real and synthetic data.

Finally, we compare our framework to state of the art methods, namely
the ADMM method from Ramani \etal \cite{ramani2012splitting} and
the OS-SQS method (with and without momentum method) from Kim \etal
\cite{kim2015combining}, and show that our framework gives superior
reconstruction quality. Please consult \cite{aly2016tomography} for
further details, expanded experiments, and more results.

In summary, we provide the following contributions:
\begin{enumerate}
\item We present TRex, a flexible proximal reconstruction framework that
relies on iterative methods for directly solving the tomography proximal
operator.
\item We perform a thorough experimental comparison of famous iterative
reconstruction methods on synthetic and real datasets.
\item We derive proximal operators for SART, ART, BICAV, and OS-SQS; and
compare them.
\item We derive solvers for different data terms assuming different noise
models, namely Gaussian and Poisson models, using the derived proximal
operators, and show how to use our framework with different powerful
regularizers.
\item We compare our framework to state of the art methods and show that
it produces superior reconstructions.
\item We make our code---which is based on the ASTRA toolbox \cite{van2015astra}---and
data publicly available at \url{https://github.com/mohamedadaly/TRex}.
\end{enumerate}
This paper is organized as follows. In Sec. \ref{sec:Related-Work}
we present related work. An overview of the famous iterative algorithms
is detailed in Sec. \ref{sec:Iterative-Algorithms}. The different
proximal operators are derived in Sec. \ref{sec:Proximal-Operators}.
The TRex framework is explained in Sec. \ref{sec:Proximal-Framework},
where we show the general algorithm together with three regularizers
and two data terms. The experiments and datasets are presented in
Sec. \ref{sec:Experiments}, and finally the conclusions are given
in Sec. \ref{sec:Discussion-and-Conclusion}.

\section{Related Work\label{sec:Related-Work}}

There are two general approaches for X-ray tomography reconstruction:
transform-based methods and iterative methods \cite{kak2001principles,herman2009fundamentals}.
Transform methods rely on the Radon transform and its inverse introduced
in 1917. The most widely used reconstruction method is the Filtered
Backprojection (FBP) algorithm introduced \cite{herman2009fundamentals,kak2001principles}. Transform methods are usually viewed as much faster than iterative
methods, and have therefore been the method of choice for X-ray scanner
manufacturers \cite{pan2009commercial}.

Iterative methods, on the other hand, use algebraic techniques to
solve the reconstruction problem. They generally model the problem
as a linear system and solve it using established numerical methods
\cite{herman2009fundamentals}. ART, and its many variants, are among
the best known iterative reconstruction algorithms  \cite{gordon1970algebraic,lent1977convergent,shepp1982maximum,censor1983finite,andersen1984simultaneous,andersen1989algebraic}.
They use variations of the projection method of Kaczmarz \cite{kaczmarz1937angenaherte}
and have modest memory requirements, and have been shown to yield
better reconstruction results than transform methods. They are matrix
free, and work without having to explicitly store the system matrix.
OS-SQS and related methods \cite{depierro1994modified,erdogan1999ordered,kim2013accelerating}
are closely related to ART and have similar properties to SIRT \cite{gregor2015comparison}.
They have also been shown \cite{kim2015combining} to be accelerated
using momentum techniques.

Iterative methods provide more flexibility in incorporating prior
information into the reconstruction process. For example, instead of assuming a Gaussian noise model and minimizing
a least squares data term, one can easily use iterative methods with
other noise models, such as the Poisson noise model \cite{clinthorne1993preconditioning,depierro1994modified,elbakri2002statistical,wang2006penalized,thibault2007three}
that boils down to solving WLS problem instead. Priors are also easy to use with iterative methods. For example,
the Total Variation \cite{rudin1992nonlinear} prior has been used
for tomography reconstruction \cite{sidky2008image,mory2012ecg}.

Proximal algorithms have been widely used in many problems in machine
learning and signal processing \cite{bauschke2011convex,combettes2011proximal,boyd2011distributed,parikh2013proximal}.
They have also been used in tomography reconstruction \cite{mory2012ecg,sidky2008image}.
For example, \cite{mory2012ecg} used the Alternating Direction Method
of Multipliers (ADMM) \cite{boyd2011distributed} with total variation
prior, where the data term was optimized using CG \cite{bjorck1996numerical}.
\cite{sidky2012convex} discussed using the Chambolle-Pock algorithm
\cite{chambolle2011first} for tomography reconstruction with different
priors. \cite{ramani2012splitting} used ADMM with Preconditioned CG (PCG)
\cite{fessler1999conjugate} for optimizing the weighted least squares
data term. \cite{nien2015fast} used Linearized ADMM \cite{parikh2013proximal}
(also known as Inexact Split Uzawa \cite{esser2010general}) with
Ordered Subset-based methods \cite{erdogan1999ordered} for optimizing
the data term and FISTA \cite{beck2009fast} for optimizing the prior
term. However, none of these methods used the iterative algorithms
we study in this work as their data term solver, which provides superior
reconstruction as we will show.

There are currently a number of open source software packages for
tomography reconstruction. SNARK09 \cite{klukowska2013snark09} is one of the oldest. The Reconstruction
ToolKit (RTK) \cite{rit2014reconstruction} is a high performance
C++ toolkit focusing on 3D cone beam reconstruction that is based
on the image processing package Insight ToolKit (ITK). It includes
implementations of several algorithms, including FDK, SART, and an
ADMM TV-regularized solver with CG \cite{mory2012ecg}. The ASTRA toolbox \cite{van2015astra} is a Matlab-based GPU-accelerated
toolbox for tomography reconstruction. It includes implementations
of several algorithms, including SART, SIRT, FBP, among others. We
modify and extend ASTRA to implement our algorithms and generate the
experiments in this work.

\section{Iterative Algorithms\label{sec:Iterative-Algorithms}}

\begin{algorithm}[h]
\protect\protect\caption{\label{alg:Iterative-Algorithm}Outline of Iterative Algorithms}

\begin{algorithmic}[1]

\Require $A\in\mathbb{R}{}^{m\times n}$, $\alpha\in\mathbb{R}$,
$p\in\mathbb{R}{}^{m}$

\State Initialize: $x^{(0)}$

\ForAll {$t=1\ldots T$ }

\ForAll {subsets $S\in\mathcal{S}$ } 

\State$x^{(t+1)}=x^{(t)}+\alpha\Delta x^{(t)}$

\State$x^{(t+1)}=\mbox{clip}(x^{(t+1)})$

\EndFor

\EndFor

\Return volume reconstruction $x\in\mathbb{R}{}^{n}$

\end{algorithmic} 
\end{algorithm}

The tomography problem can be represented as solving a linear system
\cite{kak2001principles,herman2009fundamentals} 
\begin{equation}
Ax=p,\label{eq:linear-system}
\end{equation}
where $x\in\mathbb{R}^{n}$ is the unknown volume in vector form,
$A\in\mathbb{R}^{m\times n}$ is the projection system matrix, and
$p\in\mathbb{R}^{m}$ represents the measured line projections (sinogram).
The iterative algorithms that we study in this work all have the same
general outline in Alg. \ref{alg:Iterative-Algorithm}, but differ
in the update formula in step 4. The subset $S$ in step 3 can be
only 1 projection ray as in ART i.e. there are $m$ subsets $S_{i}=\{i\,|\, i=1\ldots m\}$;
can contain all the rays in a projection view as in SART i.e. there
are $m/s$ subsets where $s$ is the number of projection views; or
can contain the whole projection rays as in SIRT i.e. there is only
one subset $S=\{1,\ldots,m\}$. Step 5 clips the negative values of
the volume, which is assumed to be non-negative.

The update step $\Delta x^{(t)}$ is typically a function of (a subset
of) the forward projection error $p_{S}-A_{S}x^{(t)}$ that is then
back projected with some normalization procedure. It can take the form 
\[
\Delta x^{(t)}=\Phi\left(A_{S}^{T},p_{S}-A_{S}x^{(t)}\right)
\]
where the function $\Phi(\cdot)$ computes the required update, $A_{S}$
contains a subset of the rows of $A$, similarly for $p_{S}$--please
see below. This can be seen as an approximation to the actual gradient
$A^{T}(p-Ax)$ of the least square objective 
\[
\argmin_{x}\Vert Ax-p\Vert_{2}^{2}
\]
 and so these algorithms can be viewed as variations of (stochastic)
gradient descent \cite{kim2015combining} where they differ on how
they approximate the gradient. We also notice that the inner loop
in step 3 for all these algorithms takes roughly the same time, since
it involves one full sweep over the rows of $A$.

Below we quickly review the different methods, and Table \ref{tab:Iterative-Methods}
provides a summary of their important properties. 

{\footnotesize{}}
\begin{table*}
{\footnotesize{}\center
\renewcommand*\arraystretch{1.5}}%
\begin{tabular}{cc>{\centering}m{0.03\textwidth}cc}
\toprule 
Method & Update Step & Subset & Solved Problem & Converges\tabularnewline
\midrule
\midrule 
ART \cite{gordon1970algebraic} & $\begin{aligned}x_{j}^{(t+1)} & =x_{j}^{(t)}+\alpha\frac{p_{i}-\sum_{k}a_{ik}x_{k}^{(t)}}{\sum_{k}a_{ik}^{2}}a_{ij}\\
x^{(t+1)} & =x^{(t)}+\alpha A_{i}^{T}R^{-1}\left(p_{i}-A_{i}x^{(t)}\right)
\end{aligned}
$ & one ray & $\begin{aligned}x^{\star}= & \argmin_{x}\Vert x\Vert_{2}^{2}\\
 & \st Ax=p
\end{aligned}
$ & Yes\tabularnewline
\midrule 
SIRT \cite{gilbert1972iterative} & $\begin{aligned}x_{j}^{(t+1)} & =x_{j}^{(t)}+\alpha\frac{1}{\sum_{i=1}^{m}a_{ij}}\sum_{i=1}^{m}\frac{p_{i}-\sum_{k=1}^{n}a_{ik}x_{k}^{(t)}}{\sum_{k=1}^{n}a_{ik}}a_{ij}\\
x^{(t+1)} & =x^{(t)}+\alpha C^{-1}A^{T}R^{-1}\left(p-Ax^{(t)}\right)
\end{aligned}
$ & all rays & $x^{\star}=\argmin_{x}\Vert Ax-p\Vert_{R^{-1}}^{2}$ & Yes\tabularnewline
\midrule 
SART \cite{andersen1984simultaneous} & $\begin{aligned}x_{j}^{(t+1)} & =x_{j}^{(t)}+\alpha\frac{1}{\sum_{i\in S}a_{ij}}\sum_{i\in S}\frac{p_{i}-\sum_{k=1}^{n}a_{ik}x_{k}^{(t)}}{\sum_{k=1}^{n}a_{ik}}a_{ij}\\
x^{(t+1)} & =x^{(t)}+\alpha C_{S}^{-1}A_{S}^{T}R^{-1}\left(p-A_{S}x^{(t)}\right)
\end{aligned}
$ & one view & $\begin{aligned}x^{\star}\approx & \argmin_{x}\Vert x\Vert_{2}^{2}\\
 & \st Ax=p
\end{aligned}
$ & No\tabularnewline
\midrule 
BSSART \cite{censor2002block} & $\begin{aligned}x_{j}^{(t+1)} & =x_{j}^{(t)}+\alpha\frac{1}{\sum_{i=1}^{m}a_{ij}}\sum_{i\in S}\frac{p_{i}-\sum_{k=1}^{n}a_{ik}x_{k}^{(t)}}{\sum_{k=1}^{n}a_{ik}}a_{ij}\\
x^{(t+1)} & =x^{(t)}+\alpha C^{-1}A_{S}^{T}R^{-1}\left(p_{S}-A_{S}x^{(t)}\right)
\end{aligned}
$ & one view & $\begin{aligned}x^{\star}= & \argmin_{x}\Vert x\Vert_{2}^{2}\\
 & \st Ax=p
\end{aligned}
$ & Yes\tabularnewline
\midrule 
BICAV \cite{censor2001bicav} & $\begin{aligned}x_{j}^{(t+1)} & =x_{j}^{(t)}+\alpha\frac{1}{\sum_{i\in S}\{a_{ij}\ne0\}}\sum_{i\in S}\frac{p_{i}-\sum_{k=1}^{n}a_{ik}x_{k}^{(t)}}{\sum_{k=1}^{n}a_{ik}^{2}}a_{ij}\\
x^{(t+1)} & =x^{(t)}+\alpha C_{S}^{-1}A_{S}^{T}R^{-1}\left(p_{S}-A_{S}x^{(t)}\right)
\end{aligned}
$ & one view & $\begin{aligned}x^{\star}= & \argmin_{x}\Vert x\Vert_{2}^{2}\\
 & \st Ax=p
\end{aligned}
$ & Yes\tabularnewline
\midrule 
OS-SQS \cite{erdogan1999ordered} & $\begin{aligned}x_{j}^{(t+1)} & =x_{j}^{(t)}+\frac{\alpha s}{\left(\sum_{k=1}^{m}a_{kj}\sum_{i=1}^{n}a_{ki}\right)}\sum_{i\in S}\left(p_{i}-\sum_{k=1}^{n}a_{ik}x_{k}^{(t)}\right)a_{ij}\\
x^{(t+1)} & =x^{(t)}+\alpha sC^{-1}A_{S}^{T}\left(p_{S}-A_{S}x^{(t)}\right)
\end{aligned}
$ & one view & $x^{\star}\approx\Vert Ax-p\Vert_{2}^{2}$ & No\tabularnewline
\midrule 
CGLS \cite{bjorck1996numerical,fessler1999conjugate} & $x^{(t+1)}=x^{(t)}+\alpha_{t}\Phi\left(A^{T}(p-Ax^{(t)})\right)$ & all rays & $x^{\star}=\Vert Ax-p\Vert_{2}^{2}$ & Yes\tabularnewline
\bottomrule
\end{tabular}

{\footnotesize{}\renewcommand*\arraystretch{1}}{\footnotesize \par}

{\footnotesize{}\protect\caption{Summary of iterative methods and their properties.\label{tab:Iterative-Methods}
The first line in the update step is voxel-based, while the second
is the matrix formulation. See Sec. \ref{sec:Iterative-Algorithms}
 for details.}
}
\end{table*}
{\footnotesize \par}

\subsubsection{ART}

\cite{gordon1970algebraic,gordon1971reconstruction} is the first
algebraic method, and is based on Kaczmarz alternating projection
algorithm \cite{kaczmarz1937angenaherte}. ART treats each row of
$A$ in turn, and updates the current estimate according to 
\[
x_{j}^{(t+1)}=x_{j}^{(t)}+\alpha\frac{p_{i}-\sum_{k}a_{ik}x_{k}^{(t)}}{\sum_{k}a_{ik}^{2}}a_{ij}\mbox{ for }i=1\ldots m,
\]
 where $x_{j}^{(t)}$ is the $j$th voxel at time $t$, $a_{ij}$
is the entry in the $i$th row and $j$th column of $A$ and $\alpha\in\mathbb{R}$
is the relaxation parameter. This update is performed once for each
row of $A$, and one iteration includes a full pass over all the $m$
rows. The term $\sum_{k}a_{ik}x_{k}^{(t)}$ is the forward projection of
the volume estimate for the $i$th ray (equation or row), the difference
in the numerator is the projection error, that is then back projected
by multiplying the transpose of the $i$th row. It has been shown that ART converges to a least-norm solution to the
consistent system of equations \cite{tanabe1971projection} i.e. it
solves 
\begin{equation}
x^{\star}=\argmin_{x}\Vert x\Vert_{2}^{2}\st Ax=p.\label{eq:least-norm-problem}
\end{equation}
 In matrix notation, this can be also expressed as  
\[
x^{(t+1)}=x^{(t)}+\alpha A_{i}^{T}R^{-1}\left(p_{i}-A_{i}x^{(t)}\right)
\]
 where $A_{i}\in\mathbb{R}^{n}$ is the $i$th row of $A$ and $R\in\mathbb{R}^{m\times m}=\diag(r_{i})$
is a diagonal matrix where $r_{i}=\sum_{j}a_{ij}^{2}=\Vert A_{i}\Vert_{2}^{2}$
is the squared-norm of the $i$th row $A_{i}$.

\subsubsection{SIRT}

\cite{gilbert1972iterative} performs the updates \emph{simultaneously
}i.e. updates the volume once instead of updating it per each row
$A_{i}$. The update equation becomes 
\begin{eqnarray*}
x_{j}^{(t+1)} & = & x_{j}^{(t)}+\alpha\frac{1}{\sum_{i=1}^{m}a_{ij}}\sum_{i=1}^{m}\frac{p_{i}-\sum_{k=1}^{n}a_{ik}x_{k}^{(t)}}{\sum_{k=1}^{n}a_{ik}}a_{ij}.
\end{eqnarray*}
 In matrix form this becomes  
\[
x^{(t+1)}=x^{(t)}+\alpha C^{-1}A^{T}R^{-1}\left(p-Ax^{(t)}\right)
\]
 where $C\in\mathbb{R}^{n\times n}=\diag(c_{j})$ is a diagonal matrix
where $c_{j}=\sum_{i}a_{ij}$ is the sum of column $j$ of $A$ and
$R=\diag(r_{i})$ where $r_{i}=\sum_{j}a_{ij}$ is the sum of row
$i$ of $A$. In each iteration, SIRT performs a full forward projection
$Ax^{(t)}$, computes the residual, and then back projects it. The
diagonal matrices $R$ and $C$ perform scaling for the relevant entries.
It has been shown \cite{censor2002block,jiang2003convergence} that
SIRT converges to a solution of the WLS problem 
\begin{equation}
x^{\star}=\argmin_{x}\Vert Ax-p\Vert_{R^{-1}}^{2}=\min_{x}(Ax-p)^{T}R^{-1}(Ax-p)\label{eq:SIRT-WLS-problem}
\end{equation}
for $0<\alpha<2$. SIRT has been shown \cite{gregor2015comparison}
to be closely related, and in fact quite equivalent in terms of convergence
properties, to the OS-SQS method. It has also been shown to converge
best for $\alpha=2-\epsilon$ for a small $0<\epsilon\ll1$.

\subsubsection{SART}

\cite{andersen1984simultaneous} is a tradeoff between ART and SIRT,
in that it updates the volume after processing all the rows in a particular
projection view. The update equation becomes  
\begin{eqnarray*}
x_{j}^{(t+1)} & = & x_{j}^{(t)}+\alpha\frac{1}{\sum_{i\in S}a_{ij}}\sum_{i\in S}\frac{p_{i}-\sum_{k=1}^{n}a_{ik}x_{k}^{(t)}}{\sum_{k=1}^{n}a_{ik}}a_{ij}
\end{eqnarray*}
for $S\in\mathcal{S}$ where the summation $i\in S$ is across all
rows (rays) in projection view $S$ for all views $\mathcal{S}$.
This has been shown to provide faster convergence than ART and better
reconstruction results than SIRT \cite{mueller2000rapid,mueller1999fast}.
In matrix form it becomes  
\[
x^{(t+1)}=x^{(t)}+\alpha C_{S}^{-1}A_{S}^{T}R^{-1}\left(p_{S}-A_{S}x^{(t)}\right)
\]
 where $A_{S}\in\mathbb{R}^{s\times n}$ contains the $s$ rows in
projection $S$, $p_{S}$ contains the corresponding $s$ rays from
the projection measurements, $R$ contains the row sums as in SIRT,
while $C_{S}=\diag(c_{j}^{S})$ contains the column sums restricted
to the rows in $S$ i.e. $c_{j}^{S}=\sum_{i\in S}a_{ij}$. There is
still no proof of convergence for SART in the literature, but there
are proofs for variants of SART, such as BSSART and BICAV below, that
converge to a minimum-norm solution like ART. This motivates us to
assume that SART solves approximately the least norm problem in Eq.
\ref{eq:least-norm-problem}.

\subsubsection{BSSART}

\cite{censor2002block} is a slight simplification of SART, where
the column sums in the update equation are done over \emph{all} the
rows of $A$ instead of just over the rows in the current view, which
is quite similar to SIRT. The update equation becomes  
\begin{eqnarray*}
x_{j}^{(t+1)} & = & x_{j}^{(t)}+\alpha\frac{1}{\sum_{i=1}^{m}a_{ij}}\sum_{i\in S}\frac{p_{i}-\sum_{k=1}^{n}a_{ik}x_{k}^{(t)}}{\sum_{k=1}^{n}a_{ik}}a_{ij}
\end{eqnarray*}
 for $S\in\mathcal{S}$ , which provides a slight speedup since the
column sums are now independent of the iteration. The matrix formulation becomes 
\[
x^{(t+1)}=x^{(t)}+\alpha C^{-1}A_{S}^{T}R^{-1}\left(p_{S}-A_{S}x^{(t)}\right)
\]
 where the diagonal matrices are both independent of the projection
view $S$ as in SIRT. BSSART has been shown \cite{censor2002block} to converge to the minimum
norm solution  
\[
x^{\star}=\argmin_{x}\Vert x\Vert_{2}^{2}\st Ax=p
\]
 as ART for $0<\alpha<2$.

\subsubsection{BICAV}

\cite{censor2001bicav,censor2002block} is another closely-related
algorithm to SART. It updates the volume after each projection view
according to  
\begin{eqnarray*}
x_{j}^{(t+1)} & = & x_{j}^{(t)}+\alpha\frac{1}{c_{j}^{S}}\sum_{i\in S}\frac{p_{i}-\sum_{k=1}^{n}a_{ik}x_{k}^{(t)}}{\sum_{k=1}^{n}a_{ik}^{2}}a_{ij}
\end{eqnarray*}
for $S\in\mathcal{S}$ where $c_{j}^{S}=\sum_{i\in S}\{a_{ij}\ne0\}$
and $\{a_{ij}\ne0\}=1$ when $a_{ij}$ is non-zero is 0 otherwise.
The difference from SART is that it computes the squared norm of the
rows of $A$ and counts the number of non-zero entries in the columns
of $A$. The matrix formulation is 
\[
x^{(t+1)}=x^{(t)}+\alpha C_{S}^{-1}A_{S}^{T}R^{-1}\left(p_{S}-A_{S}x^{(t)}\right)
\]
 where now $r_{i}=\sum_{j}a_{ij}^{2}=\Vert A_{i}\Vert_{2}^{2}$ and
$C_{S}=\diag(c_{j}^{S})$. It is shown \cite{censor2002block} that BICAV converges to the minimum-norm
solution  
\[
\min_{x}\Vert x\Vert^{2}\st Ax=p
\]
 for $0<\alpha<2$.

\subsubsection{OS-SQS}

\cite{erdogan1999ordered,kim2015combining} is closely related to
SART. It is usually derived from a majorization-minimization perspective
\cite{depierro1994modified,erdogan1999ordered,kim2013accelerating,kim2015combining},
but with a specific choice of surrogate functions and parameters \cite{kim2013accelerating}
the update equation becomes 
\begin{eqnarray*}
x_{j}^{(t+1)} & = & x_{j}^{(t)}+\alpha\frac{s}{c_{j}}\sum_{i\in S}\left(p_{i}-\sum_{k=1}^{n}a_{ik}x_{k}^{(t)}\right)a_{ij}
\end{eqnarray*}
for $S\in\mathcal{S}$ where $s$ is the number of subsets $S$ in
$\mathcal{S}$ (number of inner iterations), $c_{j}=\left(\sum_{k=1}^{m}a_{kj}\sum_{i=1}^{n}a_{ki}\right)$,
and in general the set $S$ can contain more than one projection view.
In matrix form it becomes  
\[
x^{(t+1)}=x^{(t)}+\alpha sC^{-1}A_{S}^{T}\left(p_{S}-A_{S}x^{(t)}\right)
\]
 where the matrix $C=\diag(A^{T}A\mathbf{1}_{m})=\diag(c_{j})$ where
$\mathbf{1}_{m}$ is the vector of $m$ ones. OS-SQS is a special
case of the SQS method, which processes all the rows of $A$ at once
like SIRT. Simultaneous SQS, i.e. without ordered subsets, has been
shown \cite{erdogan1999ordered} to converge to a least square solution
\begin{equation}
x^{\star}=\argmin_{x}\Vert Ax-p\Vert_{2}^{2},\label{eq:least-square-problem}
\end{equation}
 and a special case of \emph{relaxed} OS-SQS converges, where the
relaxation parameter becomes iteration-dependent and decreases over
time \cite{ahn2003globally}. However, OS-SQS with fixed $\alpha$
is not known to converge. Therefore, like SART, we assume that it
solves the LS problem in Eq. \ref{eq:least-square-problem} above
approximately.

\subsubsection{CGLS}

\cite{bjorck1996numerical,fessler1999conjugate} is a type of Conjugate
Gradient that solves the least squares normal equations directly.
Like SIRT, it updates the constraint once per full sweep over the
projection rays. The update equation in matrix notation is  
\[
x^{(t+1)}=x^{(t)}+\alpha_{t}\Phi(A^{T}(p-Ax^{(t)}))
\]
 where the update step is a function of the backprojection of the
projection error, and the parameter $\alpha_{t}$ depends on the specific
version of CGLS (here we use the Fletcher-Reeves update rule\cite{bjorck1996numerical}).
CGLS is proven to be convergent to the solution of the LS problem
in Eq. \ref{eq:least-square-problem}. 
\[
x^{\star}=\argmin_{x}\Vert Ax-p\Vert_{2}^{2}.
\]

Note that the function $\Phi(\cdot)$ is more complicated than other
iterative algorithms, and involves several steps with a couple of
auxiliary variables \cite{bjorck1996numerical}.

\section{Tomography Proximal Operators\label{sec:Proximal-Operators}}

Proximal algorithms are a class of optimization algorithms that are
quite flexible and powerful \cite{combettes2011proximal,boyd2011distributed,parikh2013proximal}.
They are generally used to efficiently solve non-smooth, constrained,
distributed, or large scale optimization problems. They are more modular
than other optimization problems, in the sense that they provide a
few lines of code that depend on solving smaller conventional, and
usually simpler, optimization problems called \emph{proximal operator}.
The proximal operator \cite{bauschke2011convex,combettes2011proximal,parikh2013proximal}
for a function $h(\cdot)$ is a generalization of projections on convex
sets, and can be thought of intuitively as getting closer to the optimal
solution while staying close to the current estimate. Formally it
is defined as 
\begin{equation}
\prox_{\lambda h}(u)=\argmin_{x}h(x)+\frac{1}{2\lambda}\Vert x-u\Vert_{2}^{2},\label{eq:prox-operator}
\end{equation}
where $x,u\in\mathbb{R}^{n}$ and $\lambda$ is a regularization parameter.
Many proximal operators of common functions are easy to compute, and
often admit a closed form solution. Computing the proximal operator
of a certain function opens the way to solving hard optimization problems
involving this function and other regularization terms e.g. smoothing
norms or sparsity inducing norms, which otherwise is not generally
easy. We will derive tomography proximal operators for SART, ART,
BICAV, and OS-SQS, where the objective is to solve 
\begin{equation}
\mbox{prox}_{\lambda h}(u)=\argmin_{x}\Vert Ax-p\Vert_{2}^{2}+\frac{1}{2\lambda}\Vert x-u\Vert_{2}^{2}.\label{eq:tomography-prox-operator}
\end{equation}

\begin{table*}
\centering%
\begin{tabular}{ccc}
\toprule 
Method & Update Step & Converges\tabularnewline
\midrule
\midrule 
ART \cite{gordon1970algebraic} & $\begin{aligned}y_{i}^{(t+1)} & =y_{i}^{(t)}+\alpha\frac{\sqrt{2\lambda}p_{i}-\sqrt{2\lambda}\sum_{k}a_{ik}x_{k}^{(t)}-y_{i}^{(t)}}{2\lambda\sum_{k}a_{ik}^{2}+1}\mbox{ for }i\in S\\
x_{j}^{(t+1)} & =x_{j}^{(t)}+\alpha\frac{\sqrt{2\lambda}p_{i}-\sqrt{2\lambda}\sum_{k}a_{ik}x_{k}^{(t)}-y_{i}^{(t)}}{2\lambda\sum_{k}a_{ik}^{2}+1}\sqrt{2\lambda}a_{ij}\mbox{ for }j=1\ldots n
\end{aligned}
$ & Yes\tabularnewline
\midrule 
SART \cite{andersen1984simultaneous} & $\begin{aligned}y_{i}^{(t+1)} & =y_{i}^{(t)}+\alpha\frac{\sqrt{2\lambda}p_{i}-\sqrt{2\lambda}\sum_{k}a_{ik}x_{k}^{(t)}-y_{i}^{(t)}}{\sqrt{2\lambda}\sum_{k}a_{ik}+1}\mbox{ for }i\in S\\
x_{j}^{(t+1)} & =x_{j}^{(t)}+\alpha\frac{\sum_{i\in S}\frac{\sqrt{2\lambda}p_{i}-\sqrt{2\lambda}\sum_{k}a_{ik}x_{k}^{(t)}-y_{i}^{(t)}}{\sqrt{2\lambda}\sum_{k}a_{ik}+1}\sqrt{2\lambda}a_{ij}}{\sqrt{2\lambda}\sum_{i\in S}a_{ij}}\mbox{ for }j=1\ldots n
\end{aligned}
$ & No\tabularnewline
\midrule 
BICAV \cite{censor2001bicav} & $\begin{aligned}y_{i}^{(t+1)} & =y_{i}^{(t)}+\alpha\frac{\sqrt{2\lambda}p_{j}-\sqrt{2\lambda}\sum_{k}a_{jk}x_{k}^{(t)}-y_{j}^{(t)}}{2\lambda\sum_{k}a_{jk}^{2}+1}\mbox{ for }i\in S\\
x_{j}^{(t+1)} & =x_{j}^{(t)}+\alpha\frac{\sum_{i\in S}\frac{\sqrt{2\lambda}p_{i}-\sqrt{2\lambda}\sum_{k}a_{ik}x_{k}^{(t)}-y_{i}^{(t)}}{2\lambda\sum_{k}a_{ik}^{2}+1}\sqrt{2\lambda}a_{ij}}{\sum_{i\in S}\{a_{ij}\ne0\}}\mbox{ for }j=1\ldots n
\end{aligned}
$ & Yes\tabularnewline
\midrule 
OS-SQS \cite{erdogan1999ordered} & $\begin{aligned}x_{j}^{(t+1)} & =x_{j}^{(t)}+\alpha\frac{s}{2\lambda c_{j}+1}\left(2\lambda\sum_{i\in S}(p_{j}-\sum_{k}a_{ik}x_{k}^{(t)})a_{ij}+u_{j}-x_{j}^{(t)}\right)\mbox{ for }j=1\ldots n\end{aligned}
$ & No\tabularnewline
\bottomrule
\end{tabular}

\protect\caption{Summary of the proximal operators update steps.\label{tab:Prox-Operators}
See Sec. \ref{sec:Proximal-Operators} for details.}
\end{table*}

\subsection{SART, ART, and BICAV}

They (approximately) solve the least-norm problem  
\[
x^{\star}=\argmin_{x}\Vert x\Vert_{2}^{2}\st Ax=p.
\]
 What we want is a solver for Eq. \ref{eq:tomography-prox-operator}.
This is equivalent to solving 
\begin{equation}
\min_{x}2\lambda\Vert Ax-p\Vert_{2}^{2}+\Vert x-u\Vert_{2}^{2}.\label{eq:SART-prox-2}
\end{equation}
Introduce new variables $y=\sqrt{2\lambda}(p-Ax)$ and $z=x-u$. The problem becomes 
\begin{eqnarray}
\min_{y,z} & \Vert y\Vert_{2}^{2}+\Vert z\Vert_{2}^{2}\nonumber \\
\st & y+\sqrt{2\lambda}Az=\sqrt{2\lambda}(p-Au).\label{eq:SART-prox-3}
\end{eqnarray}
Rewriting Eq. \ref{eq:SART-prox-3} we arrive at 
\begin{eqnarray*}
\mbox{\ensuremath{\min}}_{y,z} & \left\Vert \left[\begin{array}{c}
y\\
z
\end{array}\right]\right\Vert _{2}^{2}\\
\mbox{subject to} & \left[\begin{array}{cc}
I & \sqrt{2\lambda}A\end{array}\right]\left[\begin{array}{c}
y\\
z
\end{array}\right]=\sqrt{2\lambda}\left(p-Au\right)
\end{eqnarray*}
which can be written as 
\begin{eqnarray*}
\mbox{\ensuremath{\min}}_{\tilde{x}} & \left\Vert \tilde{x}\right\Vert _{2}^{2}\\
\st & \tilde{A}\tilde{x}=\tilde{p},
\end{eqnarray*}
where $\tilde{x}\in\mathbb{R}^{m+n}$, $\tilde{A}\in\mathbb{R}^{m\times m+n}$,
and $\tilde{p}\in\mathbb{R}^{m}$. This is now a consistent under-determined
linear system, and can be solved using either ART, SART, or BICAV. 

Although we introduced new variables $y$ and $z$ and increased the
dimensionality of the problem from $n$ to $n+m$, we can solve the
modified algorithm efficiently with very little computational overhead.
Instead of solving explicitly for the optimal $y^{\star}$ and $z^{\star}$,
we can manipulate the algorithm to solve directly for the optimal
$x^{\star}.$ For example, for SART, the initialization and update
equation for $\tilde{x}$ become 
\begin{eqnarray}
\tilde{x}_{j}^{(0)} & = & \mathbf{0},\nonumber \\
\tilde{x}_{j}^{(t+1)} & = & \tilde{x}_{j}^{(t)}+\alpha\frac{\sum_{i\in\mathcal{S}}\frac{\tilde{p}_{i}-\sum_{k}\tilde{a}_{ik}\tilde{x}_{k}^{(t)}}{\sum_{k}\tilde{a}_{ik}}\tilde{a}_{ij}}{\sum_{i\in S}\tilde{a}_{ij}},\label{eq:SART-proximal-update-1}
\end{eqnarray}
which can be expanded in terms of $y$, $z$, and $A$ as 
\begin{eqnarray*}
y^{(0)} & = & \mathbf{0}_{m}\\
z^{(0)} & = & \mathbf{0}_{n}\\
y_{j}^{(t+1)} & = & y_{j}^{(t)}+\frac{\alpha\sum_{i\in\mathcal{S}}\frac{\tilde{p}_{i}-\sqrt{2\lambda}\sum_{k}a_{ik}z_{k}^{(t)}-y_{i}^{(t)}}{\sqrt{2\lambda}\sum_{k}a_{ik}+1}\delta_{ij}}{1},\\
z_{j}^{(t+1)} & = & z_{j}^{(t)}+\alpha\frac{\sum_{i\in\mathcal{S}}\frac{\tilde{p}_{i}-\sqrt{2\lambda}\sum_{k}a_{ik}z_{k}^{(t)}-y_{i}^{(t)}}{\sqrt{2\lambda}\sum_{k}a_{ik}+1}\sqrt{2\lambda}a_{ij}}{\sqrt{2\lambda}\sum_{i\in S}a_{ij}},
\end{eqnarray*}
where $\delta_{ij}=1$ when $i=j$ and 0 otherwise. Using the fact
that $z=x-u$ and $\tilde{p}_{i}=\sqrt{2\lambda}p_{i}-\sqrt{2\lambda}\sum_{k}a_{ik}u_{k}$
and simplifying we arrive at 
\begin{eqnarray*}
y^{(0)} & = & \mathbf{0}_{m}\\
x^{(0)} & = & u\\
y_{j}^{(t+1)} & = & y_{j}^{(t)}+\alpha\sum_{i\in S}\frac{\sqrt{2\lambda}p_{i}-\sqrt{2\lambda}\sum_{k}a_{ik}x_{k}^{(t)}-y_{i}^{(t)}}{\sqrt{2\lambda}\sum_{k}a_{ik}+1}\delta_{ij},\\
x_{j}^{(t+1)} & = & x_{j}^{(t)}+\alpha\frac{\sum_{i\in S}\frac{\sqrt{2\lambda}p_{i}-\sqrt{2\lambda}\sum_{k}a_{ik}x_{k}^{(t)}-y_{i}^{(t)}}{\sqrt{2\lambda}\sum_{k}a_{ik}+1}\sqrt{2\lambda}a_{ij}}{\sqrt{2\lambda}\sum_{i\in S}a_{ij}}.
\end{eqnarray*}

Following the same line of reasoning, we can arrive at similar update
formulas for both ART and BICAV. The steps are summarized in Table
\ref{tab:Prox-Operators}. Alg. \ref{alg:ART-like-Prox-Operator} provides an outline of the
proximal operator.
\begin{algorithm}[h]
\protect\protect\caption{\label{alg:ART-like-Prox-Operator}SART, ART, and BICAV Proximal Operator}

\begin{algorithmic}[1]

\Require $A\in\mathbb{R}{}^{m\times n}$, $\alpha,\lambda\in\mathbb{R}$,
$p\in\mathbb{R}{}^{m}$, $u\in\mathbb{R}^{n}$

\State Initialize: 
\begin{align*}
x^{(0)} & =u\\
y^{(0)} & =\mathbf{0}_{m}
\end{align*}

\ForAll {$t=1\ldots T$ }

\ForAll {subsets $S\in\mathcal{S}$ } 

\State Update according to Table \ref{tab:Prox-Operators}:
\begin{align*}
x_{j}^{(t+1)} & =x_{j}^{(t)}+\alpha\Delta x_{j}^{(t)}\mbox{ for }j=1\ldots n\\
y_{i}^{(t+1)} & =y_{i}^{(t)}+\alpha\Delta y_{i}^{(t)}\mbox{ for }i\in S
\end{align*}

\State$x^{(t+1)}=\mbox{clip}(x^{(t+1)})$

\EndFor

\EndFor

\Return $x^{\star}=\argmin_{x}\Vert Ax-p\Vert_{2}^{2}+\frac{1}{2\lambda}\Vert x-u\Vert_{2}^{2}$

\end{algorithmic} 
\end{algorithm}

\subsection{OS-SQS}

We want to express the proximal operator problem  
\[
\prox_{\lambda h}(u)=\argmin_{x}h(x)+\frac{1}{2\lambda}\Vert x-u\Vert_{2}^{2}
\]
 in the form of the LS problem that can be solved (approximately)
by OS-SQS i.e.  
\[
x^{\star}=\argmin_{x}\Vert Ax-p\Vert_{2}^{2}.
\]
 Rewrite as  
\[
\argmin_{x}2\lambda\Vert Ax-p\Vert_{2}^{2}+\Vert x-u\Vert_{2}^{2},
\]
which is equivalent to 
\begin{align*}
 & \argmin_{x}\left\Vert \left[\begin{array}{c}
\sqrt{2\lambda}A\\
I
\end{array}\right]x-\left[\begin{array}{c}
\sqrt{2\lambda}p\\
u
\end{array}\right]\right\Vert _{2}^{2}\\
 & \iff\argmin_{x}\left\Vert \tilde{A}x-\tilde{p}\right\Vert _{2}^{2}
\end{align*}
 where 
\begin{eqnarray*}
\tilde{A} & = & \left[\begin{array}{c}
\sqrt{2\lambda}A\\
I
\end{array}\right]\in\mathbb{R}^{m+n\times n}\\
\tilde{p} & = & \left[\begin{array}{c}
\sqrt{2\lambda}p\\
u
\end{array}\right]\in\mathbb{R}^{m+n}.
\end{eqnarray*}
 The weighting matrix $\tilde{C}\in\mathbb{R}^{n\times n}$ now becomes
\begin{eqnarray*}
\tilde{C} & = & \mbox{diag}(\tilde{A}^{T}\tilde{A}\mathbf{1})\\
 & = & \mbox{diag}\left((2\lambda A^{T}A+I)\mathbf{1}\right)\\
 & = & \mbox{diag}\left(2\lambda A^{T}A\mathbf{1}+\mathbf{1}\right)
\end{eqnarray*}
 and its diagonal entries are 
\[
\tilde{c}_{j}=2\lambda c_{j}+1.
\]

Write the matrix update equation in terms of $\tilde{A}$ and $\tilde{p}$
as  
\begin{align*}
x^{(t+1)} & =x^{(t)}+\alpha s\tilde{C}^{-1}\tilde{A}_{S}^{T}\left(\tilde{p}_{S}-\tilde{A}_{S}x^{(t)}\right)\\
 & =x^{(t)}+\alpha s\tilde{C}^{-1}\left[\begin{smallmatrix}\sqrt{2\lambda}A_{S}^{T} & I\end{smallmatrix}\right]\left[\begin{smallmatrix}\sqrt{2\lambda}(p_{S}-A_{S}x^{(t)})\\
u-x^{(t)}
\end{smallmatrix}\right]\\
 & =x^{(t)}+\alpha s\tilde{C}^{-1}{\scriptstyle \left(2\lambda A_{S}^{T}(p_{S}-A_{S}x^{(t)})+(u-x^{(t)})\right)}.
\end{align*}
 In component form it becomes 
\begin{align*}
x_{j}^{(t+1)} & =x_{j}^{(t)}+\alpha\frac{s}{\tilde{c}_{j}}{\scriptstyle \left(2\lambda\sum_{i\in S}(p_{j}-\sum_{k}a_{ik}x_{k}^{(t)})a_{ij}+u_{j}-x_{j}^{(t)}\right)}.
\end{align*}

The steps are summarized in Table \ref{tab:Prox-Operators}. Alg.
\ref{alg:OS-SQS-Prox-Operator} gives an outline.
\begin{algorithm}[h]
\protect\protect\caption{\label{alg:OS-SQS-Prox-Operator}OS-SQS Proximal Operator}

\begin{algorithmic}[1]

\Require $A\in\mathbb{R}{}^{m\times n}$, $\alpha,\lambda\in\mathbb{R}$,
$p\in\mathbb{R}{}^{m}$, $u\in\mathbb{R}^{n}$

\State Initialize: 
\begin{align*}
x^{(0)} & =\mathbf{0}_{n}
\end{align*}

\ForAll {$t=1\ldots T$ }

\ForAll {subsets $S\in\mathcal{S}$ } 

\State $x_{j}^{(t+1)}=x_{j}^{(t)}+\alpha\Delta x_{j}^{(t)}\mbox{ for }j=1\ldots n$
according to Table \ref{tab:Prox-Operators}.

\State$x^{(t+1)}=\mbox{clip}(x^{(t+1)})$

\EndFor

\EndFor

\Return $x^{\star}=\argmin_{x}\Vert Ax-p\Vert_{2}^{2}+\frac{1}{2\lambda}\Vert x-u\Vert_{2}^{2}$

\end{algorithmic} 
\end{algorithm}

\section{TRex Proximal Framework\label{sec:Proximal-Framework}}

\subsection{Proximal Algorithm\label{sub:Proximal-Algorithm}}

The overall problem we want to solve is a regularized data fitting
problem, namely 
\begin{equation}
\argmin_{x}f(x)+g(Kx),\label{eq:full-problem}
\end{equation}
where $f(\cdot)$ is a data fitting term that measures how much the
solution fits the data and that depends on the measurement noise model
assumed, $K\in\mathbb{R}^{d\times n}$ is a matrix, and $g(\cdot)$
is a regularization term that imposes constraints on acceptable solutions.
We will use the Linearized ADMM method \cite{boyd2011distributed,parikh2013proximal}
(also known as Inexact Split Uzawa \cite{esser2010general,zhang2011unified}
or Proximal ADMM \cite{eckstein1994some,fazel2013hankel,chen2015inertial,yuan2015l0tv}), for solving this problem for different data terms and different
regularizers. 

It rewrites Eq. \ref{eq:full-problem} into the equivalent form 
\begin{align*}
\argmin_{x,z} & f(x)+g(z)\\
\st & Kx=z,
\end{align*}
 writes out the scaled augmented Lagrangian function \cite{parikh2013proximal}
\[
\mathcal{L}_{\rho}(x,z,y)=f(x)+g(z)+\frac{\rho}{2}\Vert Kx-z+y\Vert_{2}^{2},
\]
 and then applies alternating minimization for the variables $x,$
$z$, and $y$ in turn: 
\begin{align*}
x^{(t+1)} & =\argmin_{x}f(x)+\frac{\rho}{2}\Vert Kx-z^{(t)}+y^{(t)}\Vert^{2}\\
z^{(t+1)} & =\argmin_{z}g(z)+\frac{\rho}{2}\Vert Kx^{(t+1)}-z+y^{(t)}\Vert^{2}\\
y^{(t+1)} & =y^{(t)}+Kx^{(t+1)}-z^{(t+1)}.
\end{align*}

The problem with the $x$ step is that it contains the quadratic term
$\Vert Kx\Vert^{2}=x^{T}K^{T}Kx$ in the minimization makes it hard
to minimize since it's not straightforward. We can cancel out that
term by adding the following \emph{proximal }term that makes it strongly
convex and keeps the solution \emph{close} to the previous iteration
\[
\frac{1}{2}\Vert x-x^{k}\Vert_{S}^{2}=\frac{1}{2}(x-x^{k})^{T}S(x-x^{k})
\]
 to the objective fundtion where the special matrix $S$ is 
\[
S=\frac{1}{\mu}I-\rho K^{T}K
\]
 and this gives the modified $x$ step 
\begin{align*}
x^{(t+1)} & = &  & \argmin_{x}f(x)+\frac{\rho}{2}\Vert Kx-z^{(t)}+y^{(t)}\Vert^{2}+\\
 &  &  & \frac{1}{2}\Vert x-x^{(t)}\Vert_{S}^{2}\\
 & = &  & \argmin_{x}f(x)+\frac{\rho}{2}\Vert Kx\Vert^{2}-\rho\langle Kx,z^{(t)}-y^{(t)}\rangle+\\
 &  &  & \frac{1}{2}(x-x^{(t)})^{T}S(x-x^{(t)})\\
 & = &  & \argmin_{x}f(x)+\frac{\rho}{2}\Vert Kx\Vert^{2}-\rho\langle x,K^{T}(z^{(t)}-y^{(t)})\rangle+\\
 &  &  & \frac{1}{2}\Vert x\Vert_{S}^{2}-\langle x,Sx^{(t)}\rangle\\
 & = &  & \argmin_{x}f(x)+\\
 &  &  & \frac{1}{2\mu}\Vert x-\mu\rho K^{T}(z^{(t)}-y^{(t)})-\mu Sx^{(t)}\Vert^{2}\\
 & = &  & \argmin_{x}f(x)+\\
 &  &  & \frac{1}{2\mu}\Vert x-x^{(t)}-\mu\rho K^{T}\left(z^{(t)}-y^{(t)}-Kx^{(t)}\right)\Vert^{2}
\end{align*}
 which is simply the proximal operator of $f(x)$ with input $x^{(t)}+\mu\rho\left(z^{(t)}-u^{(t)}-K^{T}Kx^{(t)}\right)$
i.e. the iterations now become

\begin{align*}
x^{(t+1)} & =\prox_{\mu f}\left(x^{(t)}+\mu\rho K^{T}\left(z^{(t)}-y^{(t)}-Kx^{(t)}\right)\right)\\
z^{(t+1)} & =\prox_{\rho^{-1}g}\left(Kx^{(t+1)}+y^{(t)}\right)\\
y^{(t+1)} & =y^{(t)}+Kx^{(t+1)}-z^{(t+1)}.
\end{align*}

The algorithm is convergent for any $\rho>0$ and $\mu>\nicefrac{1}{\rho\Vert K\Vert^{2}}$
\cite{parikh2013proximal,nien2015fast}. The steps are summarized
in Alg. \ref{alg:LADMM-algorithm}. This framework is very flexible,
and we will show how to solve for different data terms and different
regularizers. 

\begin{algorithm}[h]
\protect\protect\caption{\label{alg:LADMM-algorithm}Linearized ADMM}

\begin{algorithmic}[1]

\Require $K\in\mathbb{R}{}^{d\times n}$, $\rho,\mu\in\mathbb{R}$
such that $\mu\rho\Vert K\Vert^{2}<1$, initial values $x^{(0)}\in\mathbb{R}^{n}$
and $z^{(0)}\in\mathbb{R}^{d}$

\State Initialize $y^{(0)}=\mathbf{0}_{d}$

\ForAll {$t=1\ldots T$ }

\State $x^{(t+1)}=\mbox{prox}_{\mu f}\left(x^{(t)}-\rho\mu K^{T}(Kx^{(t)}-z^{(t)}+y^{(t)})\right)$

\State $z^{(t+1)}=\mbox{prox}_{\rho^{-1}g}\left(Kx^{(t+1)}+y^{(t)}\right)$

\State $y^{(t+1)}=y^{(t)}+Kx^{(t+1)}-z^{(t+1)}$

\EndFor

\Return $x^{(T)}=\argmin_{x}f(x)+g(Kx)$

\end{algorithmic} 
\end{algorithm}

\subsection{Data Terms\label{sub:Data-Terms}}

We will consider the following data fidelity terms, which correspond
to specific noise models:

\subsubsection{Gaussian Noise}

Assume the measurements $p_{i}\forall i=1,\ldots m$ follow the model
\begin{equation}
p_{i}=a_{i}^{T}x+\varepsilon_{i}\label{eq:Gaussian-Noise-Model}
\end{equation}
where the noise $\varepsilon\sim\mathbb{N}(0,\sigma^{2})$ follows
a Gaussian distribution. Maximizing the projection data log-likelihood 
\begin{equation}
\mathcal{L}_{\text{G}}(p)\propto-\sum_{i}\left(p_{i}-a_{i}^{T}x\right)^{2}\label{eq:Gaussian-Log-Likelihood}
\end{equation}
is equivalent to minimizing the LS $\ell_{2}$ norm data term
\begin{equation}
f_{\text{G}}(x)=\Vert Ax-p\Vert_{2}^{2}=\sum_{i=1}^{m}(A_{i}^{T}x-p_{i})^{2}.\label{eq:gaussian-data-term}
\end{equation}
We can solve proximal operator $\prox_{\lambda f_{\text{G}}}(\cdot)$
directly using any of the algorithms from Table \ref{tab:Prox-Operators}.

\subsubsection{Poisson Noise}

It can be shown that assuming an \emph{approximated} Poisson noise
model leads to a WLS data term, where the weights are proportional
to the detector measurements \cite{clinthorne1993preconditioning,depierro1994modified,elbakri2002statistical,thibault2007three}.
Indeed, the actual measurements produced by the X-ray CT scanner represent
X-ray photon energy reaching the detector as compared to the energy
leaving the X-ray gun. These are related to each other and to the
linear attenuation coefficient according to Beer-Lambert law \cite{hsieh2009computed}:
\[
I_{t}=I_{o}e^{-\int\mu(l)dl}
\]
where $I_{t}$ is the transmitted intensity as measured by the detector,
$I_{o}$ is the emitted intensity from the source, $\mu(l)$ is the
linear attenuation coefficient of the material as a function of length
$l$. The exponent represents the line integrals (projection data)
we are dealing with. In particular, assuming that the X-ray photons
are monochromatic (have only one single energy) i.e. ignoring \emph{beam
hardening}, the projection line integral data at detector $i$ is
obtained from the physical measurements as 
\begin{equation}
p_{i}=-\ln\frac{I_{t}^{i}}{I_{o}^{i}}\label{eq:p_i}
\end{equation}

where $I_{t}^{i}$ is the intensity measured by detector $i$ and
$I_{o}^{i}$ is the emitted intensity. The detector measurements are
stochastic in nature, and assuming a Poisson distribution with mean
$I_{o}^{i}\exp(-p_{i})$ we get  
\[
I_{t}^{i}\sim\mathbb{P}(I_{o}^{i}e^{-p_{i}})\approx\mathbb{P}(I_{o}^{i}e^{-A_{i}^{T}x}).
\]

Using the ML approach, we maximize the log-likelihood of the measured
data: 
\begin{equation}
\mathcal{L}_{P}(x)=\sum_{i}I_{t}^{i}\ln\left(I_{o}^{i}e^{-A_{i}^{T}x}\right)-I_{o}^{i}e^{-A_{i}^{T}x}=\sum_{i}\phi_{i}\left(A_{i}^{T}x\right)\label{eq:Poisson-log-likelihood-full}
\end{equation}
 where 
\[
\phi_{i}(q)=I_{t}^{i}\ln\left(I_{o}^{i}e^{-q}\right)-I_{o}^{i}e^{-q}.
\]

Applying a second-order Taylor's expansion for $\phi_{i}(q)$ around
an estimate of the $i$th line integral $p_{i}$ from Equation \ref{eq:p_i}
\cite{elbakri2002statistical}: 
\begin{align*}
\phi_{i}(q) & \approx\phi_{i}(p_{i})+\frac{d\phi_{i}}{dq}(p_{i})(q-p_{i})+\frac{1}{2}\frac{d^{2}\phi_{i}}{dq^{2}}(p_{i})(q-p_{i})^{2}\\
 & =(I_{t}^{i}\ln I_{t}^{i}-I_{t}^{i})-\frac{I_{t}^{i}}{2}(q-p_{i})^{2}
\end{align*}
The first term is independent of $q$ and can be dropped (since we
are interested in minimizing $\mathcal{L}_{P}(x)$). Substituting
in Equation \ref{eq:Poisson-log-likelihood-full}, we end up with
the approximated log-likelihood 
\begin{equation}
\mathcal{L}_{\text{G}}(x)\approx-\sum_{i}\frac{I_{t}^{i}}{2}\left(A_{i}^{T}x-p_{i}\right)^{2}=-\sum w_{i}\left(A_{i}^{T}x-p_{i}\right){}^{2}\label{eq:Poisson-log-likelihood-approx}
\end{equation}
 where $w_{i}$ is the weight for projection measurement $i$ and
is proportional to the measurement of the incident X-ray intensity
on detector $i$ i.e. $w_{i}\propto I_{t}^{i}.$ Typically, the weights
$w_{i}$ are normalized to have a maximum of 1, and we could apply
any non-decreasing mapping on $w_{i}$, e.g. the square root, before
feeding into the optimization problem, see Sec. \ref{sub:Poisson-Mapping-Functions-Comparison}. 

Maximizing the likelihood is equivalent to minimizing the WLS data
term 

\begin{equation}
f_{\text{P}}(x)=\Vert Ax-p\Vert_{W}^{2}=\sum_{i=1}^{m}w_{i}(a_{i}^{T}x-p_{i})^{2}\label{eq:Poisson-data-term}
\end{equation}
 where $\Vert x\Vert_{W}=x^{T}Wx$ and $W=\diag(w_{i})\in\mathbb{R}^{m\times m}$
is a diagonal matrix containing weights for each measurement. 

We can solve the proximal operator 
\begin{equation}
\prox_{\lambda f_{\text{P}}}(u)=\min_{x}\Vert Ax-p\Vert_{W}^{2}+\frac{1}{2\lambda}\Vert x-u\Vert^{2}\label{eq:Poisson-data-term-proximal-operator-initial}
\end{equation}

as follows. Define $\tilde{p}\in\mathbb{R}^{m}$ and $\tilde{A}\in\mathbb{R}^{m\times n}$
as 
\begin{eqnarray*}
\tilde{p} & = & W^{\frac{1}{2}}p\\
\tilde{A} & = & W^{\frac{1}{2}}A
\end{eqnarray*}
 where $W^{\frac{1}{2}}=\diag\left(\sqrt{w_{i}}\right).$ We get 
\begin{eqnarray*}
\Vert\tilde{A}x-\tilde{p}\Vert_{2}^{2} & = & (\tilde{A}x-\tilde{p})^{T}(\tilde{A}x-\tilde{p})\\
 & = & \left(W^{\frac{1}{2}}(Ax-p)\right)^{T}\left(W^{\frac{1}{2}}(Ax-p)\right)\\
 & = & (Ax-p)^{T}W^{\frac{1}{2}}W^{\frac{1}{2}}(Ax-p)\\
 & = & \Vert Ax-p\Vert_{W}^{2}.
\end{eqnarray*}
Now this is in the form that can be solved with the algorithms in
Table \ref{tab:Prox-Operators}  
\[
\prox_{\lambda f_{\text{P}}}(u)=\min_{x}\Vert\tilde{A}x-\tilde{p}\Vert^{2}+\frac{1}{2\lambda}\Vert x-u\Vert^{2}
\]
with input matrix $\tilde{A}$ and projections $\tilde{p}$.

\subsection{Regularizers\label{sub:Regularizers}}

The regularizers impose constraints on the reconstruction volume.
We consider the following regularizers:

\subsubsection{Isotropic Total Variation (ITV)}

It is the sum of the gradient magnitude at each voxel \cite{rudin1992nonlinear,sidky2012convex,chambolle2011first}
i.e. 
\begin{equation}
h_{\text{ITV}}(x)=g_{\text{ITV}}(Kx)=\sigma\Vert x\Vert_{\text{TV}}=\sigma\sum_{i}\Vert\nabla x_{i}\Vert_{2}\label{eq:h-ITV}
\end{equation}
 where $\nabla x_{i}=\left[\begin{smallmatrix}\nabla x_{i}^{1} & \nabla x_{i}^{2}\end{smallmatrix}\right]^{T}$
is the discrete gradient at voxel $i$ containing the horizontal forward different $\nabla x_{i}^{1}$ and
the vertical forward difference $\nabla x_{i}^{2}$. It can be represented in the form of Eq. \ref{eq:full-problem} 
\[
h_{\text{ITV}}(x)=g_{\text{ITV}}(Kx)
\]
 by defining the matrix $K=D\in\mathbb{R}^{2n\times n}$ to be the
forward difference matrix that produces the discrete gradient $\nabla x\in\mathbb{R}^{2n}$ 
\[
\nabla x=\begin{bmatrix}\nabla x_{i}\\
\vdots\\
\nabla x_{n}
\end{bmatrix}=Dx
\]
 and defining for $u\in\mathbb{R}^{2n}=\begin{bmatrix}u_{1}^{T} & \cdots & u_{n}^{T}\end{bmatrix}^{T}$
\[
g_{\text{ITV}}(u)=\sigma\sum_{i}\Vert u_{i}\Vert_{2}
\]
 The proximal operator $\prox_{\lambda g_{\text{ITV}}}(u)$ is \cite{esser2010general,chambolle2011first}
\begin{equation}
\prox_{\lambda g_{\text{ITV}}}(u_{i})=u-\frac{\lambda\sigma u_{i}}{\max(\lambda\sigma,\Vert u_{i}\Vert_{2})}\label{eq:ITV-prox}
\end{equation}
 where $u_{i}\in\mathbb{R}^{2}$ is the $i$th component of $u$.
Intuitively it projects back the vector $u_{i}$ to be on the Euclidean
ball of radius $\sigma$.

\subsubsection{Anisotropic Total Variation (ATV)}

It is a simplification of ITV \cite{sidky2008image}, and is defined
as 
\begin{equation}
h_{\text{ATV}}(x)=\sigma\Vert\nabla x\Vert_{1}\label{eq:h-ATV}
\end{equation}
 which is the $\ell_{1}$ norm of the gradient $\nabla x$ of the
volume. It can be written in the form of Eq. \ref{eq:full-problem}
\[
h_{\text{ATV}}(x)=g_{\text{ATV}}(Kx)
\]
 by defining $K=D$ as in the ITV case and defining for $u\in\mathbb{R}^{2n}$
\[
g_{\text{ATV}}(u)=\sigma\Vert u\Vert_{1}=\sigma\sum_{i}\Vert u_{i}\Vert_{1}.
\]
The proximal operator $\prox_{\lambda g_{\text{ATV}}}(u)$ is \cite{esser2010general,chambolle2011first}
\begin{equation}
\prox_{\lambda g_{\text{ITV}}}(u_{i})=\sign(u_{i})\odot\max(0,\vert u_{i}\vert-\sigma)\label{eq:ATV-prox}
\end{equation}
 which is the soft thresholding function \cite{sidky2012convex},
where the max and product are component-wise operations.

\subsubsection{Sum of Absolute Differences (SAD)}

It is an extension to the ATV by adding more forward differences around
each voxel \cite{gregson2012stochastic}. In particular, it sums the
differences of the voxels in the $3\times3$ neighborhood around each
voxel 
\begin{equation}
h_{\text{SAD}}(x)=\sigma\sum_{i}\sum_{k\in\mathcal{N}(i)}\vert x_{i}-x_{k}\vert\label{eq:h-SAD}
\end{equation}
where $\mathcal{\mathcal{N}}(i)$ contains the voxels in the neighborhood
around voxel $i$. It can be written similarly in the form 
\[
h_{\text{SAD}}(x)=g_{\text{SAD}}(Kx)
\]
 by defining $K\in\mathbb{R}^{8n\times n}$ that computes the 8 forward
differences in the $3\times3$ neighborhood and defining for $u\in\mathbb{R}^{8n}$
\[
g_{\text{SAD}}(u)=\sigma\Vert u\Vert_{1}=\sigma\sum_{i}\Vert u_{i}\Vert_{1}.
\]
 The proximal operator $\prox_{\lambda g_{\text{ATV}}}(u)$ is similar
to the ATV case: 
\begin{equation}
\prox_{\lambda g_{\text{ITV}}}(u_{i})=\sign(u_{i})\odot\max(0,\vert u_{i}\vert-\sigma).\label{eq:SAD-prox}
\end{equation}
 The SAD prior has been shown \cite{gregson2012stochastic} to produce
excellent results in stochastic tomography reconstruction.

\section{Experiments\label{sec:Experiments}}

\subsection{Datasets and Implementation Details}

\begin{figure*}
\centering\subfloat[Modified Shepp-Logan]{%
\begin{minipage}[t]{0.33\textwidth}%
\includegraphics[width=1\textwidth]{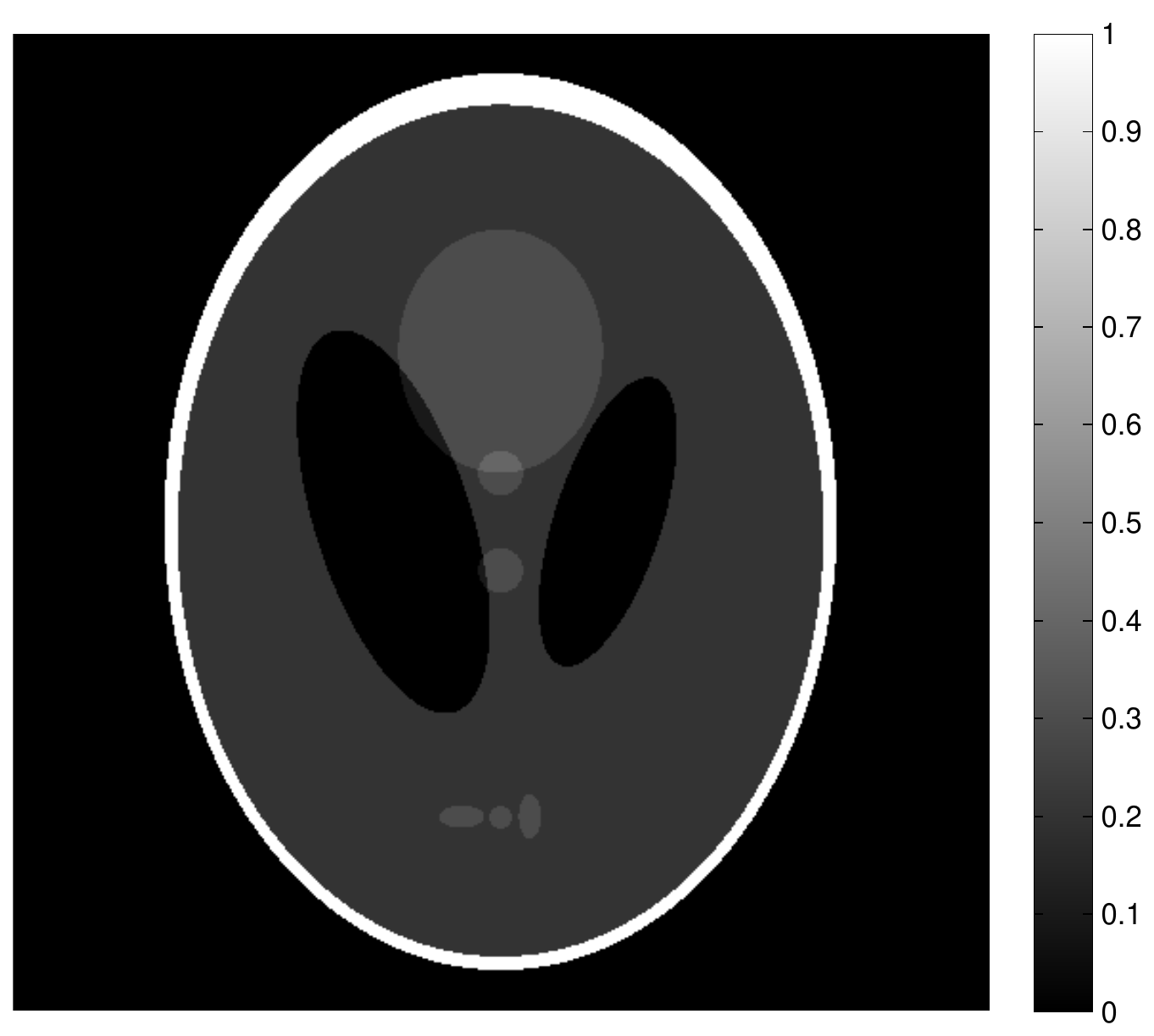}%
\end{minipage}

}\subfloat[NCAT]{%
\begin{minipage}[t]{0.33\textwidth}%
\includegraphics[width=1\textwidth]{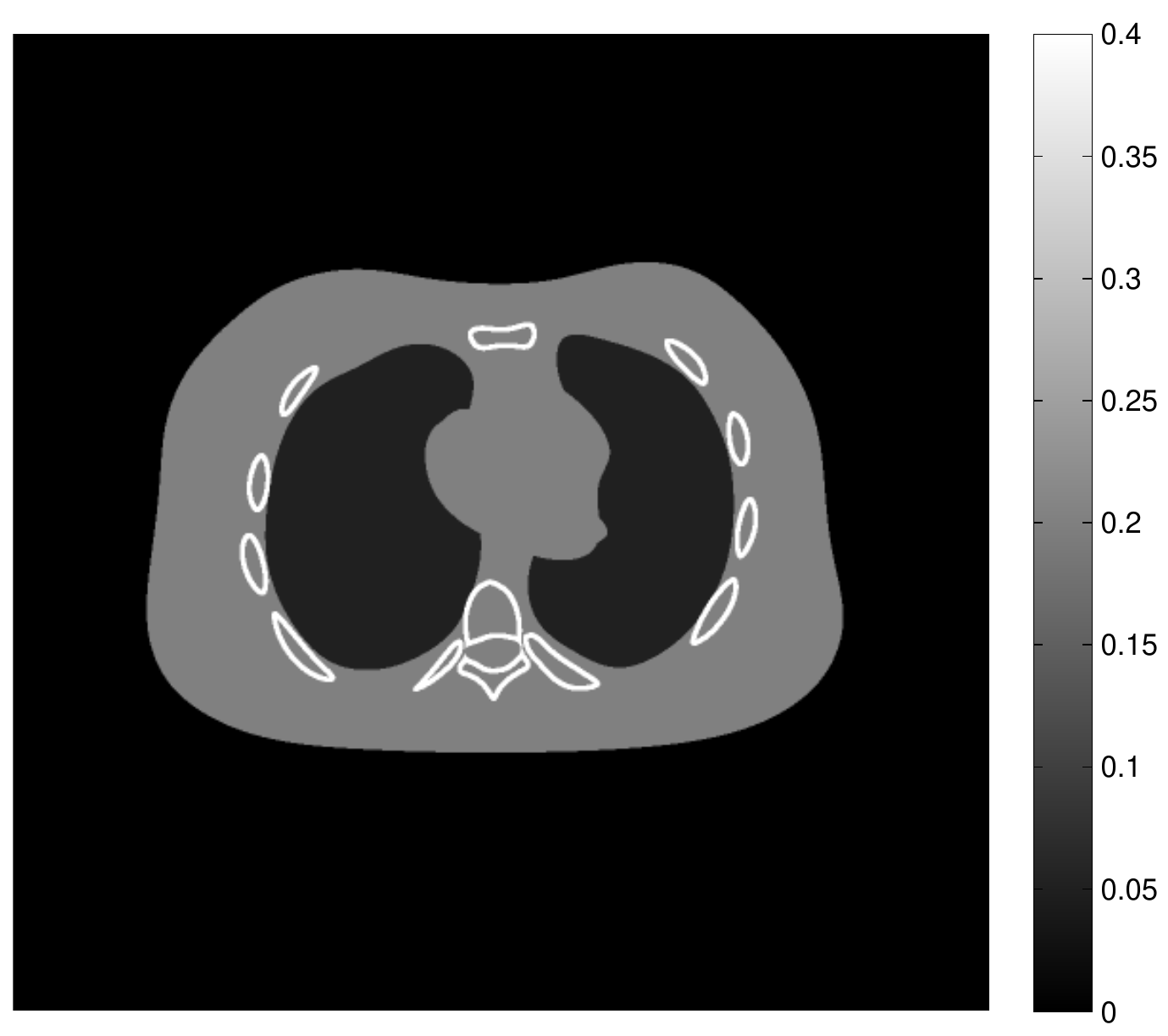}%
\end{minipage}

}\subfloat[Mouse]{%
\begin{minipage}[t]{0.33\textwidth}%
\includegraphics[width=1\textwidth]{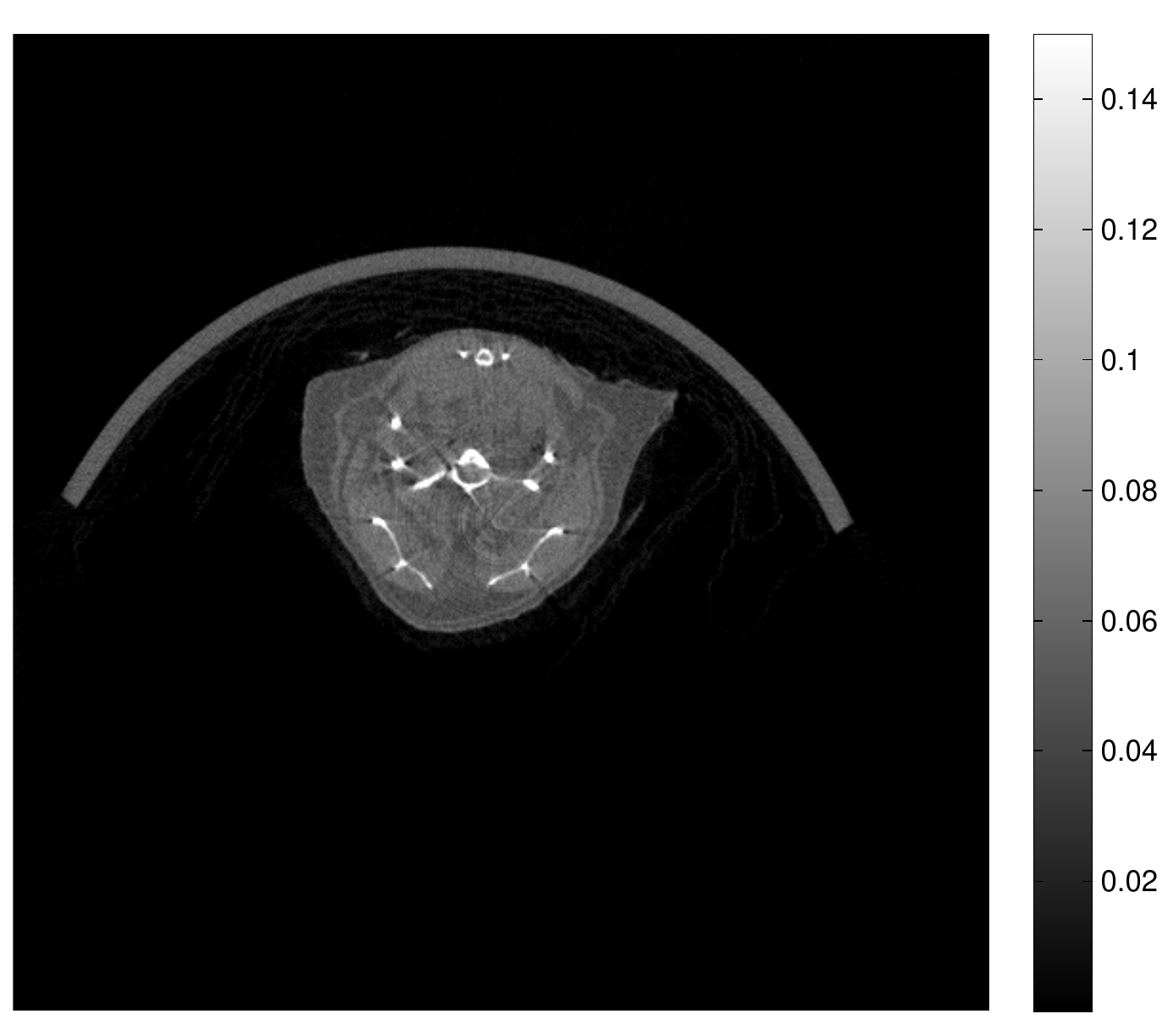}%
\end{minipage}

}

\protect\caption{The datasets used. (c)  shows the \emph{ground truth} converged
result from the projections. \label{fig:Phantoms}}
\end{figure*}

We present experiments on two simulated phantoms and one real dataset,
see Fig. \ref{fig:Phantoms}. The phantoms are: the modified 2D Shepp-Logan head phantom \cite{toft1996radon};
and a 2D slice of the NCAT phantom \cite{segars2002study}. The phantoms
were generated at a resolution of $512\times512$ pixels, and ground
truth sinograms were generated in ASTRA using a fan beam geometry
with $888$ detectors, isotropic pixels of 1 mm, isotropic detectors
of $1.0239$ mm, and source-to-detector distance of $949.075$ mm. We assumed Poisson measurement noise with emitted intensity count
$I_{0}=10^{5}$ to generate the noisy projections used.

The real dataset is a 2D slice of a 3D cone beam scan of a mouse from
the Exxim Cobra software %
\footnote{available from http://www.exxim-cc.com/ %
}. The data contains 194 projections (over 194 degrees) of a fan beam
geometry with 512 detectors of size 0.16176 mm, source-to-detector
distance of 529.29 mm, source-to-isocenter distance of 395.73 mm,
and reconstructed volume of $512\times512$ pixels of isotropic size
0.12 mm. We ran 500 iterations of BSSART with $\alpha=0.1$ to generate
the \emph{ground truth} volume, but we note that results on this dataset
should be taken with a grain of salt. We measure performance in terms
of SNR (signal-to-noise ratio) defined as 
\[
\mbox{SNR}(x,\hat{x})=10\log\frac{\sum_{j}\hat{x}_{j}^{2}}{\sum_{j}\left(x_{j}-\hat{x}_{j}\right)^{2}}
\]
 where $x\in\mathbb{R}^{n}$ is the current estimate of the volume
and $\hat{x}\in\mathbb{R}^{n}$ is the ground truth volume.

We clip the reconstruction estimate $x$ at the end of each inner
iteration (i.e. after each update step) using this function 
\[
\mbox{clip}(x)=\max(0,x)
\]
 to get rid of negative voxel values. 

We implemented all methods using ASTRA with a mix of C++ and Matlab
code. The iterative algorithms not present in ASTRA, namely BICAV, BSSART,
and OS-SQS, were implemented in C++. The proximal operators were also
implemented in C++. The Linearized ADMM was implemented in Matlab.
We also modified existing algorithms in ASTRA to suit our needs e.g.
compute SNR, report run times, etc. All experiments were run on one core of an Intel Xeon E5-280 2.7
GHz with 64 GB RAM.

\subsection{Iterative Algorithms Comparison\label{sub:Iterative-Algorithms-Comparison}}

\begin{figure*}
\center\subfloat[Modified Shepp-Logan]{%
\begin{minipage}[t]{0.33\textwidth}%
\includegraphics[width=1\textwidth]{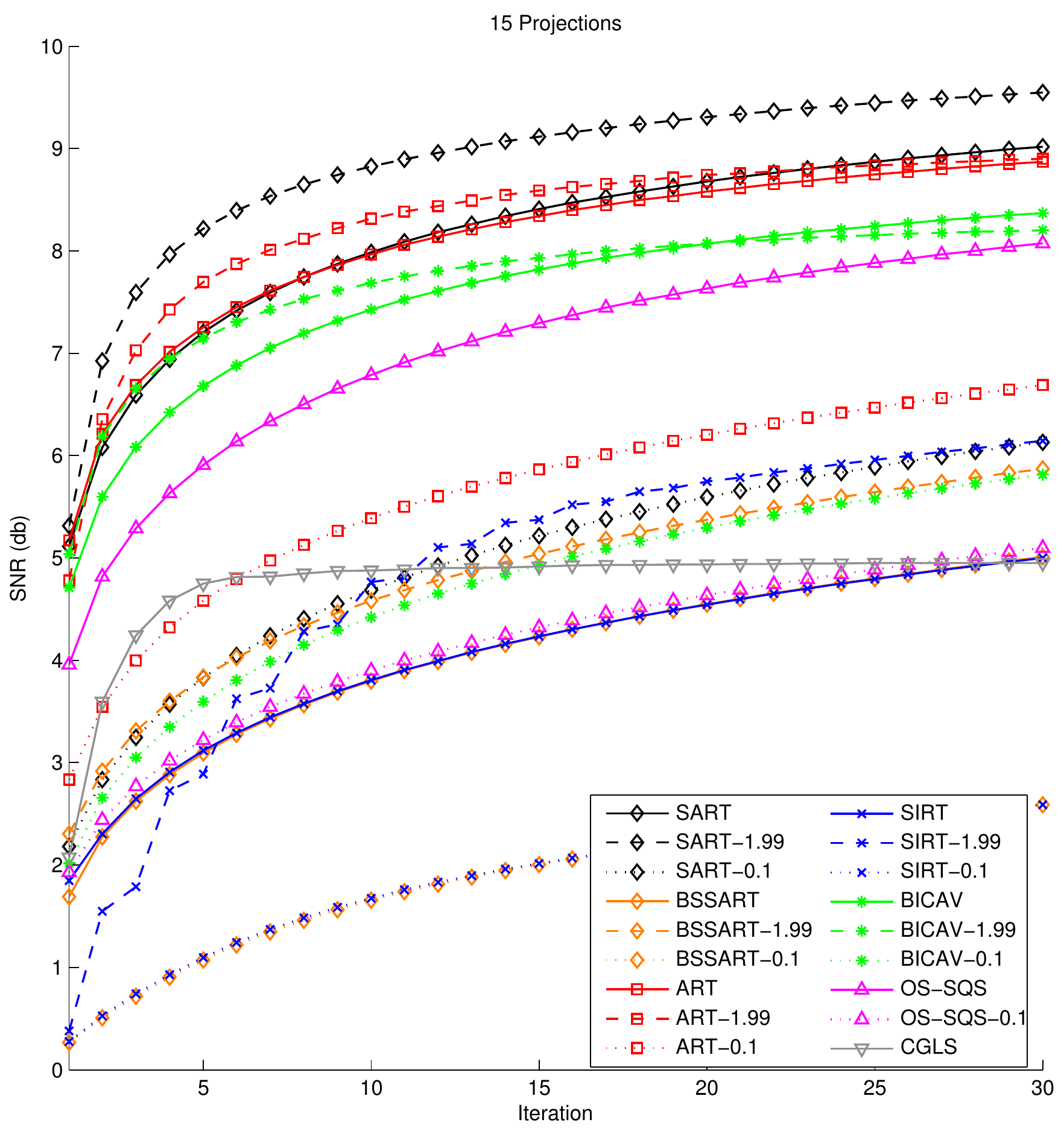}

\includegraphics[width=1\textwidth]{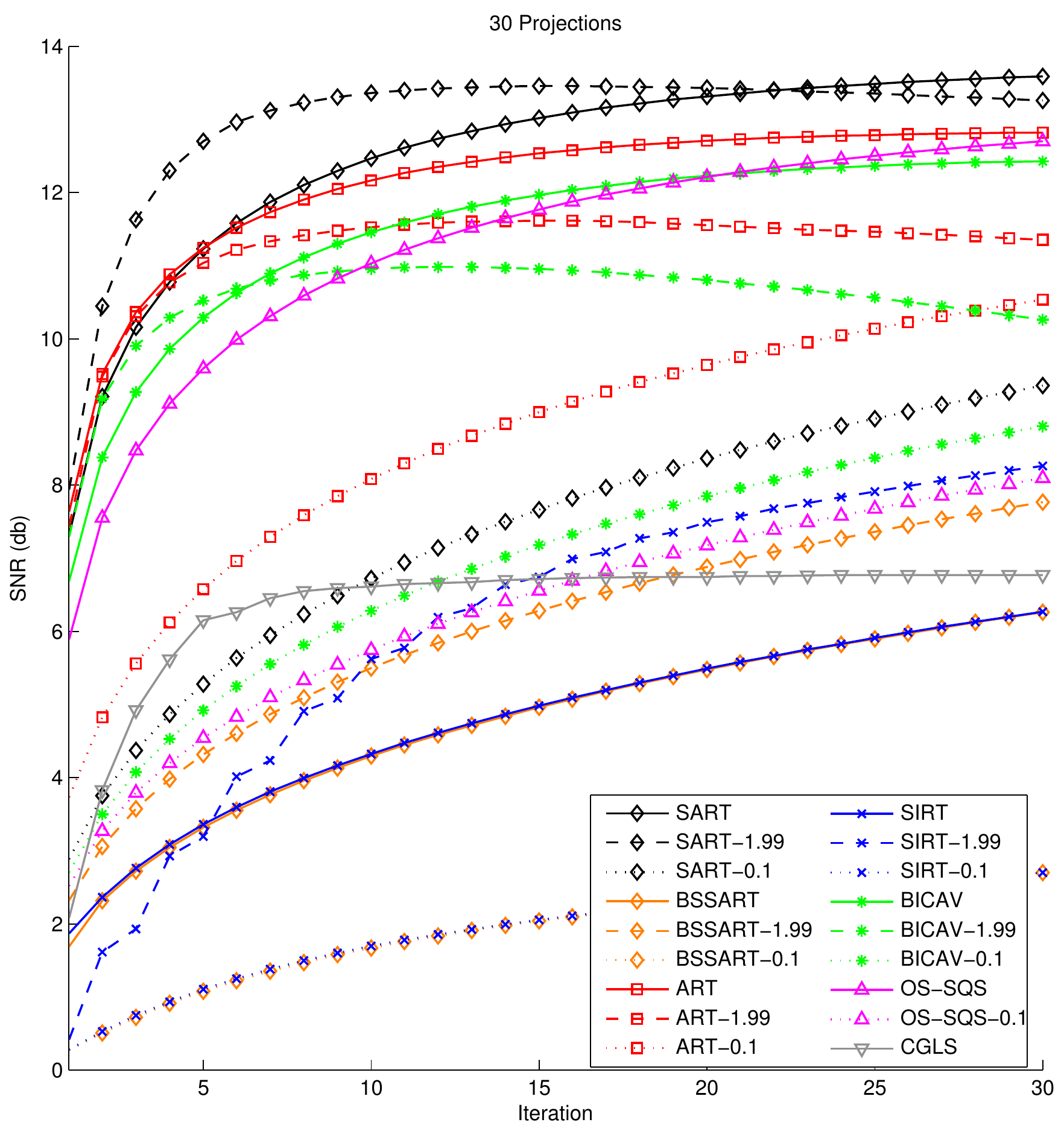}

\includegraphics[width=1\textwidth]{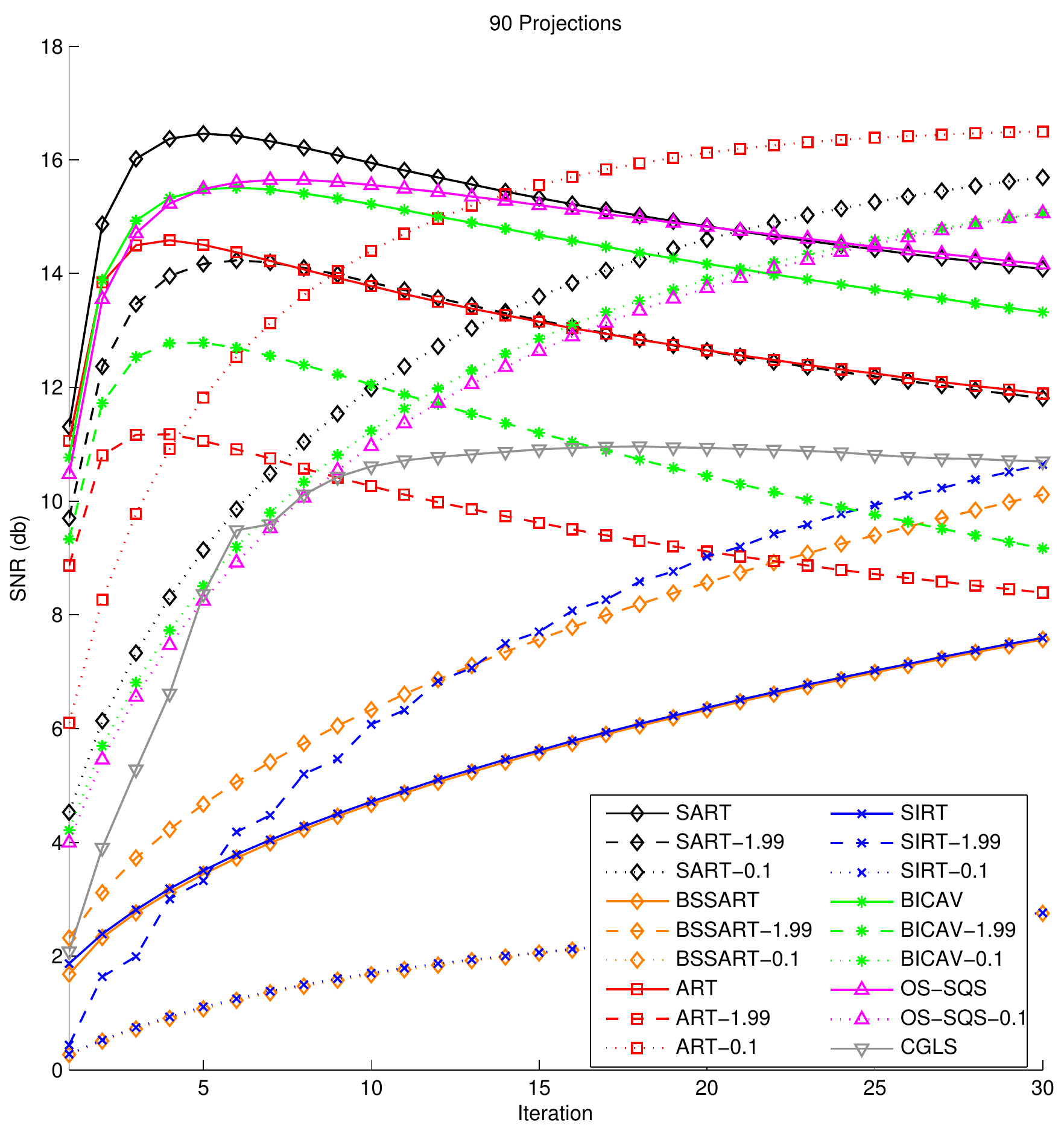}%
\end{minipage}

}\subfloat[NCAT]{%
\begin{minipage}[t]{0.33\textwidth}%
\includegraphics[width=1\textwidth]{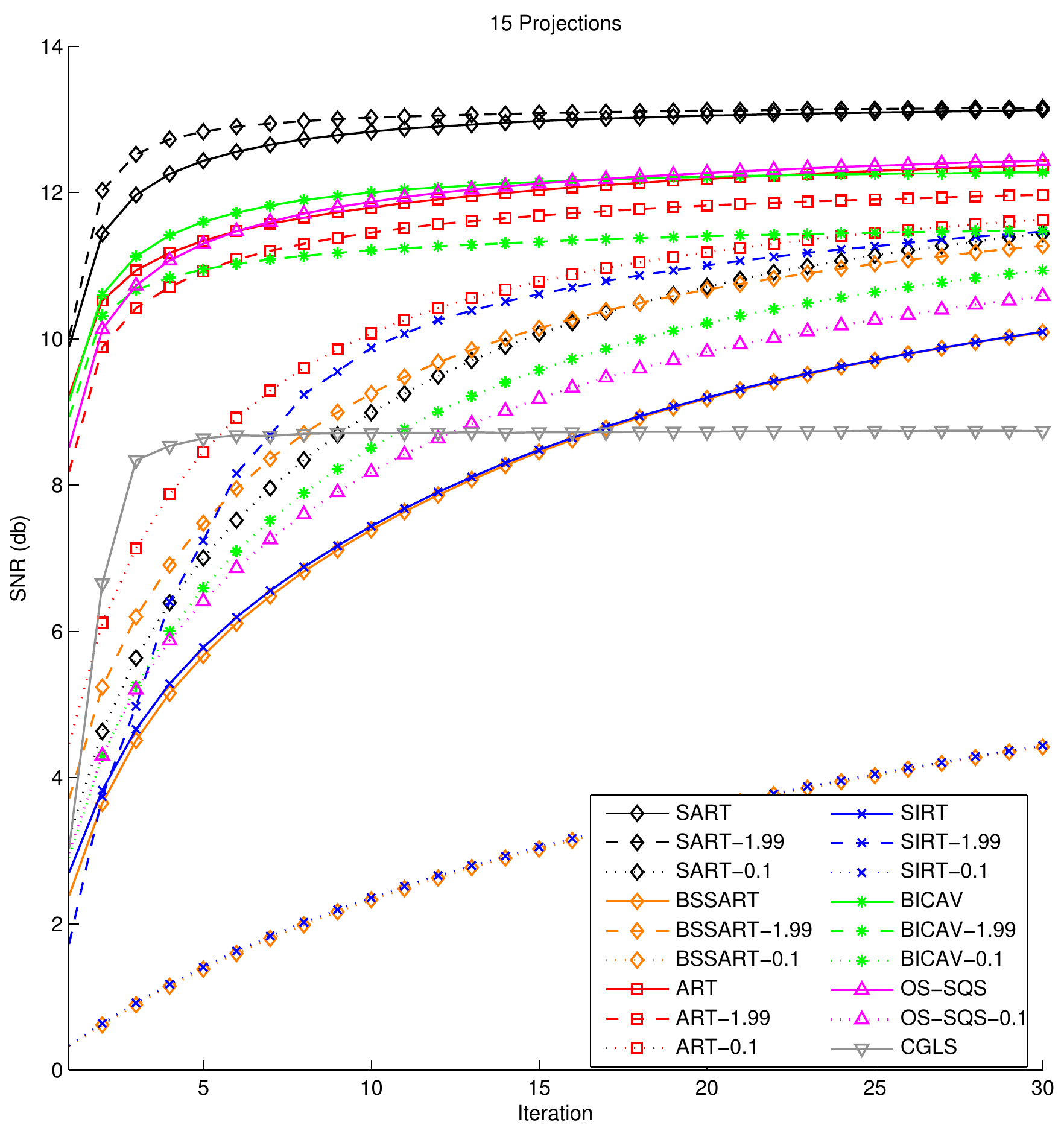}

\includegraphics[width=1\textwidth]{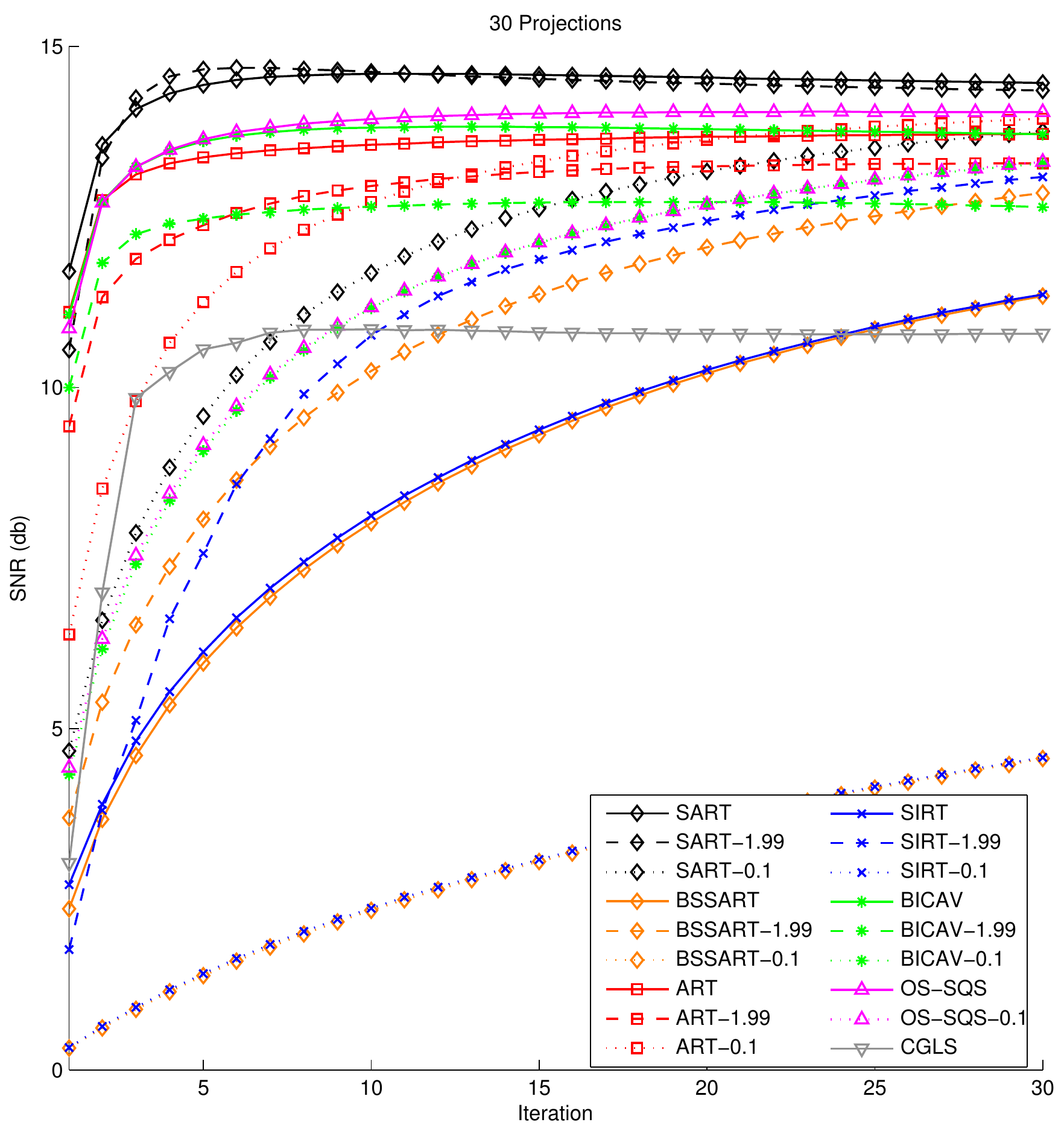}

{\small{}\includegraphics[width=1\textwidth]{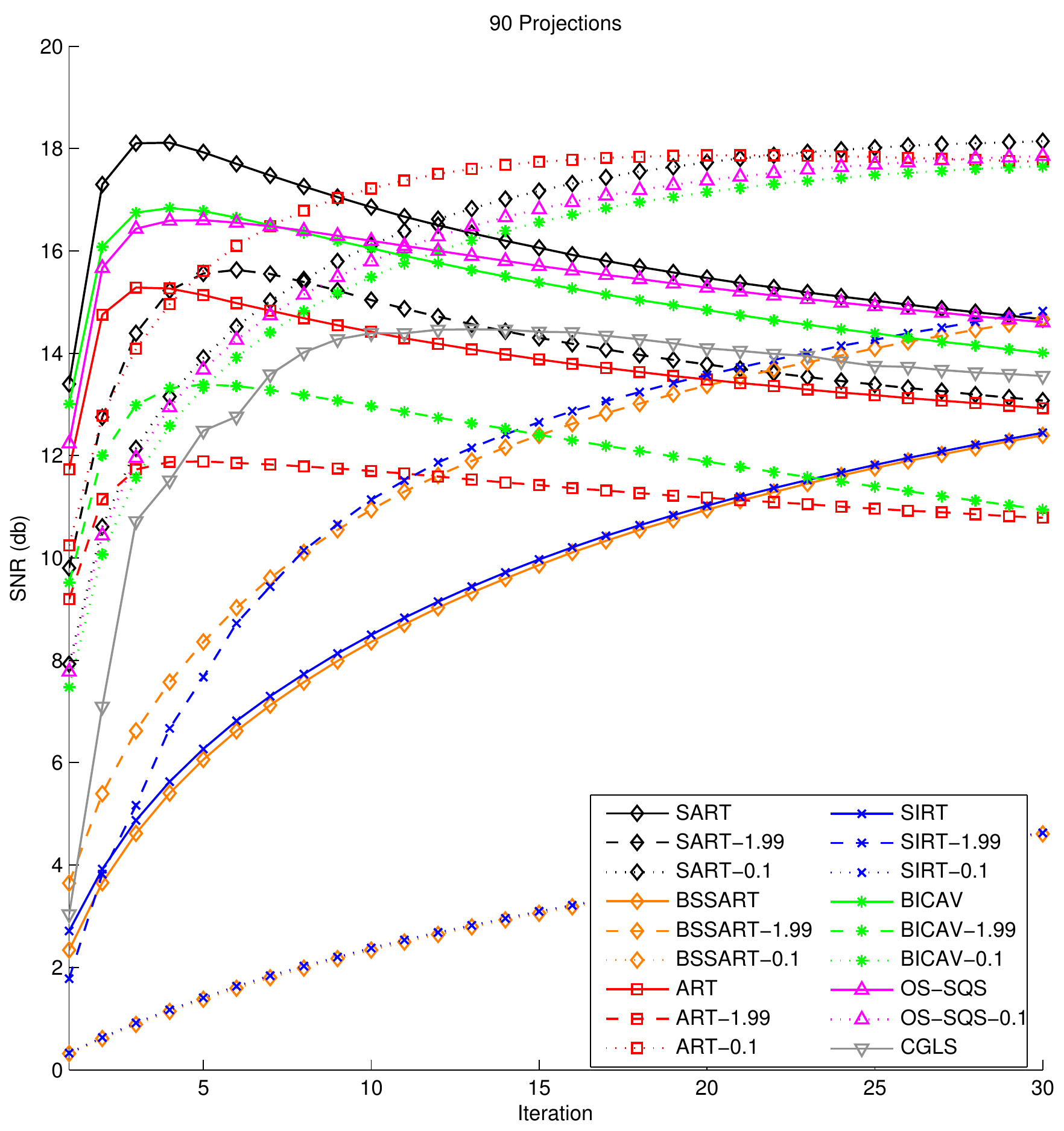}}%
\end{minipage}

}\subfloat[Mouse]{%
\begin{minipage}[t]{0.33\textwidth}%
\includegraphics[width=1\textwidth]{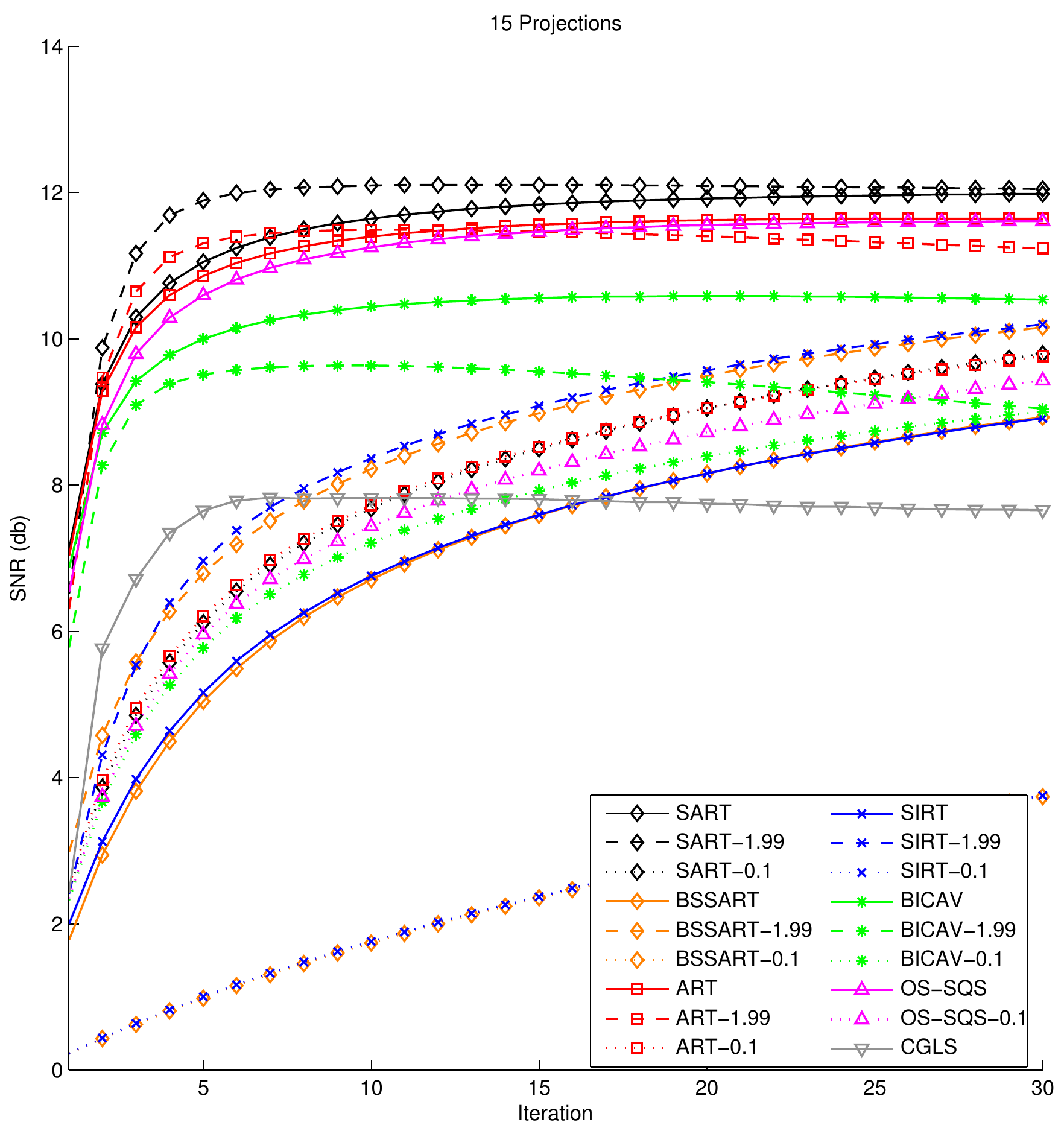}

\includegraphics[width=1\textwidth]{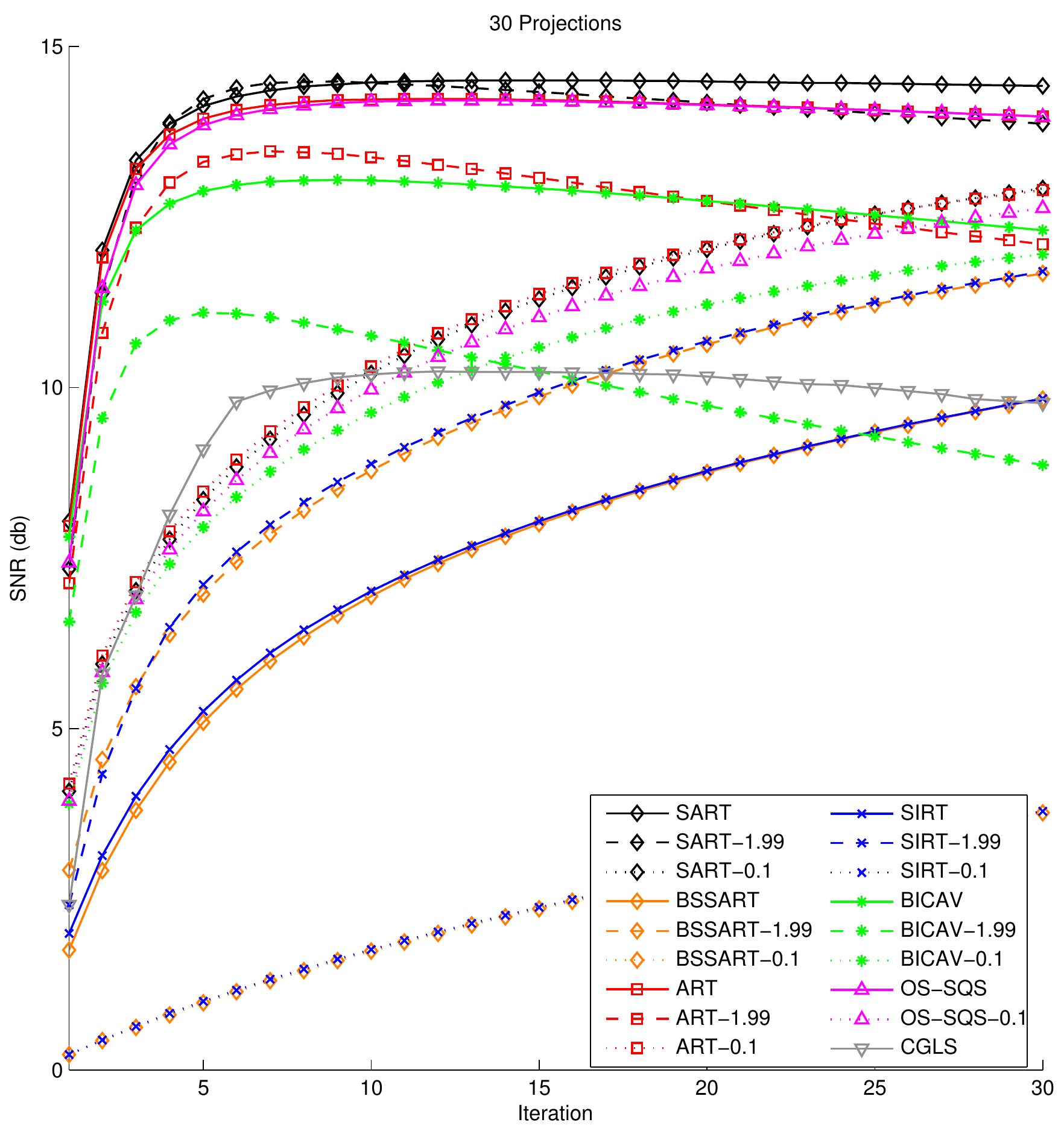}

\includegraphics[width=1\textwidth]{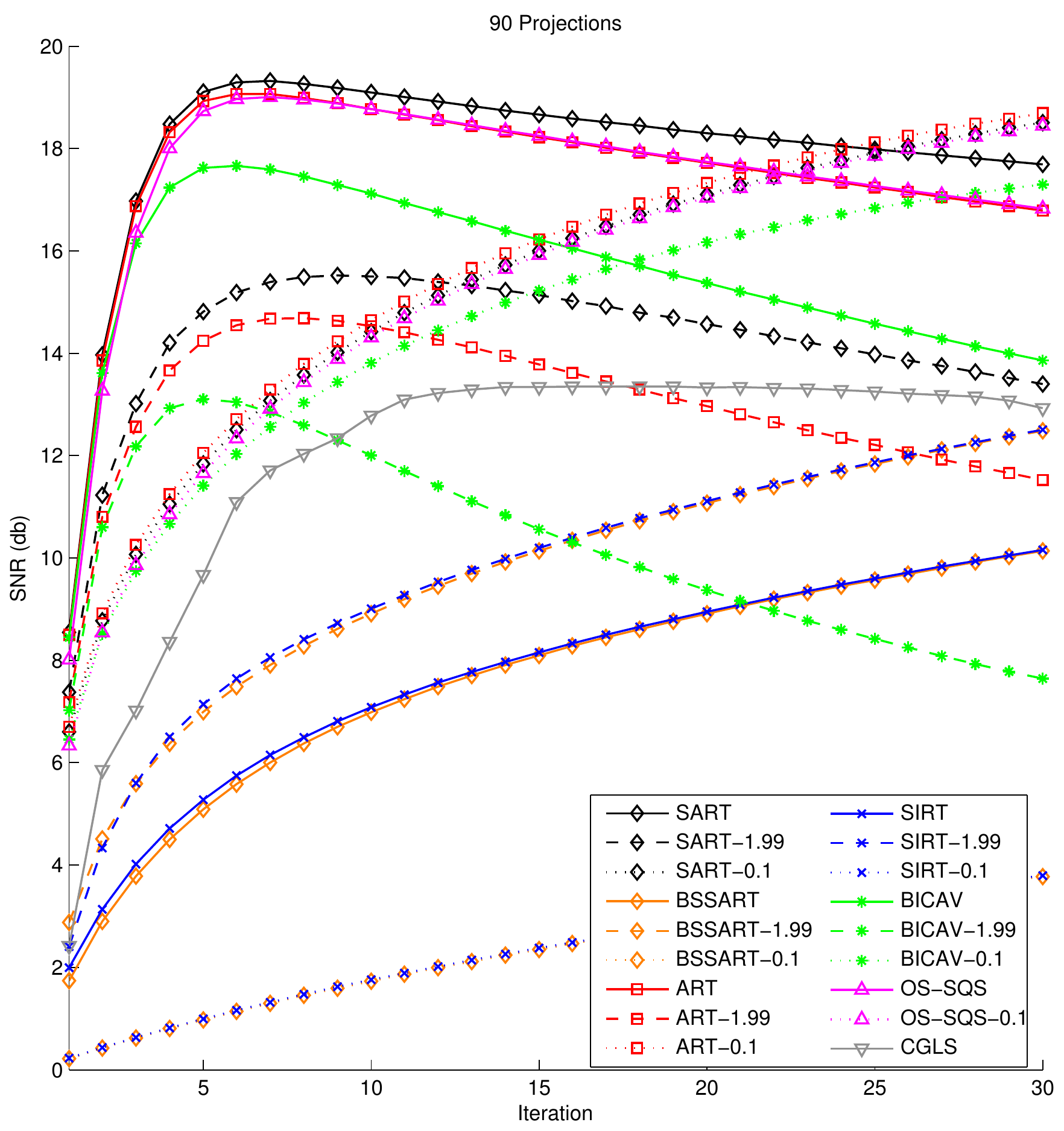}%
\end{minipage}

}

\protect\caption{\textbf{Iterative Algorithms Comparison}. Plots show SNR per iteration.
Solid lines have $\alpha=1$, dashed lines have $\alpha=1.99$, and
dotted lines have $\alpha=0.1$. \label{fig:SNR-per-iteration-iterative}}
\end{figure*}

\begin{figure*}
\center\subfloat[Modified Shepp-Logan]{%
\begin{minipage}[t]{0.33\textwidth}%
\includegraphics[width=1\textwidth]{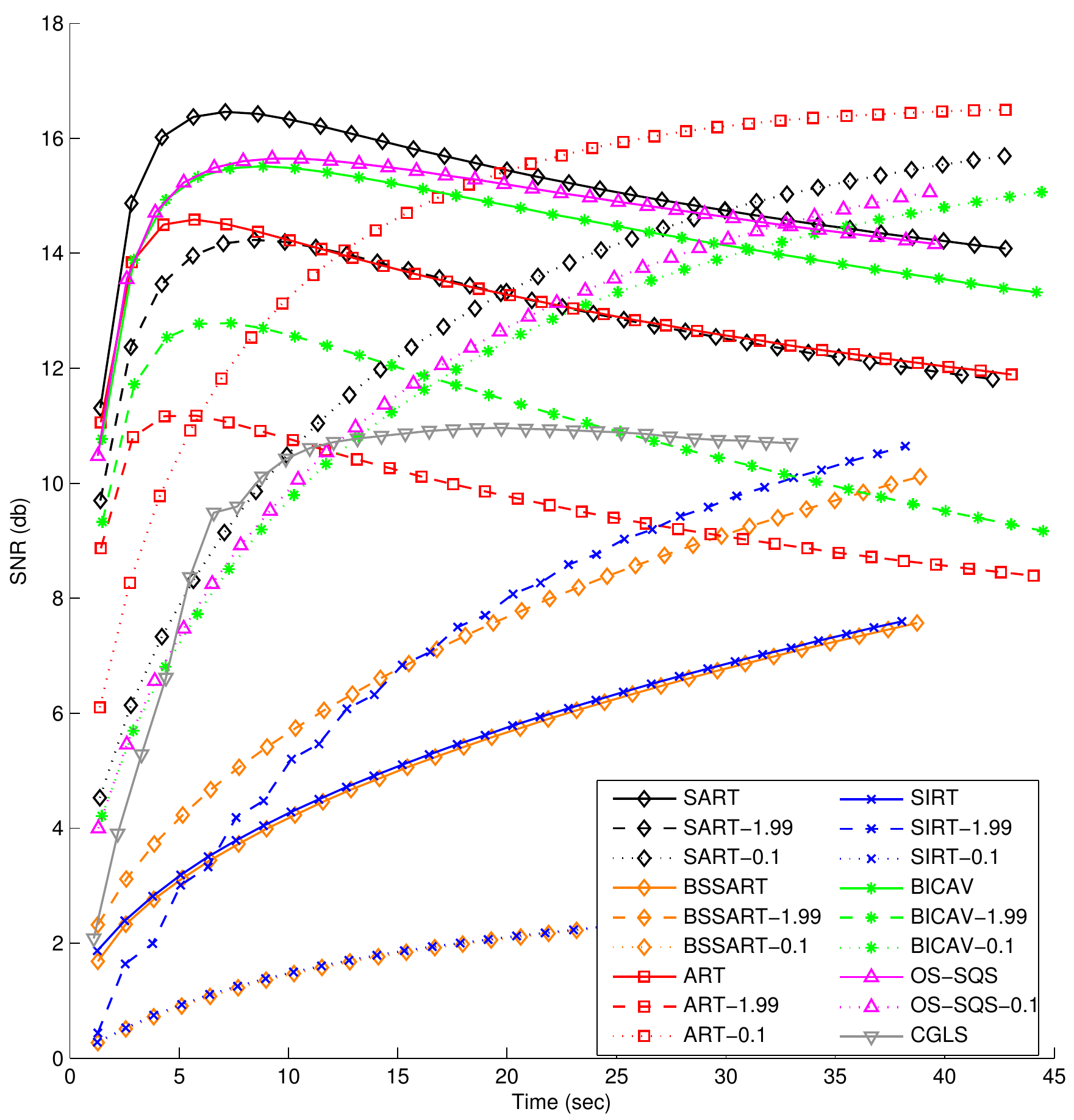}%
\end{minipage}

}\subfloat[NCAT]{%
\begin{minipage}[t]{0.33\textwidth}%
\includegraphics[width=1\textwidth]{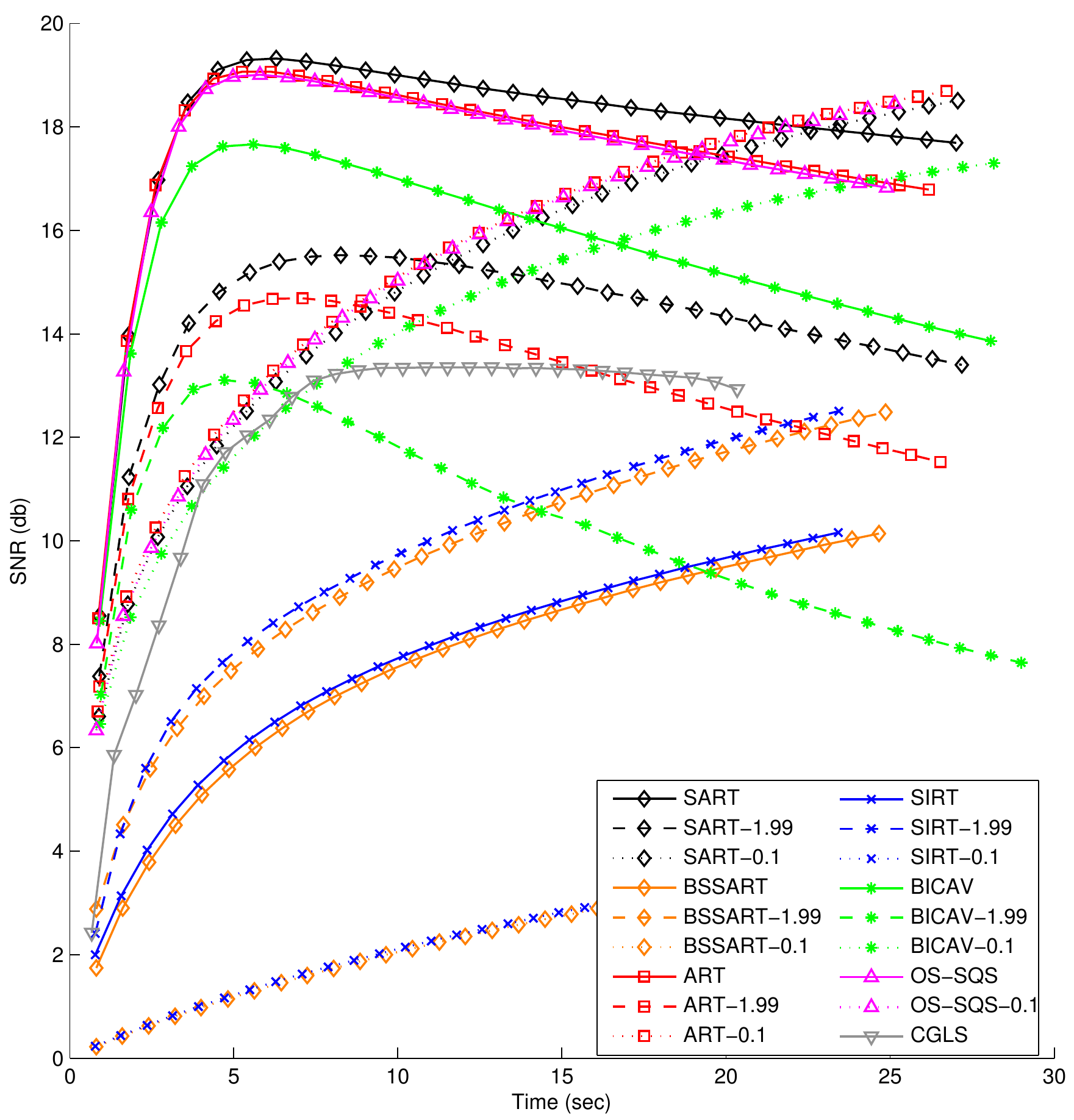}%
\end{minipage}

}\subfloat[Mouse]{%
\begin{minipage}[t]{0.33\textwidth}%
\includegraphics[width=1\textwidth]{plots/per_iter-ph_mouse-512_nt_gauss-nl_0\lyxdot 000-np_90-p_mouse_t1-snr}%
\end{minipage}

}

\protect\caption{\textbf{Running Time Comparison}. Curves show SNR per running time
for 90 projections . Compare with Fig. \ref{fig:SNR-per-iteration-iterative}
(bottom row). \label{fig:SNR-per-time-iterative}}
\end{figure*}

\begin{figure*}
\center\subfloat[Modified Shepp-Logan]{%
\begin{minipage}[t]{0.33\textwidth}%
\includegraphics[width=1\textwidth]{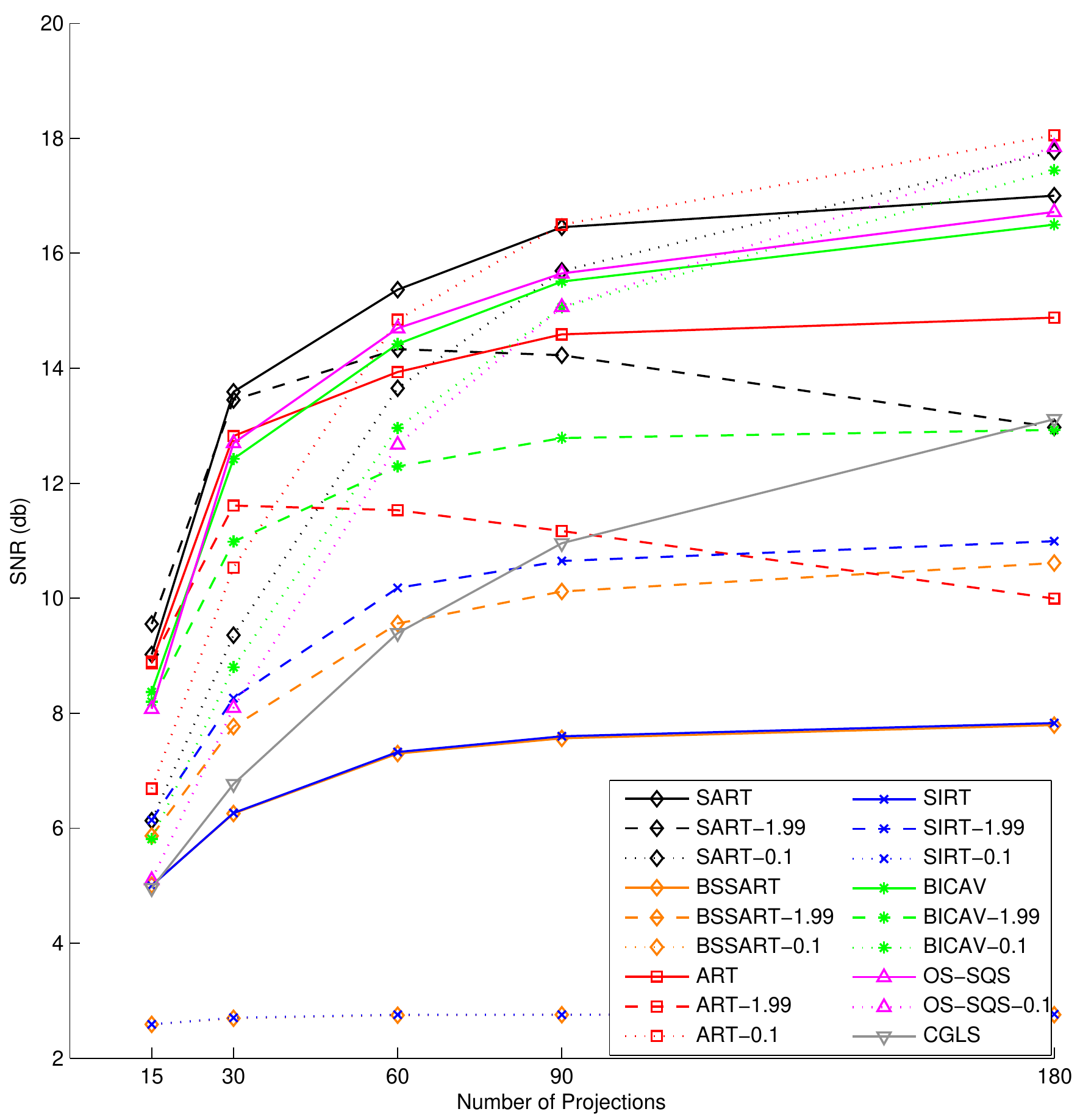}%
\end{minipage}

}\subfloat[NCAT]{%
\begin{minipage}[t]{0.33\textwidth}%
\includegraphics[width=1\textwidth]{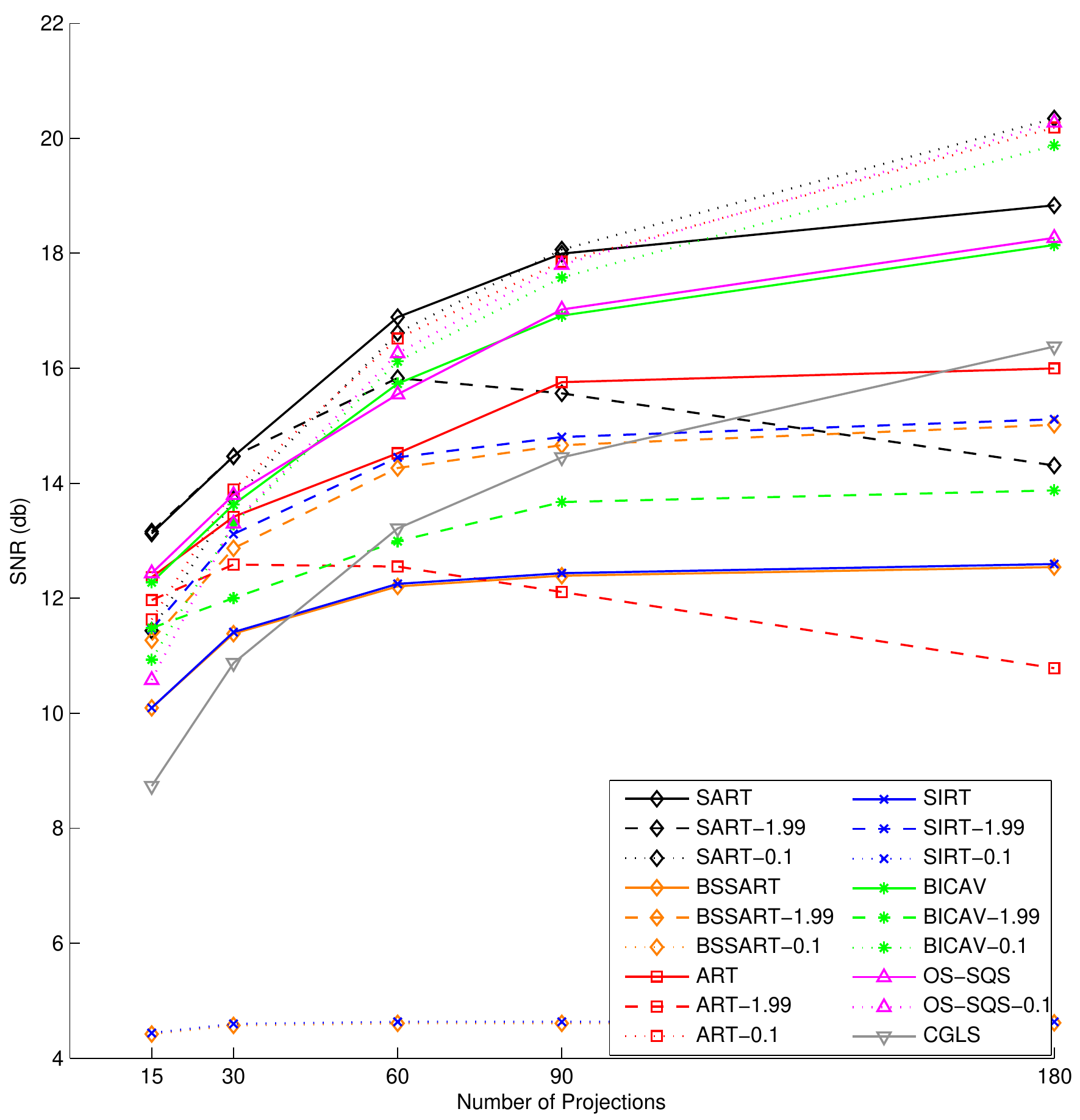}%
\end{minipage}

}\subfloat[Mouse]{%
\begin{minipage}[t]{0.33\textwidth}%
\includegraphics[width=1\textwidth]{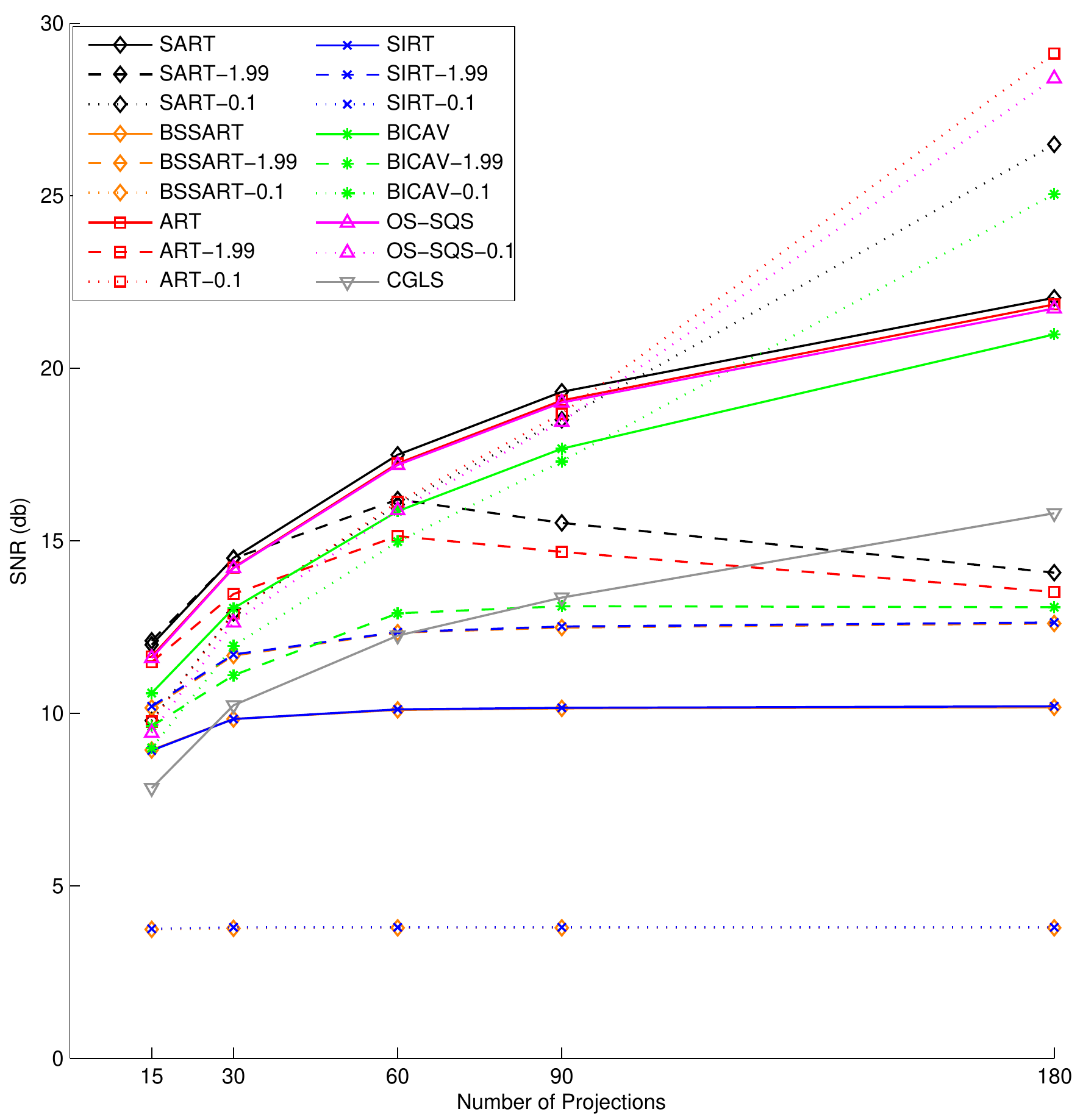}%
\end{minipage}

}

\protect\caption{\textbf{Effect of the number of projections}. Plots show the maximum
SNR achieved over 30 iterations per number of projections used. \label{fig:SNR-per-num-proj-iterative}}
\end{figure*}

We first compare the different iterative algorithms presented in Sec.
\ref{sec:Iterative-Algorithms} on the datasets. We set the number
of subsets in OS-SQS to the number of projections to have a fair comparison
with SART, since we noticed that increasing the number of subsets
increases the convergence rate. We compare different values of $\alpha$,
namely 0.1, 1, and 1.99. We compare convergence per iteration since
all methods are roughly equal in runtime, as each outer iteration
contains (roughly) one forward and one backward projection. This is
confirmed in Fig. \ref{fig:SNR-per-time-iterative}. Note that our
implementation is not optimized for any of the methods, and the processing
time is just an indication. We initialize all methods with uniform
volume $x^{(0)}=\mathbf{0}_{n}$.

Fig. \ref{fig:SNR-per-iteration-iterative} shows the SNR per iteration
for 15, 30, 90  projections for 30 iterations. Fig. \ref{fig:SNR-per-num-proj-iterative}
shows the the maximum SNR over 30 iterations for different number
of equally distributed projections from 15 to 180. From the figures,
we make the following conclusions:
\begin{itemize}
\item The simulated projections closely resemble the results from the real
dataset, which suggests that the measurement noise model is reflective
of real data.
\item With fewer projections (15 or 30 projections), using larger values
$\alpha=1.99$ gives faster and better convergence. 
\item With many projections, moderate values $\alpha=1$ produces a fast
convergence that then falls off and is overtaken by $\alpha=0.1$.
\item SART provides the fastest convergence within a handful of iterations,
and is consistently better for fewer projections. However, it is overtaken
by ART and others for many projections. This provides the motivation
to use it in the proximal framework, since typically the tomography
solver is invoked for only a few iterations per outer iteration of
ADMM for example \cite{ramani2012splitting}.
\item With more projections, e.g. 90, we notice that the SNR for a few methods
go up and then down. This doesn't mean, however, that they are not
converging. The objective function is the reprojection error not the
SNR. This can be explained by the fact of the presence of noise, and
that at some point the algorithm starts fitting the noise in the measurements
\cite{herman2009fundamentals}. Usually these kinds of algorithms
are run interactively where the user inspects the reconstruction quality
every few iterations and stops the procedure when it starts to deteriorate,
which motivates Fig. \ref{fig:SNR-per-num-proj-iterative}.
\item Even though BICAV, SIRT, and BSSART have formal proofs of convergence,
their convergence speed per iteration is in fact much lower than SART
or (this version of) OS-SQS, that lack these proofs.
\item The faster convergence and best results are achieved by SART, followed
by ART, OS-SQS, and BICAV. They work better with $\alpha=1$ for few
projections, and with $\alpha=0.1$ for more projections.
\item CGLS, that was used before for solving tomography problems \cite{ramani2012splitting,sidky2012convex},
performs quite poorly compared to the other iterative algorithms.
\item Using plain iterative methods does not give acceptable results with
fewer projections. Thus we focus next on using regularizers in the
proximal framework with SART, ART, OS-SQS, and BICAV and fewer projections,
namely 30 projections.
\end{itemize}

\subsection{Poisson Model Mapping Functions Comparison\label{sub:Poisson-Mapping-Functions-Comparison}}

\begin{figure*}
\subfloat[Modified Shepp-Logan]{%
\begin{minipage}[t]{0.33\textwidth}%
\includegraphics[width=1\textwidth]{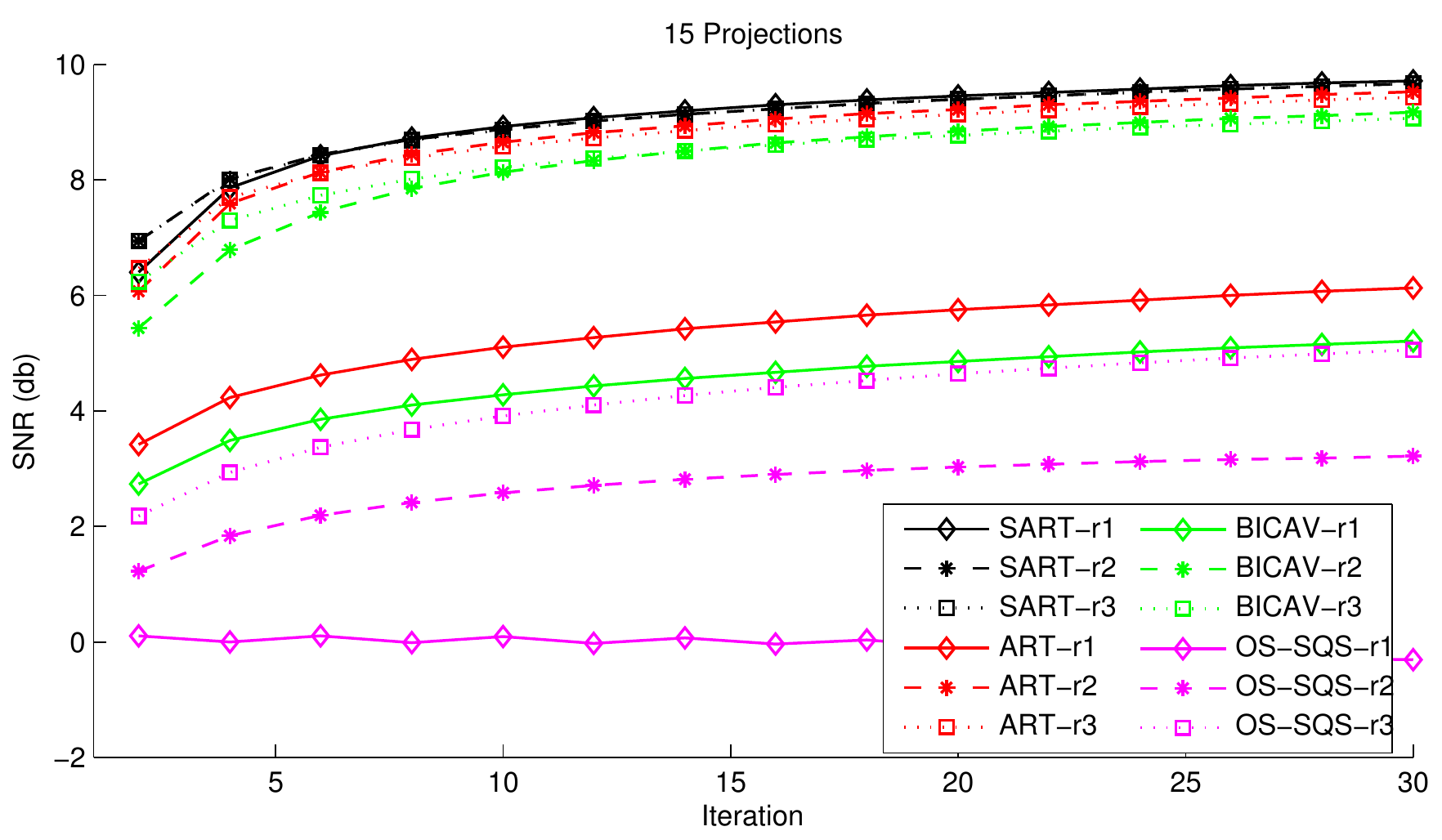}

\includegraphics[width=1\textwidth]{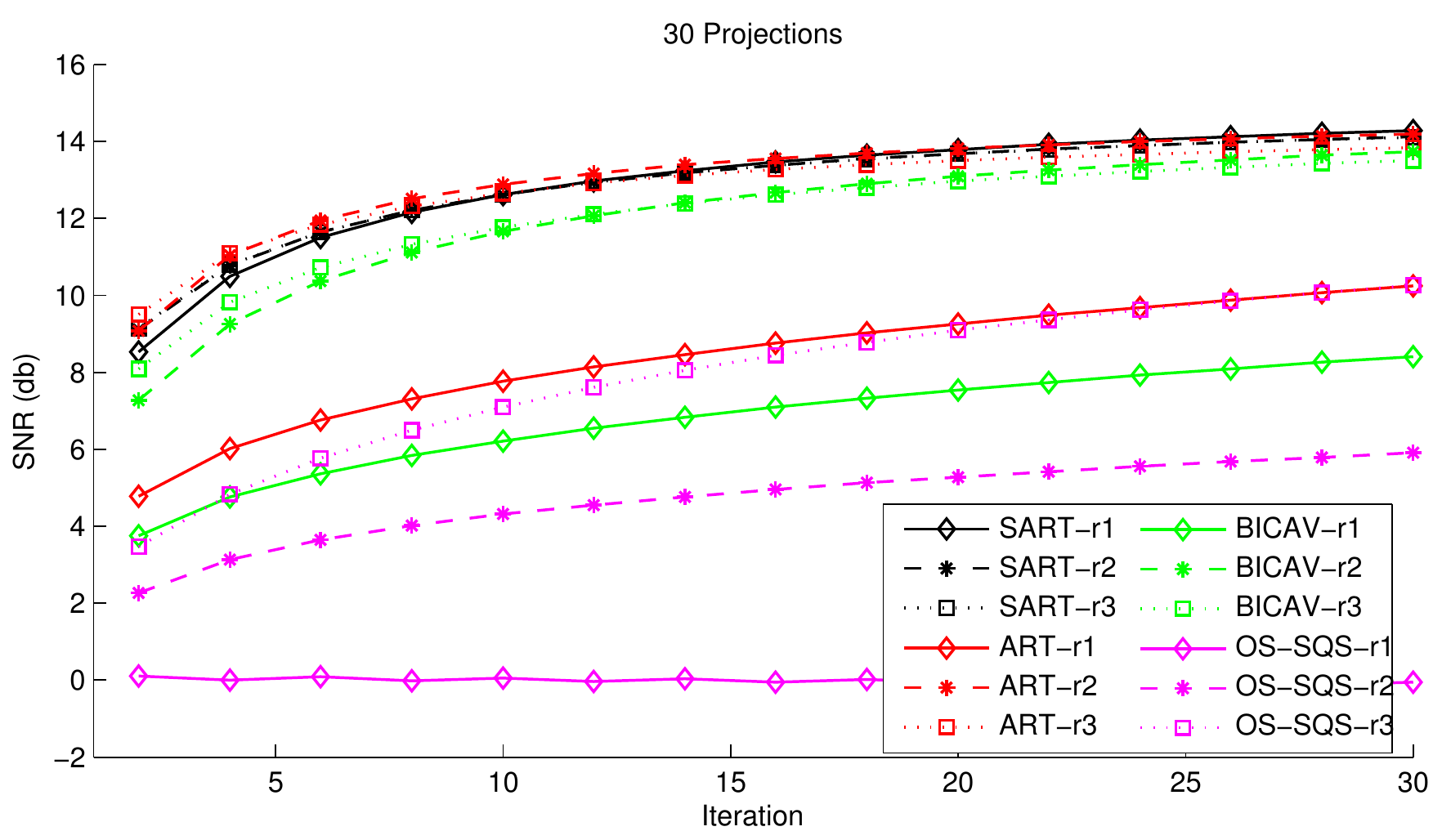}%
\end{minipage}

}\subfloat[NCAT]{%
\begin{minipage}[t]{0.33\textwidth}%
\includegraphics[width=1\textwidth]{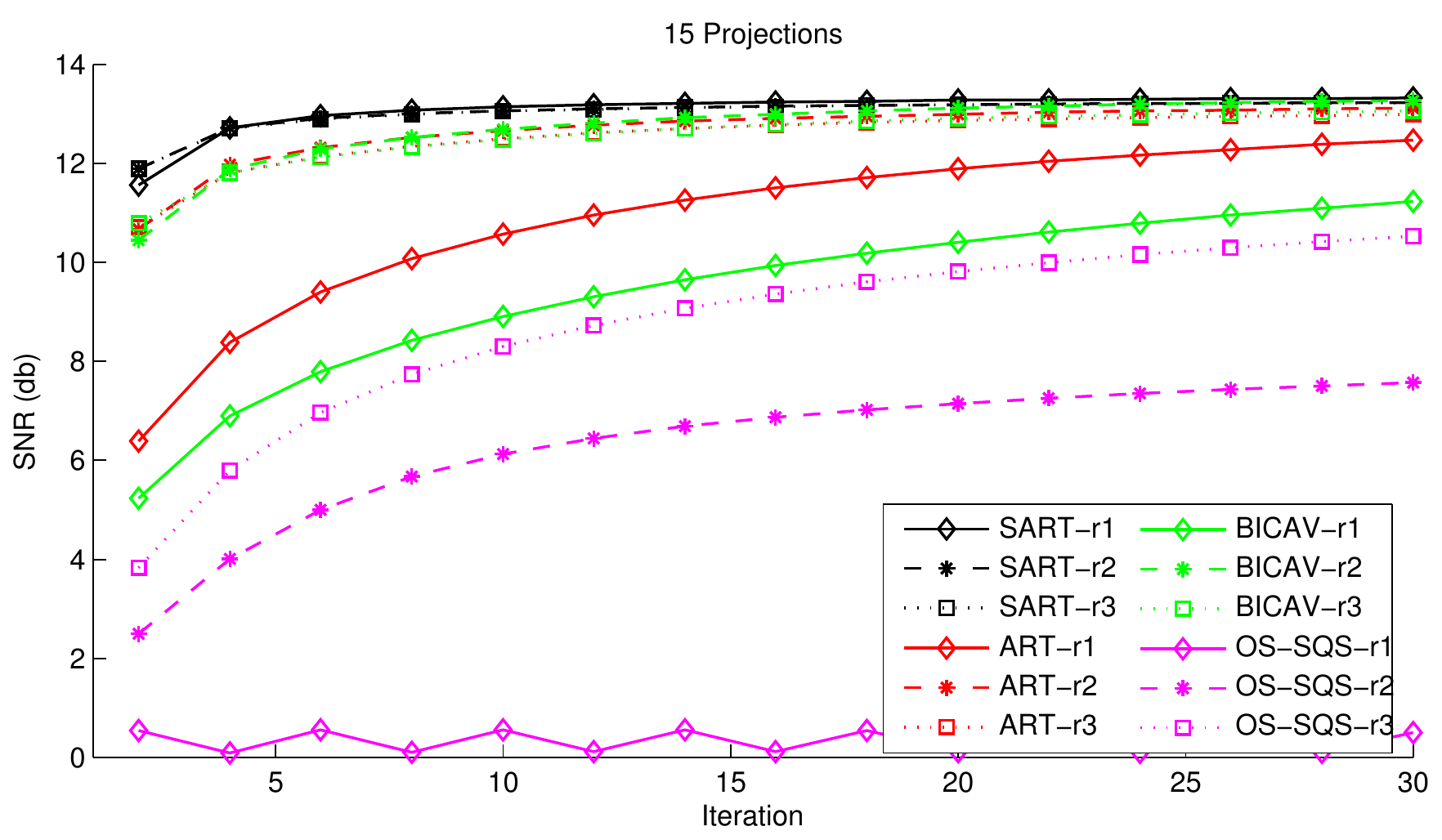}

\includegraphics[width=1\textwidth]{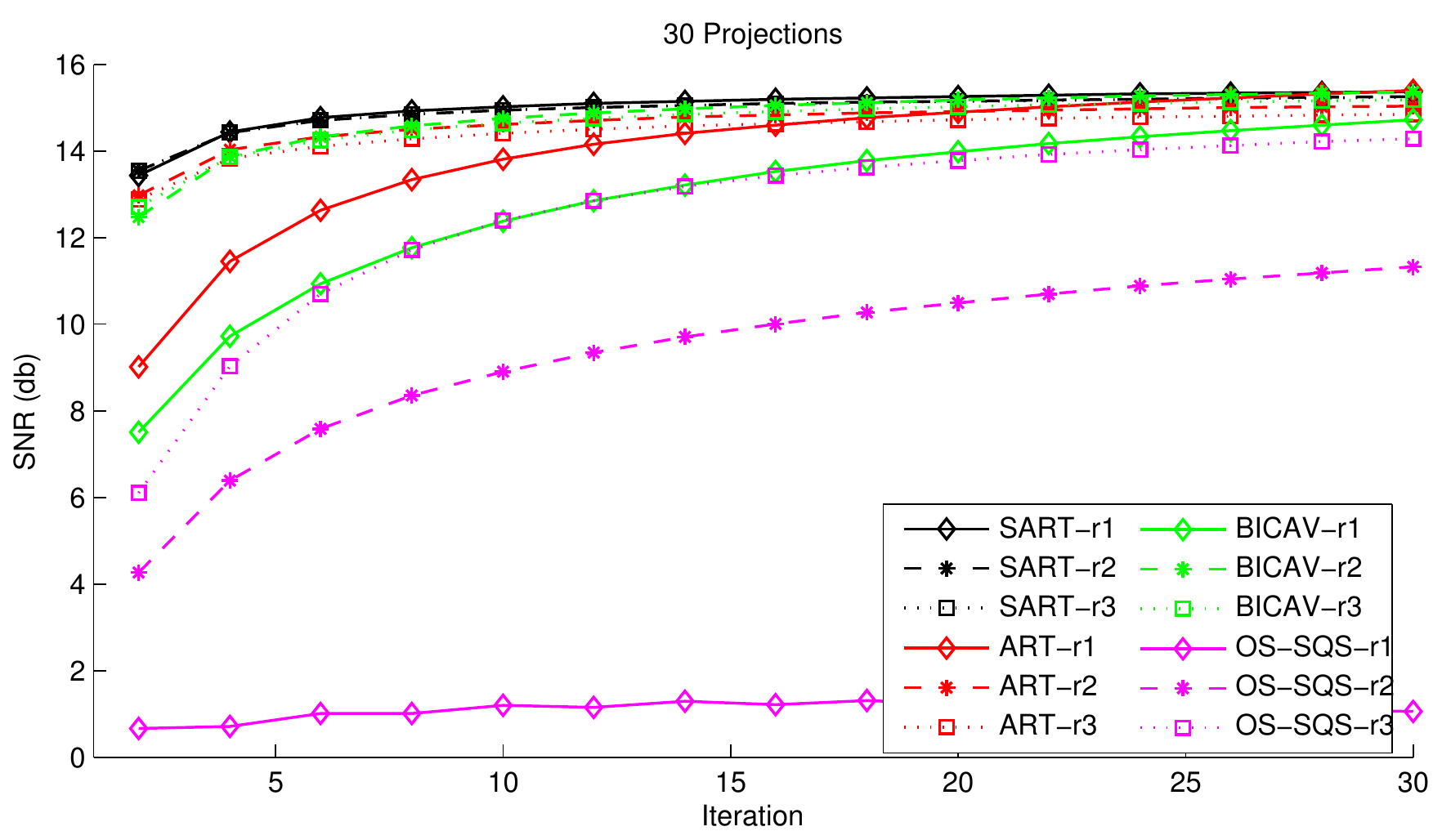}%
\end{minipage}

}\subfloat[Mouse]{%
\begin{minipage}[t]{0.33\textwidth}%
\includegraphics[width=1\textwidth]{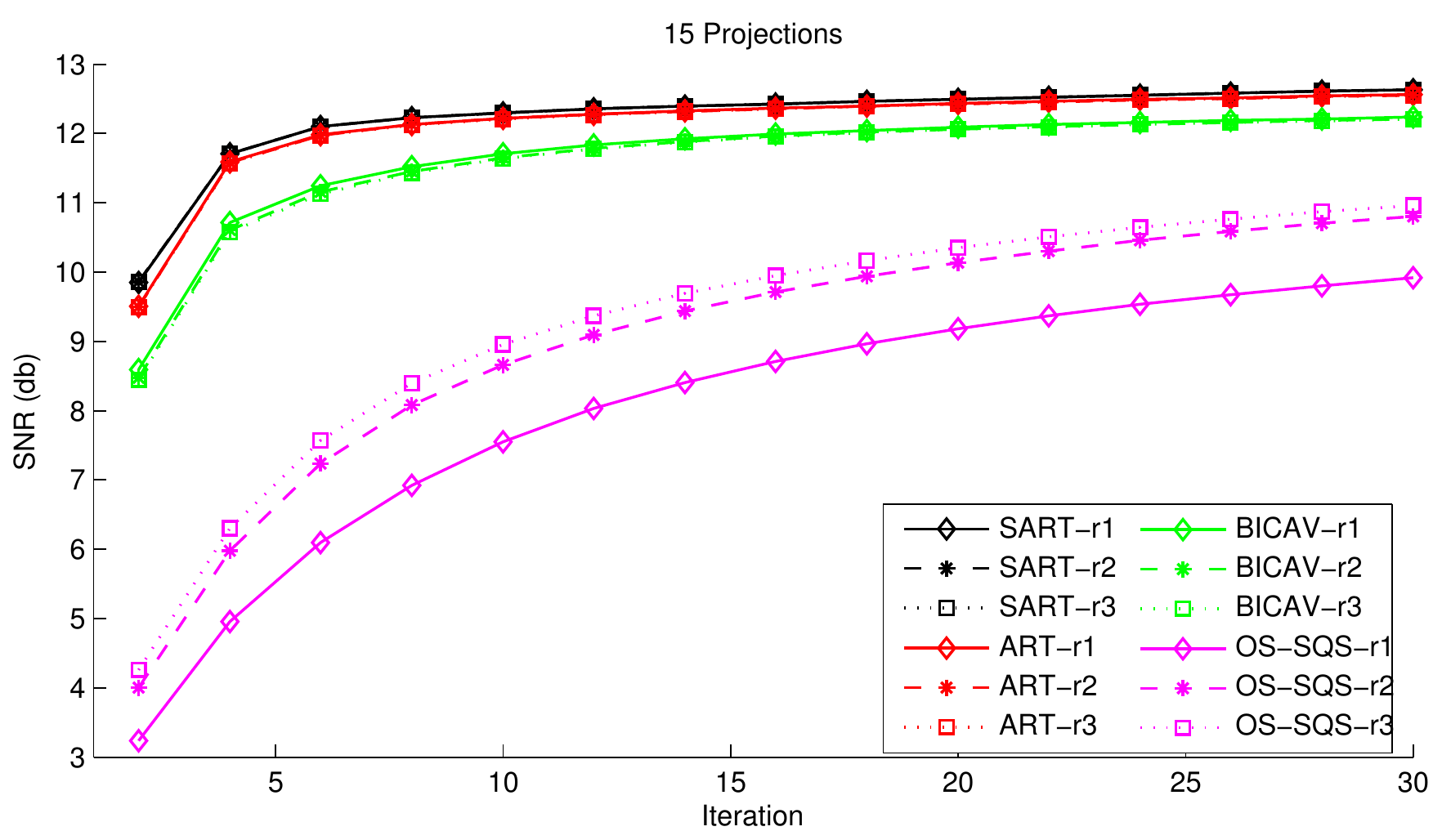}

\includegraphics[width=1\textwidth]{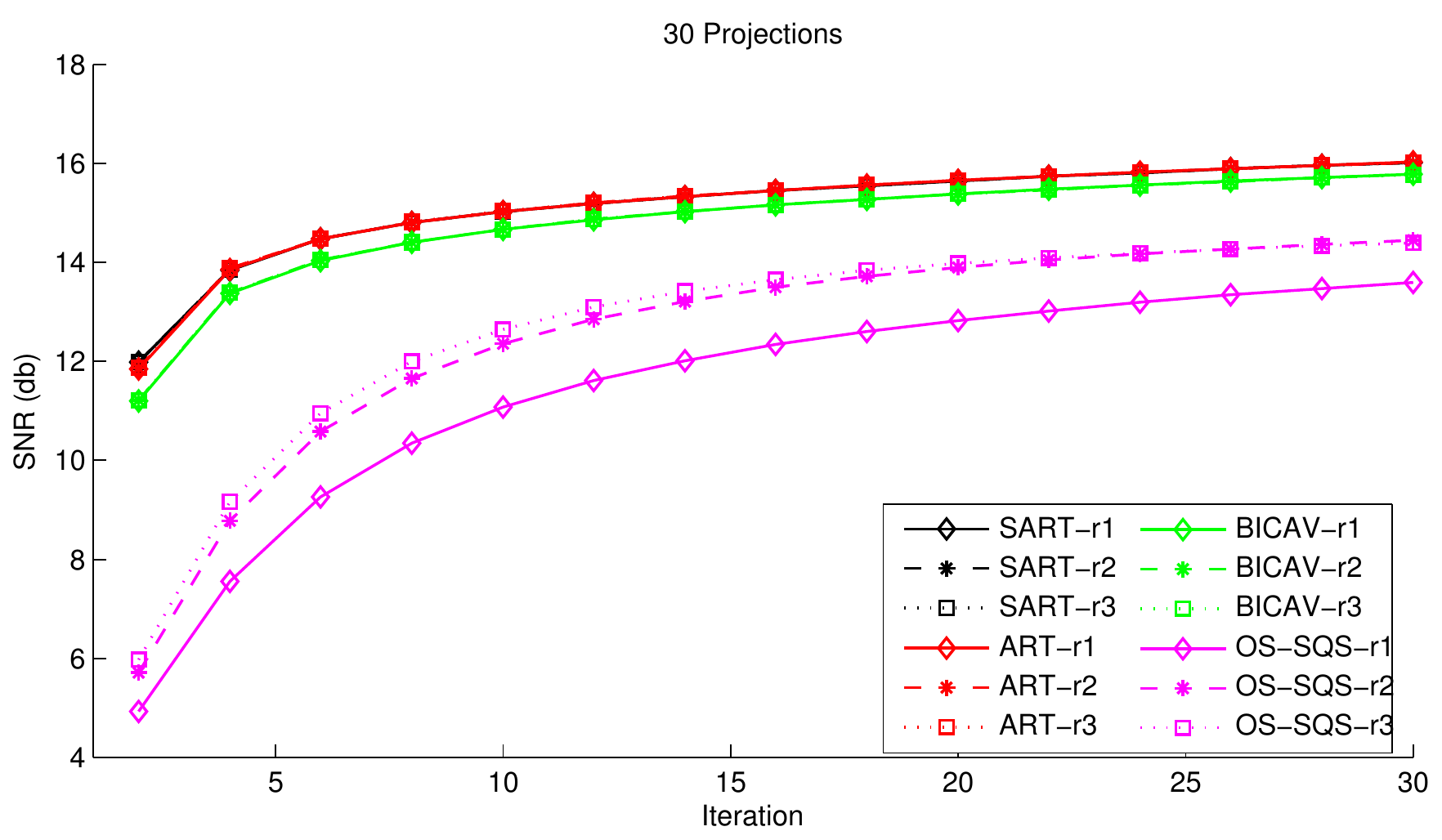}%
\end{minipage}

}

\protect\caption{\textbf{TRex Poisson Model Mapping Functions Comparison}. Curves show
SNR per iteration for the Poisson noise models with mapping functions
$r_{1}$ (\emph{diamonds}), $r_{2}$ (\emph{asterisks}), and $r_{3}$
(\emph{square}s). See Sec. \ref{sub:Poisson-Mapping-Functions-Comparison}.
\label{fig:SNR-per-iter-data-term-l2-prior}}
\end{figure*}

We investigate different mapping functions for the Poisson noise model
$f_{\text{P}}(\cdot)$ in Eq. \ref{eq:Poisson-data-term} using the
ITV regularizer from Eq. \ref{eq:h-ITV} using the ADMM algorithm.
We compare applying different mapping functions on the weights, since
we noticed it improves the performance for some proximal operators.
In particular, we try three functions: identity $r_{1}(w_{i})=w_{i}$,
the square root $r_{2}=\sqrt{w_{i}}$, and the cubic root $r_{3}=\sqrt[3]{w_{i}}$.
Figure \ref{fig:SNR-per-iter-data-term-l2-prior} shows results for
the three datasets for 15 and 30 projections with $\sigma=0.1$ in
Eq. \ref{eq:h-ITV}. We set $\rho=100$ and $\mu=\nicefrac{1}{\rho\Vert K\Vert^{2}}$
except for OS-SQS which was tuned manually as this default value didn't
provide good performance. 

We note the following that using mapping $r_{1}$ is generally worse
than $r_{2}$ and $r_{3}$. This is especially true for OS-SQS, BICAV,
and ART. We believe this is due to the normalizing matrix $C$, where
in this case it includes sum of \emph{squares} of entries in the matrix
$A$ that are typically $<1$. This makes them even smaller, and taking
the square or cubic root of the weights, which are also $<1$, makes
them bigger to counter-balance the former effect, and make the two
terms of the optimization problem in Eq. \ref{eq:prox-operator} of
the same order. This is not the case in SART where the $C$ matrix
contains sums of entries of $A$. We also note that ART and SART provide
very similar performance, closely followed by BICAV and then OS-SQS.
This is also confirmed in the comparison in Sec. \ref{sub:Proximal-Operators-Comparison}.

\subsection{Data Terms and Regularizers Comparison\label{sub:Data-Terms-and-Regularizers-Comparison}}

\begin{figure*}
\center\subfloat[Modified Shepp-Logan]{%
\begin{minipage}[t]{0.33\textwidth}%
\includegraphics[width=1\textwidth]{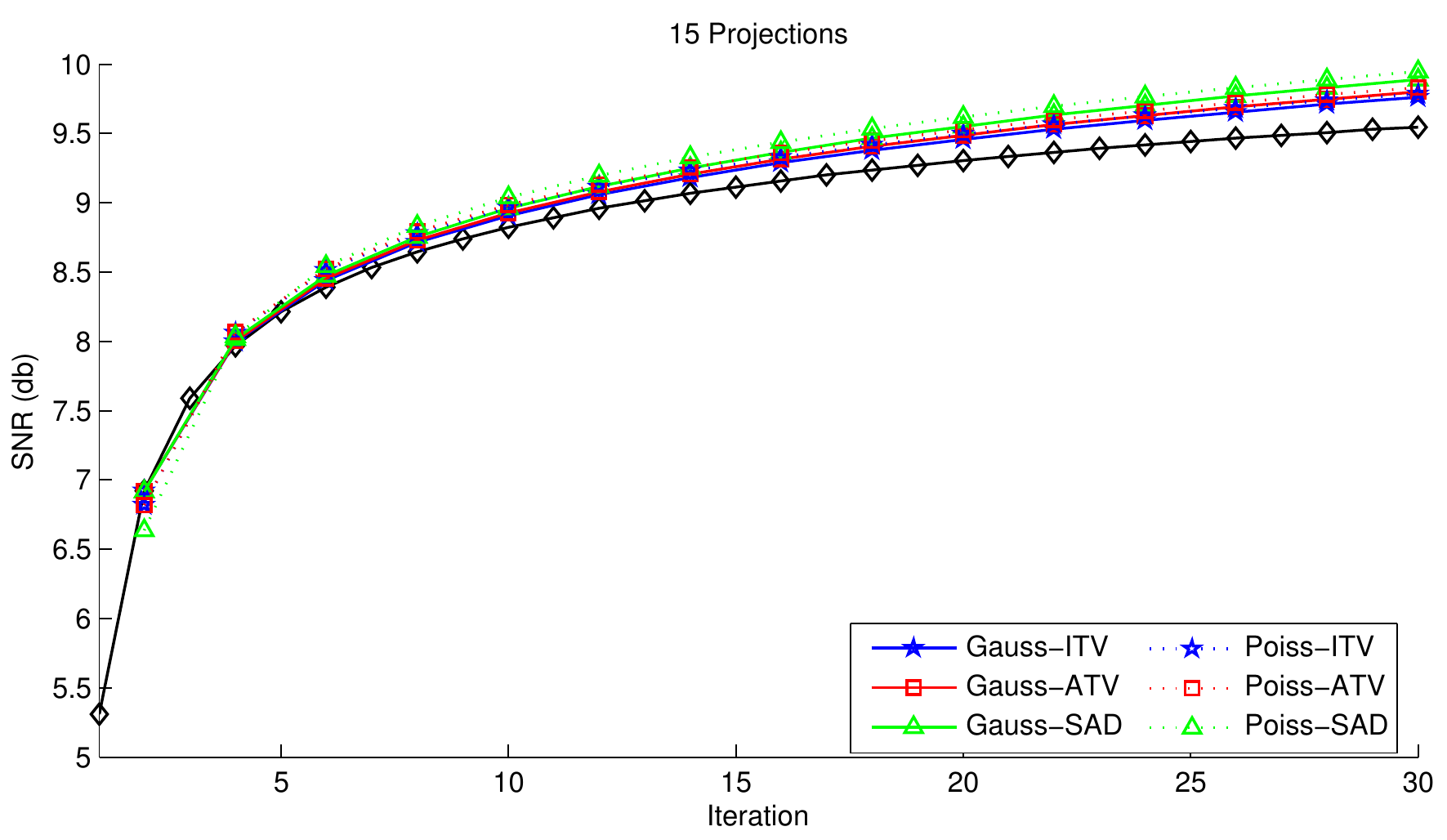}

\includegraphics[width=1\textwidth]{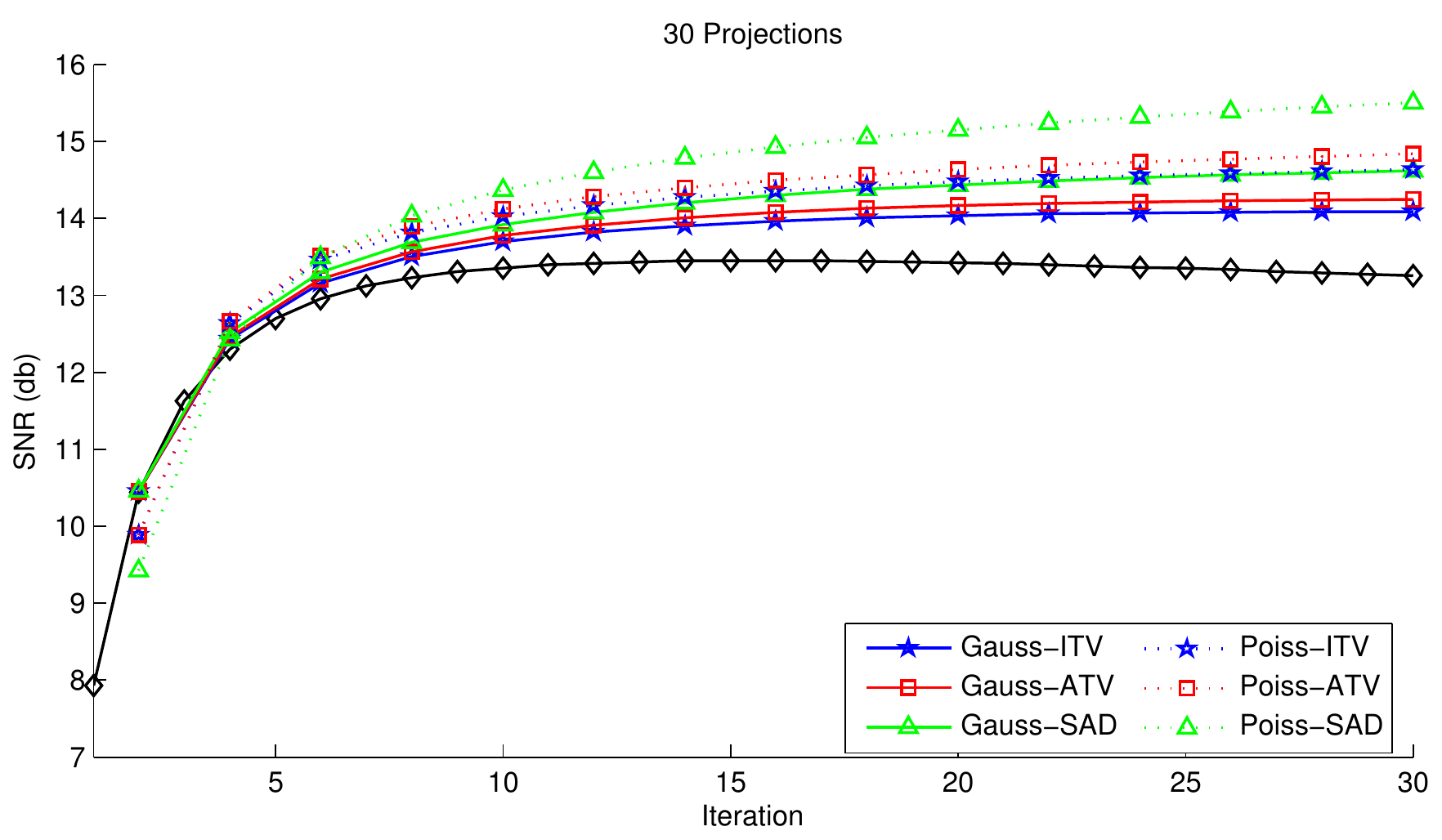}%
\end{minipage}

}\subfloat[NCAT]{%
\begin{minipage}[t]{0.33\textwidth}%
\includegraphics[width=1\textwidth]{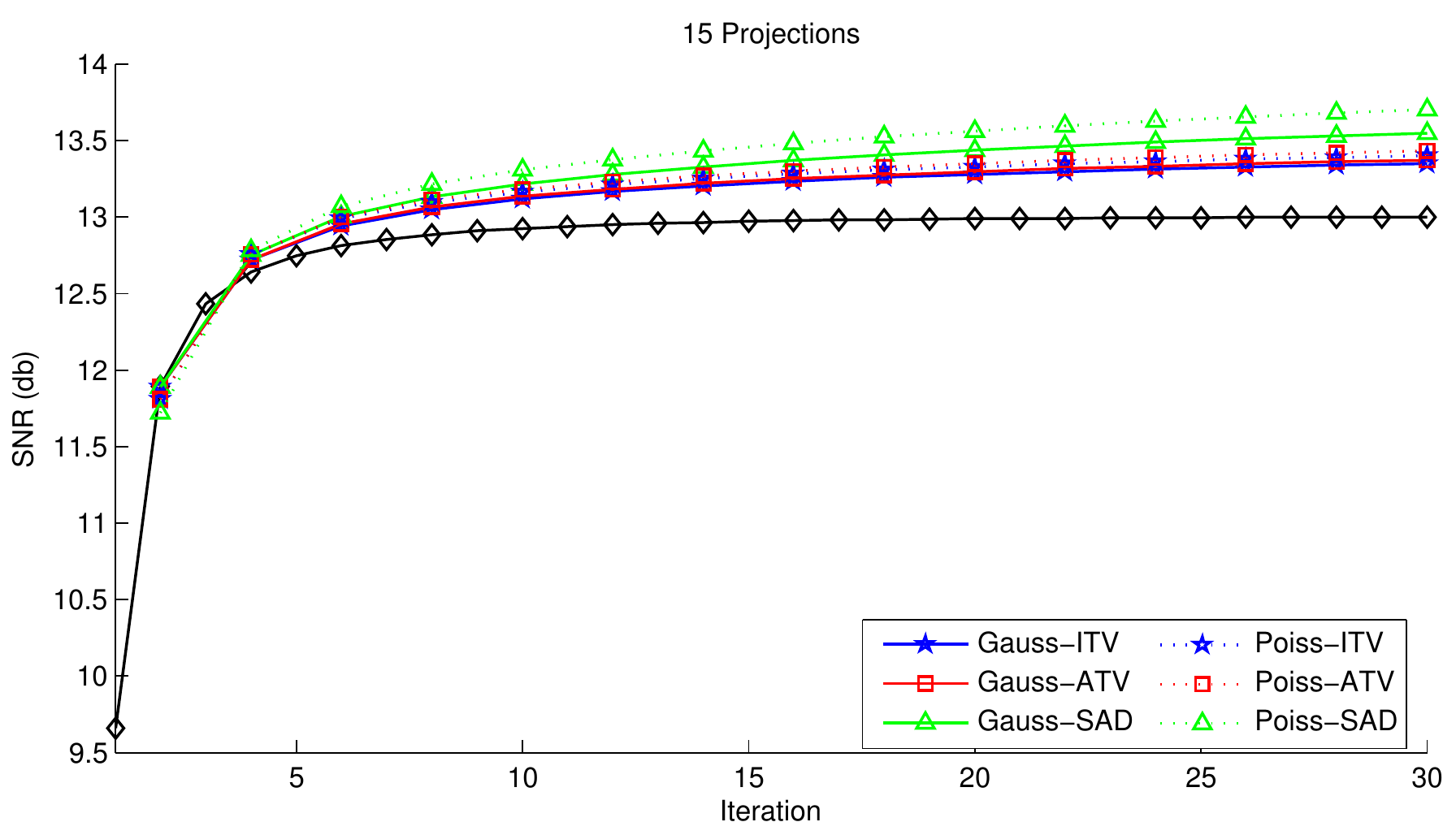}

\includegraphics[width=1\textwidth]{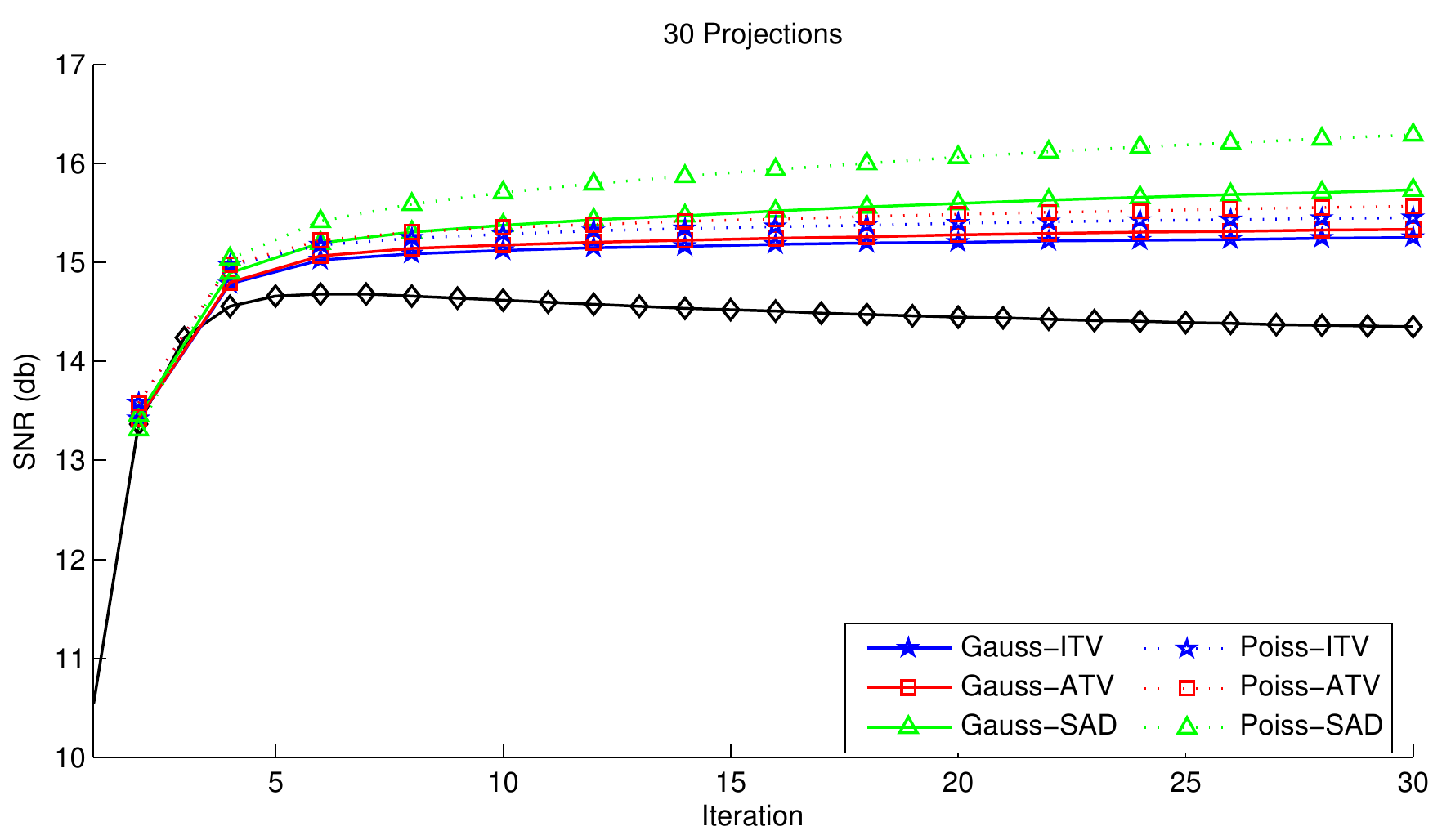}%
\end{minipage}

}\subfloat[Mouse]{%
\begin{minipage}[t]{0.33\textwidth}%
\includegraphics[width=1\textwidth]{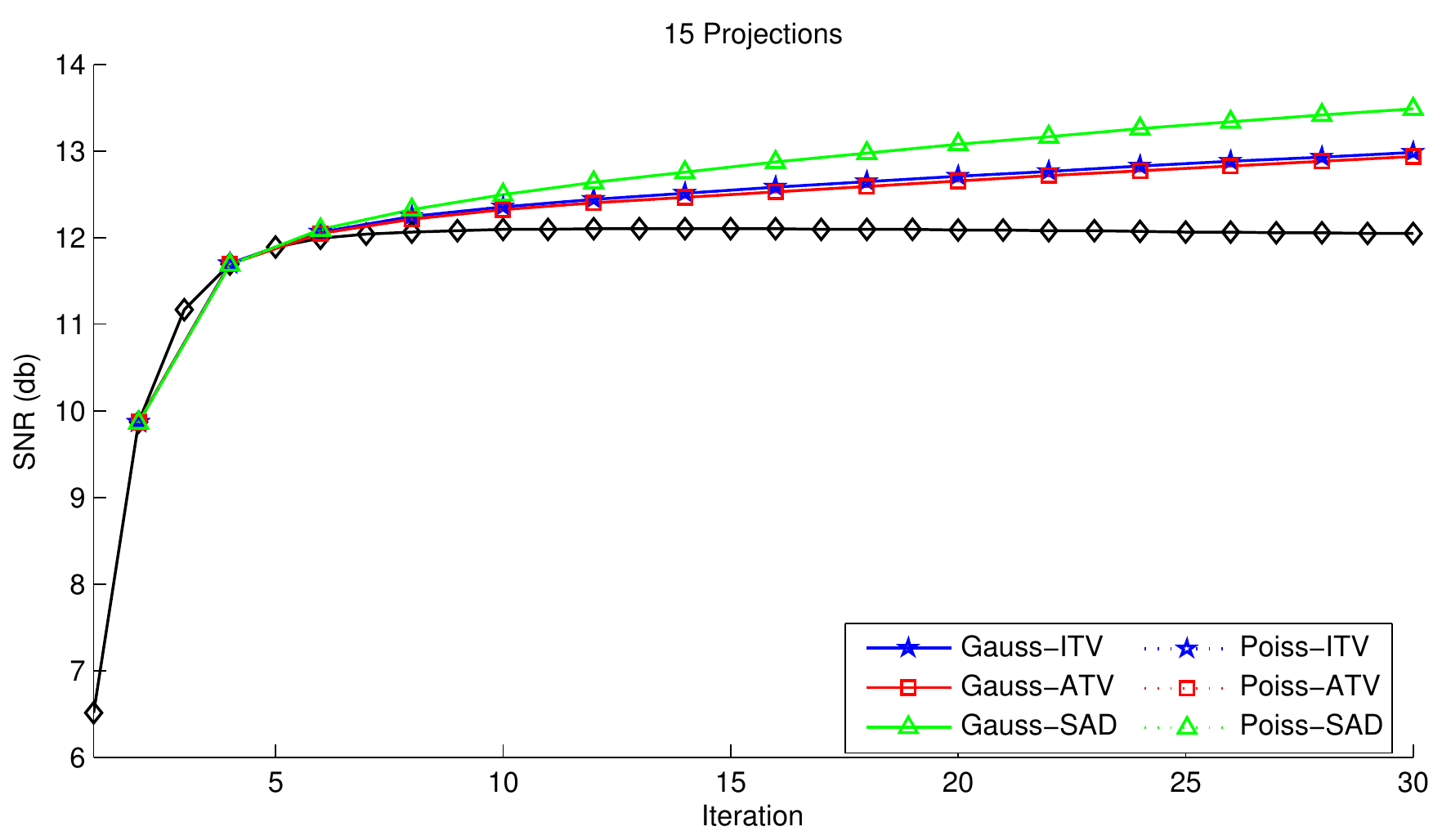}

\includegraphics[width=1\textwidth]{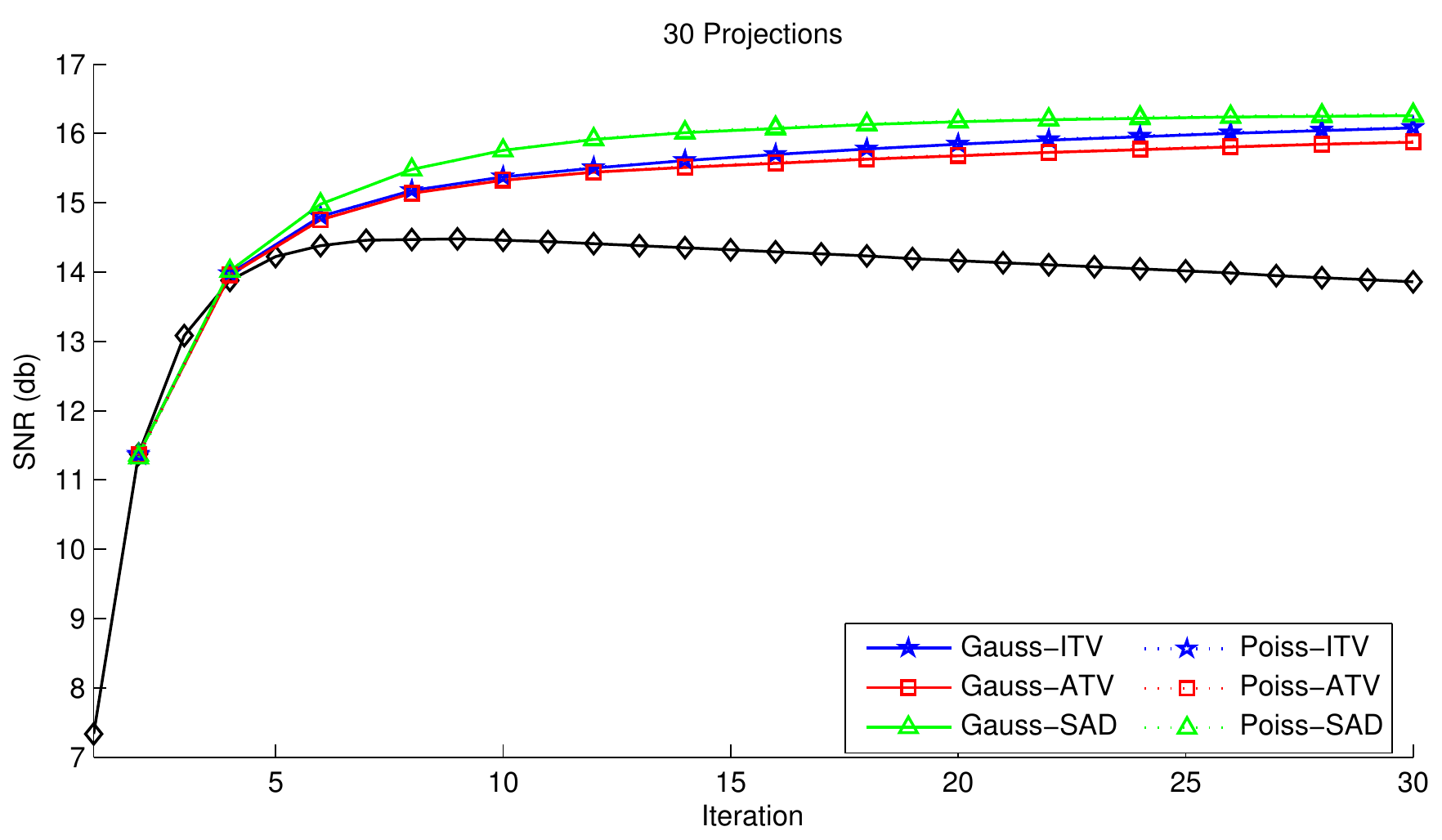}%
\end{minipage}

}

\protect\caption{\textbf{Data Terms and Regularizers Comparison}. Plots show SNR per
iteration  for the Gaussian (\emph{solid} curves) and Poisson (\emph{dashed}
curves)noise models with ITV (\emph{blue}), ATV (\emph{red}), and
SAD (\emph{green}) regularizers. The \emph{black} curve shows the
results for SART. See Sec. \ref{sub:Data-Terms-and-Regularizers-Comparison}.
\label{fig:SNR-per-iter-data-term-and-reg-sart-prox}}
\end{figure*}

We compare the different data terms and regularizers defined in Sec.
\ref{sec:Proximal-Framework}. We solve the tomography proximal operator
(step 3 in Alg. \ref{alg:LADMM-algorithm}) using 2 iterations of
the SART proximal operator (from Table \ref{tab:Prox-Operators}),
using $\alpha=1.99$ with for 15 and 30 projections  (see Sec.
\ref{sub:Iterative-Algorithms-Comparison}). We use $\sigma=0.05$ and $\rho=25$ for 15 projections; $\sigma=0.1$
and $\rho=50$ for 30 projections; and set $\mu=\nicefrac{1}{\rho\Vert K\Vert^{2}}$.
We initialize all methods with uniform volume $x^{(0)}=\mathbf{0}_{n}$.
We estimated the matrix norm $\Vert K\Vert$ using the power method.
Fig. \ref{fig:SNR-per-iter-data-term-and-reg-sart-prox} shows the
results for the three datasets, where we plot against the number of
SART iterations. We note the following: 
\begin{itemize}
\item Using the proximal framework provides significantly better results
than the unregularized iterative methods in Sec. \ref{sub:Iterative-Algorithms-Comparison}.
This is expected since adding a powerful regularizer constrains the
reconstruction to better resemble the ground truth.
\item The Poisson noise model $f_{\text{P}}(\cdot)$ is better than the
Gaussian noise model $f_{\text{G}}(\cdot)$ for the datasets, specially
with more projections. This is consistent with the noise model used
to generate the noisy simulated sinograms, and with the physical noise
model in the real dataset.
\item With more projections, more regularization (higher $\sigma$) produces
better results while for fewer projections less regularization is
sufficient . This is expected because using more projections adds
more constraints (rows in the projection matrix $A$) that need better
regularization to get good results.
\item The SAD regularizer is better for all datasets. 
\end{itemize}

\subsection{Proximal Operators Comparison\label{sub:Proximal-Operators-Comparison}}

We compare the different proximal operators from Sec. \ref{sec:Proximal-Operators}
(SART, ART, OS-SQS, and BICAV) and Table \ref{tab:Prox-Operators}
using our proximal framework in Alg. \ref{alg:LADMM-algorithm}. We
use the best regularizer from Sec. \ref{sub:Data-Terms-and-Regularizers-Comparison},
i.e. SAD regularizer, and both the Poisson and Gaussian noise models.
For the Poisson model, we use $r_{1}$ mapping for SART, and $r_{3}$
for ART, BICAV, and OS-SQS since $r_{1}$ does produce good results
(see Sec. \ref{sub:Poisson-Mapping-Functions-Comparison}). We set
$\sigma=0.05$ and $\rho=50$ for 15 projections; and $\sigma=0.1$
and $\rho=100$ for 30 projections. We set $\mu=\nicefrac{1}{\rho\Vert K\Vert^{2}}$
except for OS-SQS which had to be tuned manually. We note the following:
\begin{itemize}
\item The Poisson model is consistently better than the Gaussian model for
all operators and all datasets.
\item SART proximal operator is generally better than other proximal operators,
and ART and BICAV are quite competitive. 
\item OS-SQS provides the worst performance. We think this has to do with
structure of the update formula in Table \ref{tab:Prox-Operators},
where the \emph{gradient} update $A_{S}^{T}(p_{S}-A_{S}x^{(t)})$
is added to the difference between the current estimate and the input
to the proximal operator $u-x^{(t)}$, where the scaling between the
two terms has to be adjusted properly. That is the reason $\mu$ had
to be carefully tuned to get better results.
\end{itemize}
\begin{figure*}
\center\subfloat[Modified Shepp-Logan]{%
\begin{minipage}[t]{0.33\textwidth}%
\includegraphics[width=1\textwidth]{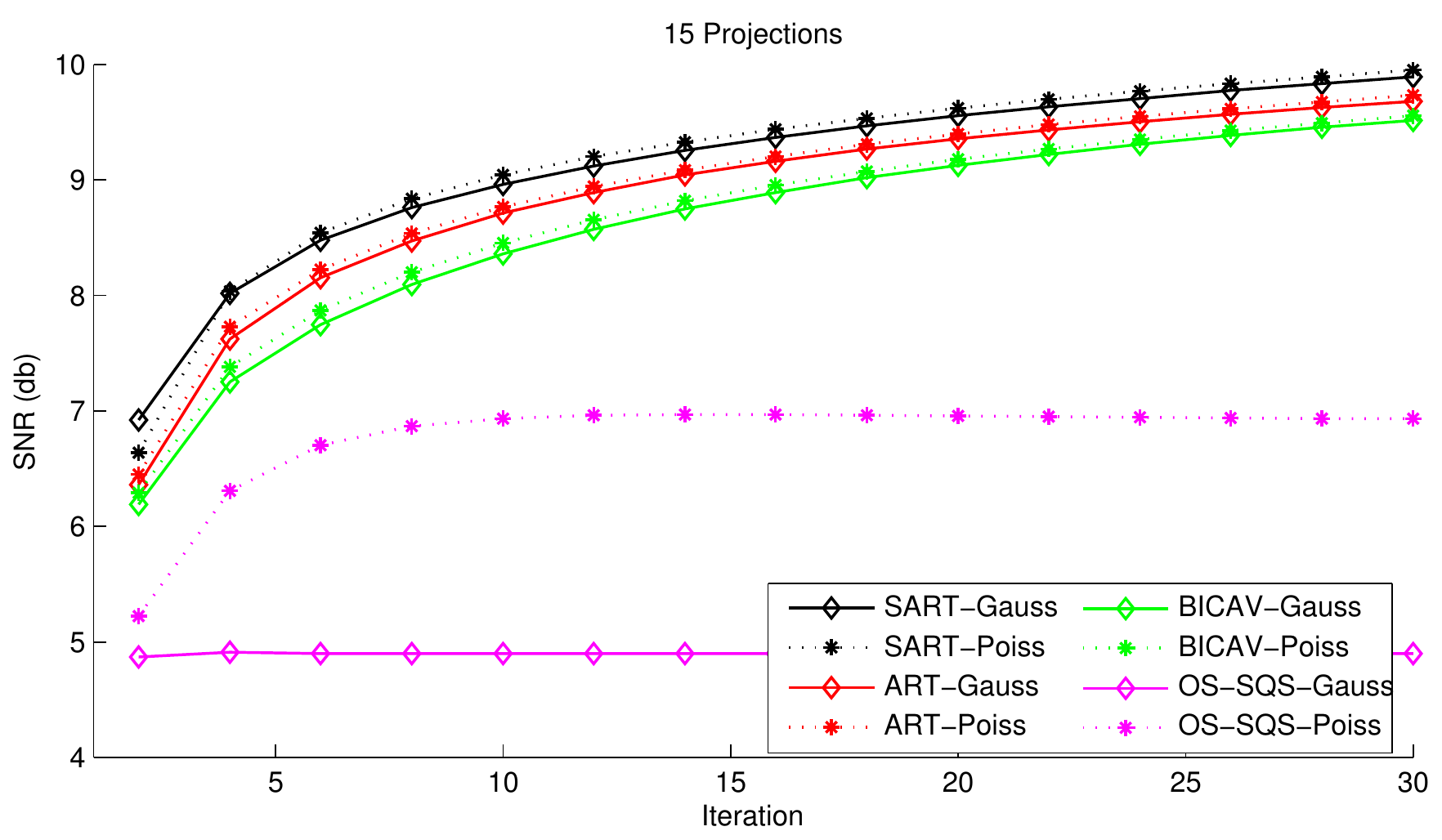}

\includegraphics[width=1\textwidth]{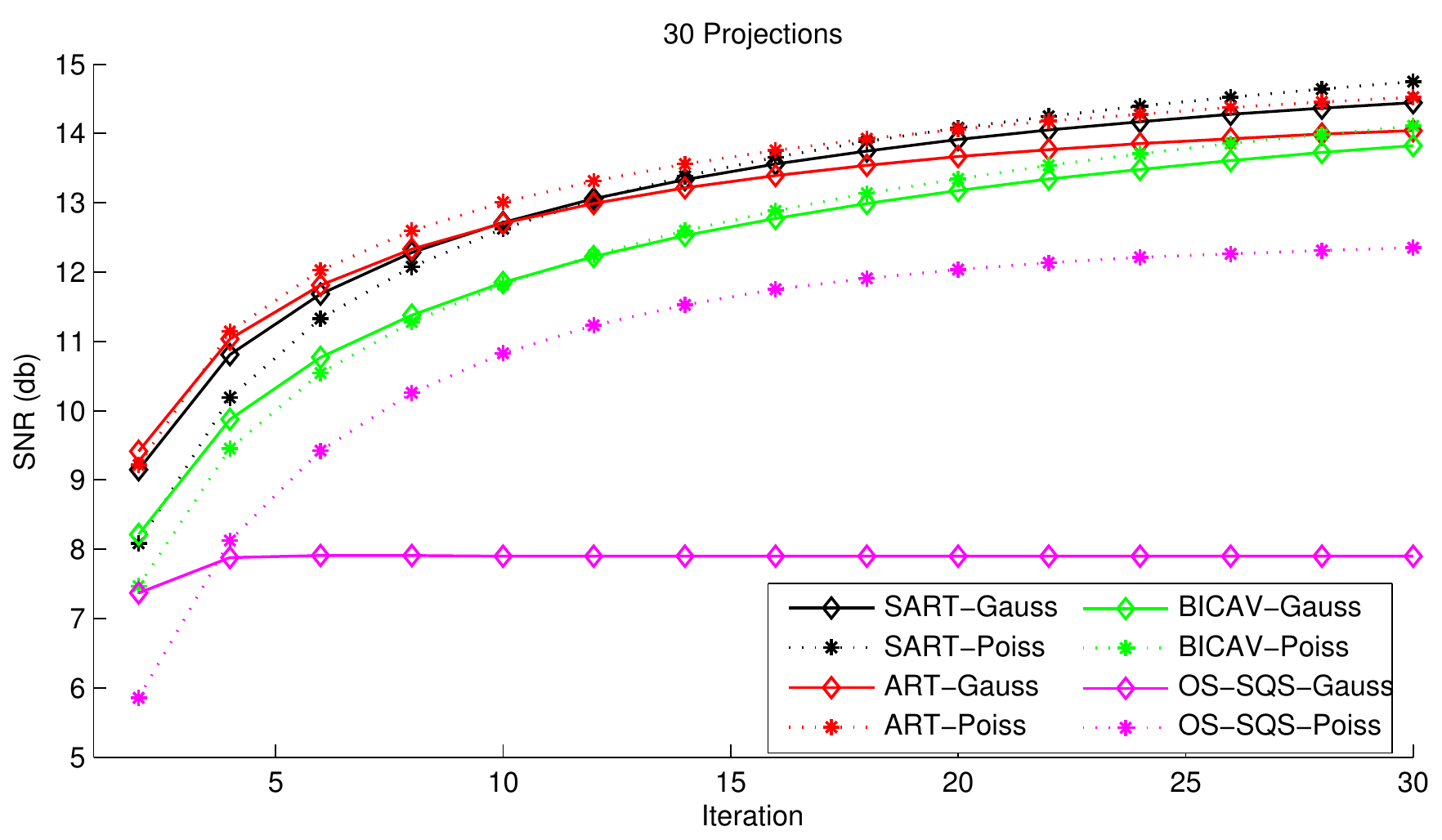}%
\end{minipage}

}\subfloat[NCAT]{%
\begin{minipage}[t]{0.33\textwidth}%
\includegraphics[width=1\textwidth]{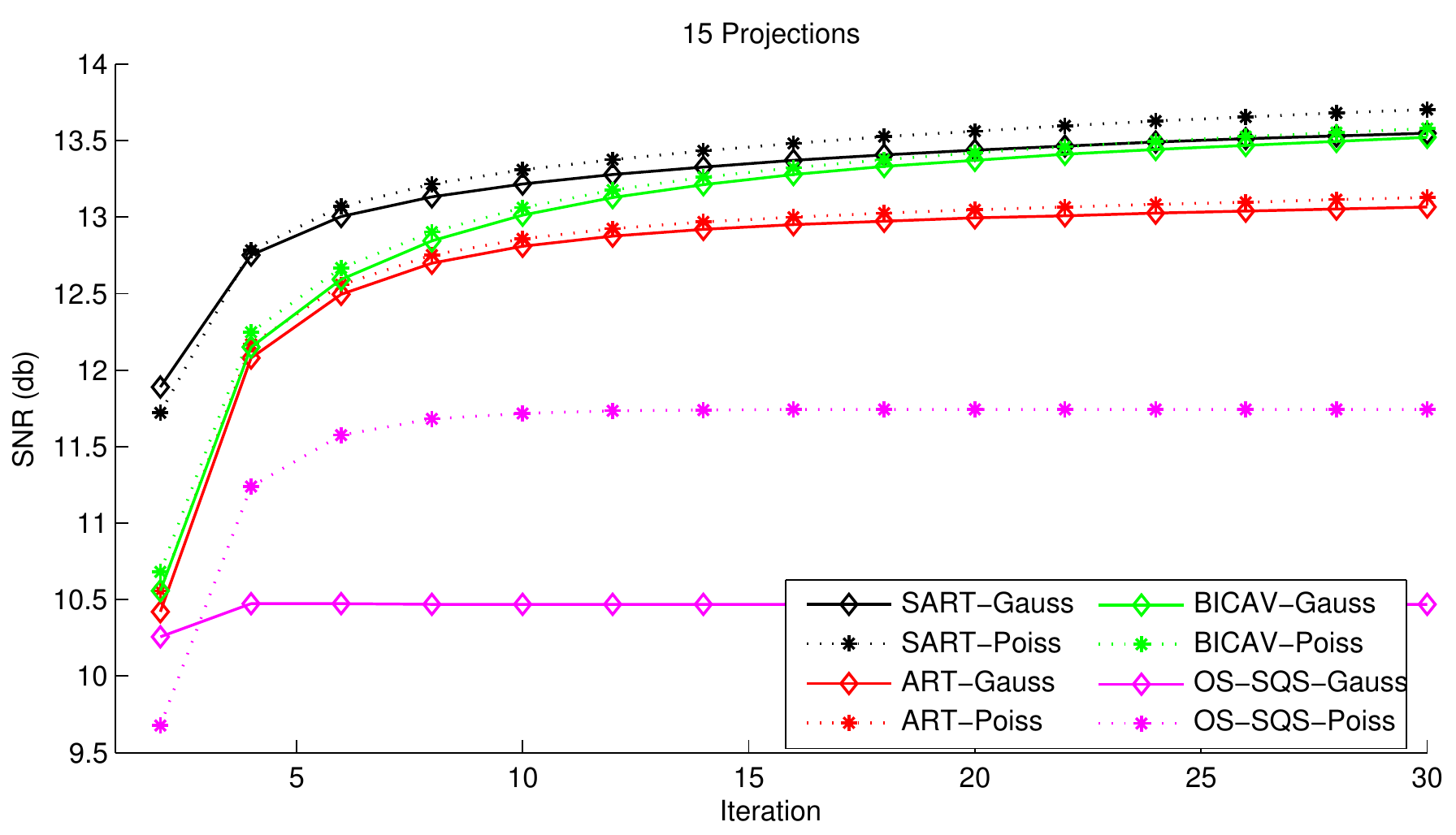}

\includegraphics[width=1\textwidth]{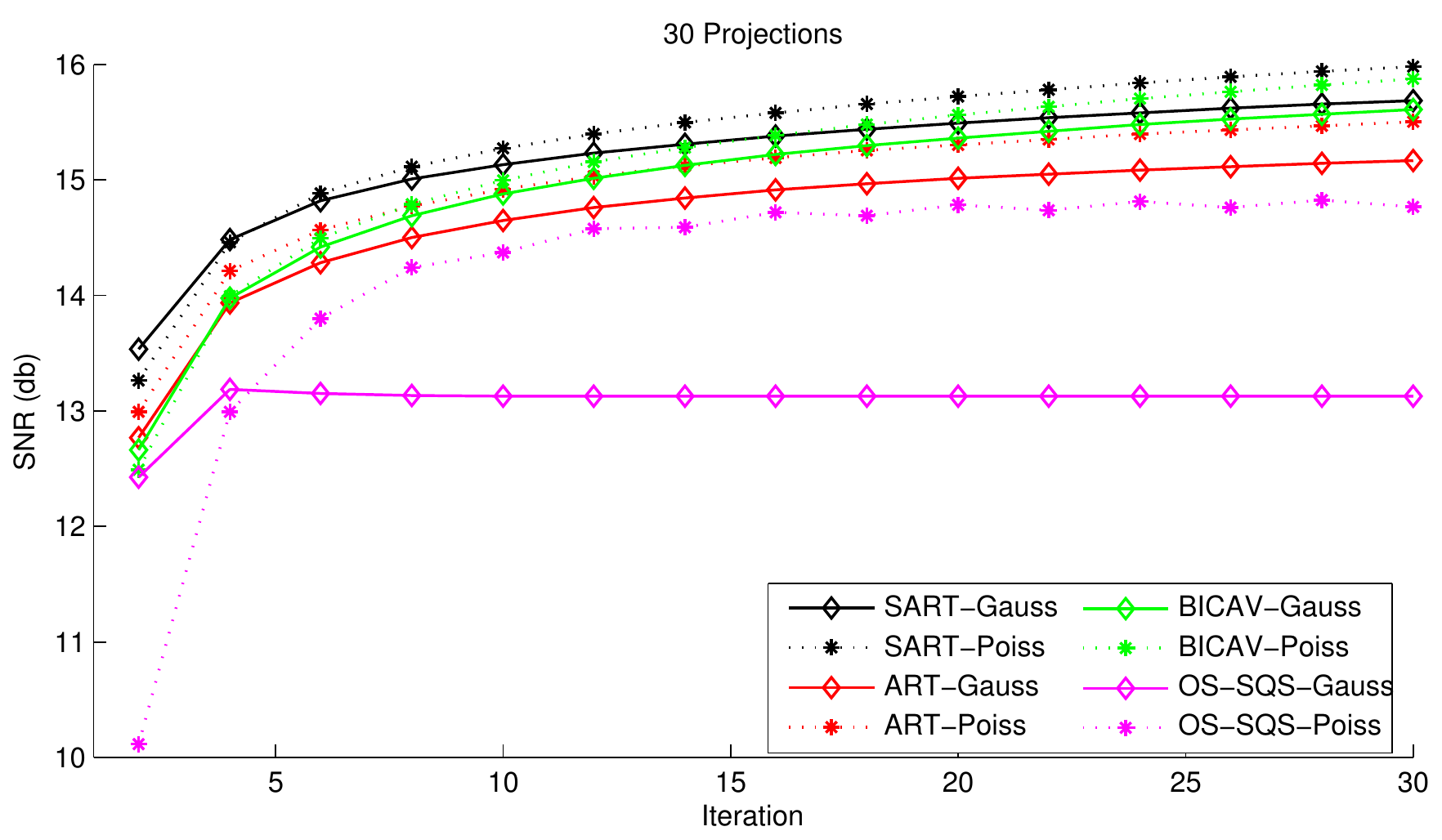}%
\end{minipage}

}\subfloat[Mouse]{%
\begin{minipage}[t]{0.33\textwidth}%
\includegraphics[width=1\textwidth]{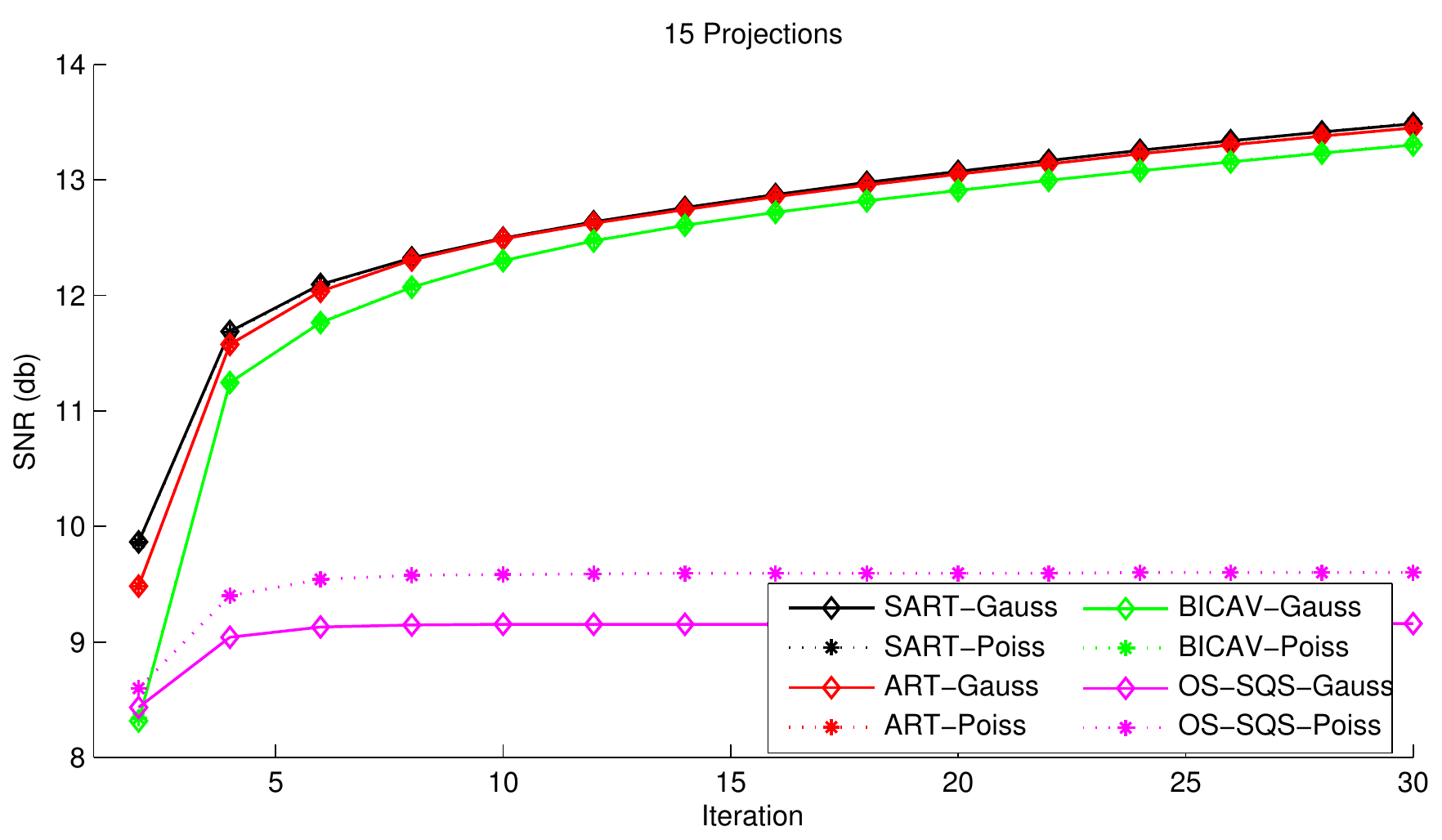}

\includegraphics[width=1\textwidth]{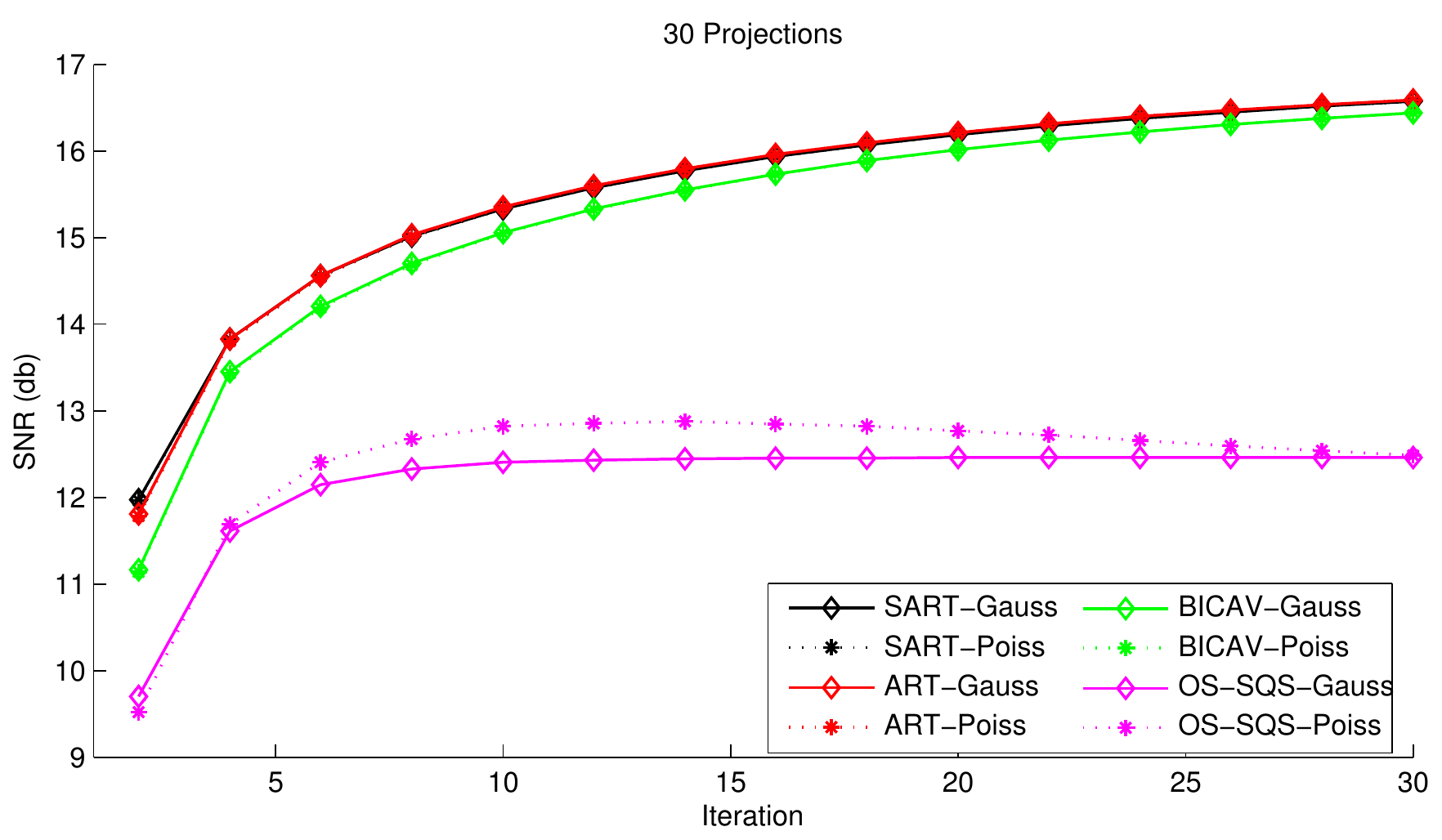}%
\end{minipage}

}

\protect\caption{\textbf{TRex Proximal Operators Comparison}. Plots show SNR per iteration
 for the Gaussian (\emph{solid} curves) and Poisson (\emph{dashed}
curves) noise models with SAD regularizer. See Sec. \ref{sub:Proximal-Operators-Comparison}.
\label{fig:SNR-per-iter-prox-operator-comparison}}
\end{figure*}

\subsection{Comparison to State of the Art\label{sub:Comparison-to-State-of-Art}}

We compare our framework to two state of the art methods: the ADMM-PCG
method of Ramani \etal \cite{ramani2012splitting} and the OS-MOM
method from Kim \etal \cite{kim2015combining} that combines ordered
subsets with momentum. We don't compare to the method of Nien \etal
\cite{nien2015fast} because the authors indicate that the performance
is closely matched by the OS-MOM method and is quite similar.

The ADMM-PCG minimizes a combination of a WLS data term and regularization
term 
\[
\min_{x}\frac{1}{2}\Vert p-Ax\Vert_{W}^{2}+\lambda\sum_{r=1}^{n}\kappa_{r}\Vert Rx_{r}\Vert_{1}
\]
 where $\kappa_{r}\in\mathbb{R}$ are spatial weights that govern
the spatial resolution in the reconstruction and $Rx_{r}\in\mathbb{R}^{d}$
is the vector of wavelet decomposition at voxel $x_{r}$. It uses 2 iterations of PCG to solve the proximal operator of the
data term, as opposed to our framework that uses SART, ART, ... etc.
We run the Matlab code available online from the author as part of
the IRT toolbox%
\footnote{available from http://web.eecs.umich.edu/\textasciitilde{}fessler/code/
}. The default procedure for choosing the parameters didn't work well
with our datasets, so we had to manually tweak the parameters, \cite{ramani2012splitting}.
We set $\nu=2\times10^{5}$, $\lambda=10^{-3}$, and $\mu=10^{-4}$ for 15 projections and $\mu=10^{-3}$
for 30 projections. We use the Wavelet decomposition basis with the
$\ell_{1}$ norm regularizer.

The OS-Mom also minimizes a WLS data term and regularization term 
\[
\min_{x}\frac{1}{2}\Vert p-Ax\Vert_{W}^{2}+\beta\sum_{r=1}^{n}\psi(\nabla x_{r})
\]
 where $\psi(\cdot)$ is an edge preserving potential function and
$\nabla x_{r}\in\mathbb{R}^{2}$ is the gradient at voxel $x_{r}$. We implemented the Momentum 2 method (Table IV in \cite{kim2015combining})
within the IRT toolbox. We use the settings from the paper for the
regularizer i.e. the Fair potential 
\[
\psi(t)=\frac{\delta^{2}}{b^{3}}\left(\frac{ab^{2}}{2}\left|\frac{t}{\delta}\right|^{2}+b(b-a)\left|\frac{t}{\delta}\right|+(a-b)\log(1+b\left|\frac{t}{\delta}\right|)\right)
\]
with $\delta=10$, $a=0.0558$, and $b=1.6395$, and bit reversal
for subset ordering. We tweaked $\beta$ and $M$ \cite{kim2015combining}
to get good performance. We set $\beta=0.05$ and $M=5$ subsets for 15 projections and
$\beta=0.1$ and $M=10$ subsets for 30 projections. We use relaxation
with parameter $10^{-3}$ to help the convergence. We also compare
to plain OS method without momentum with the same WLS data term and
regularizer as OS-Mom.

\begin{figure*}
\center\subfloat[Modified Shepp-Logan]{%
\begin{minipage}[t]{0.33\textwidth}%
\includegraphics[width=1\textwidth]{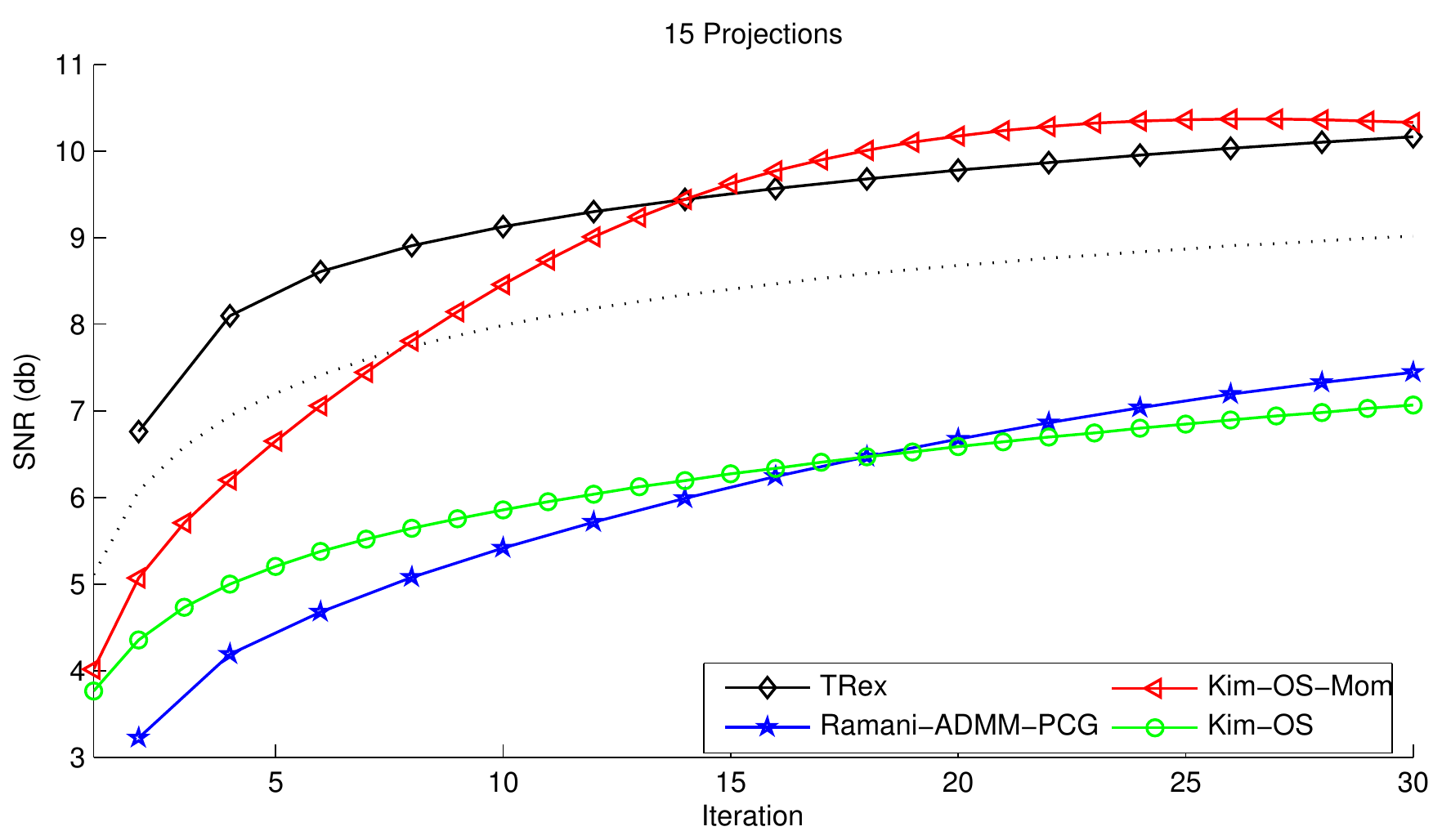}

\includegraphics[width=1\textwidth]{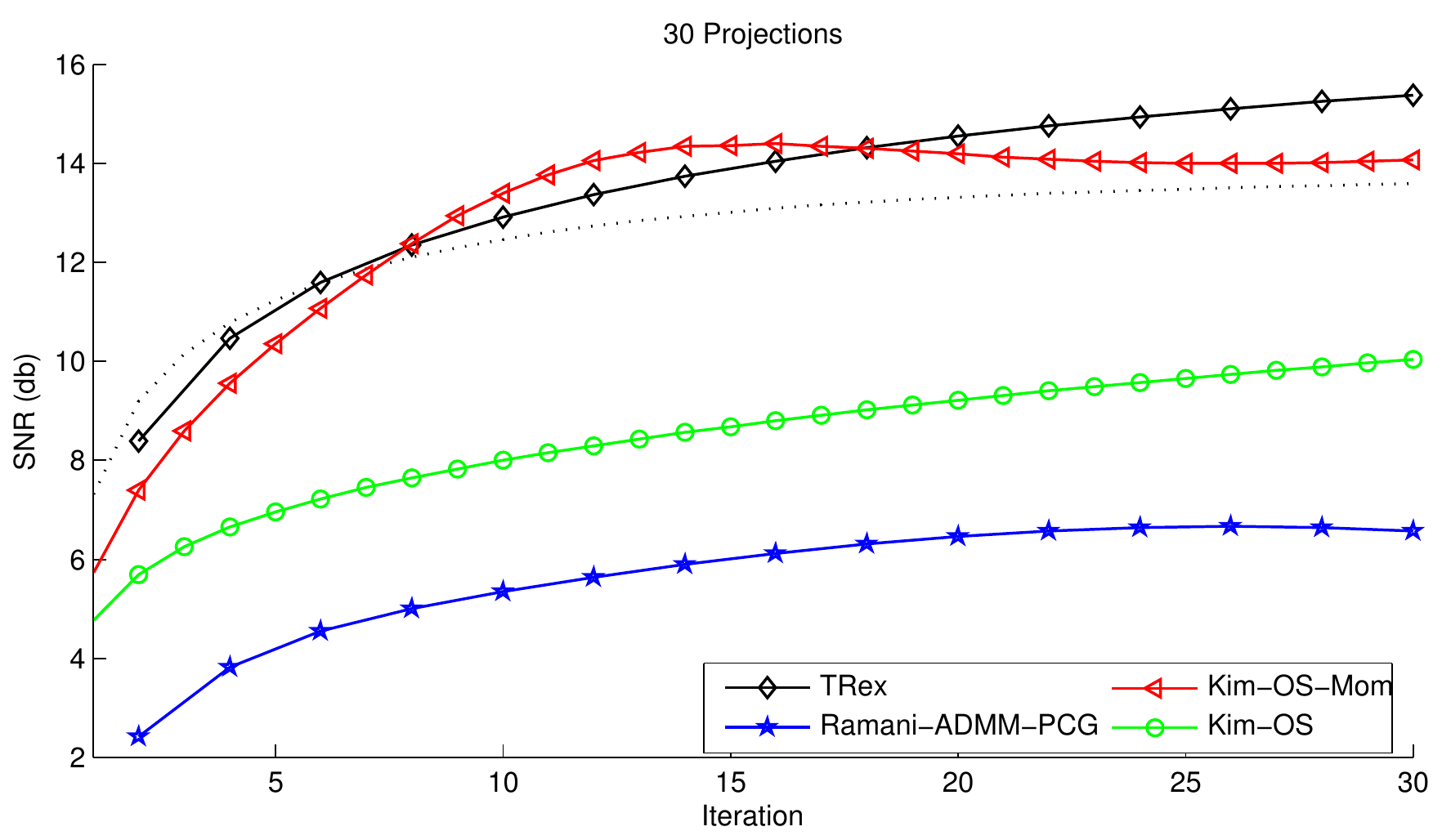}%
\end{minipage}

}\subfloat[NCAT]{%
\begin{minipage}[t]{0.33\textwidth}%
\includegraphics[width=1\textwidth]{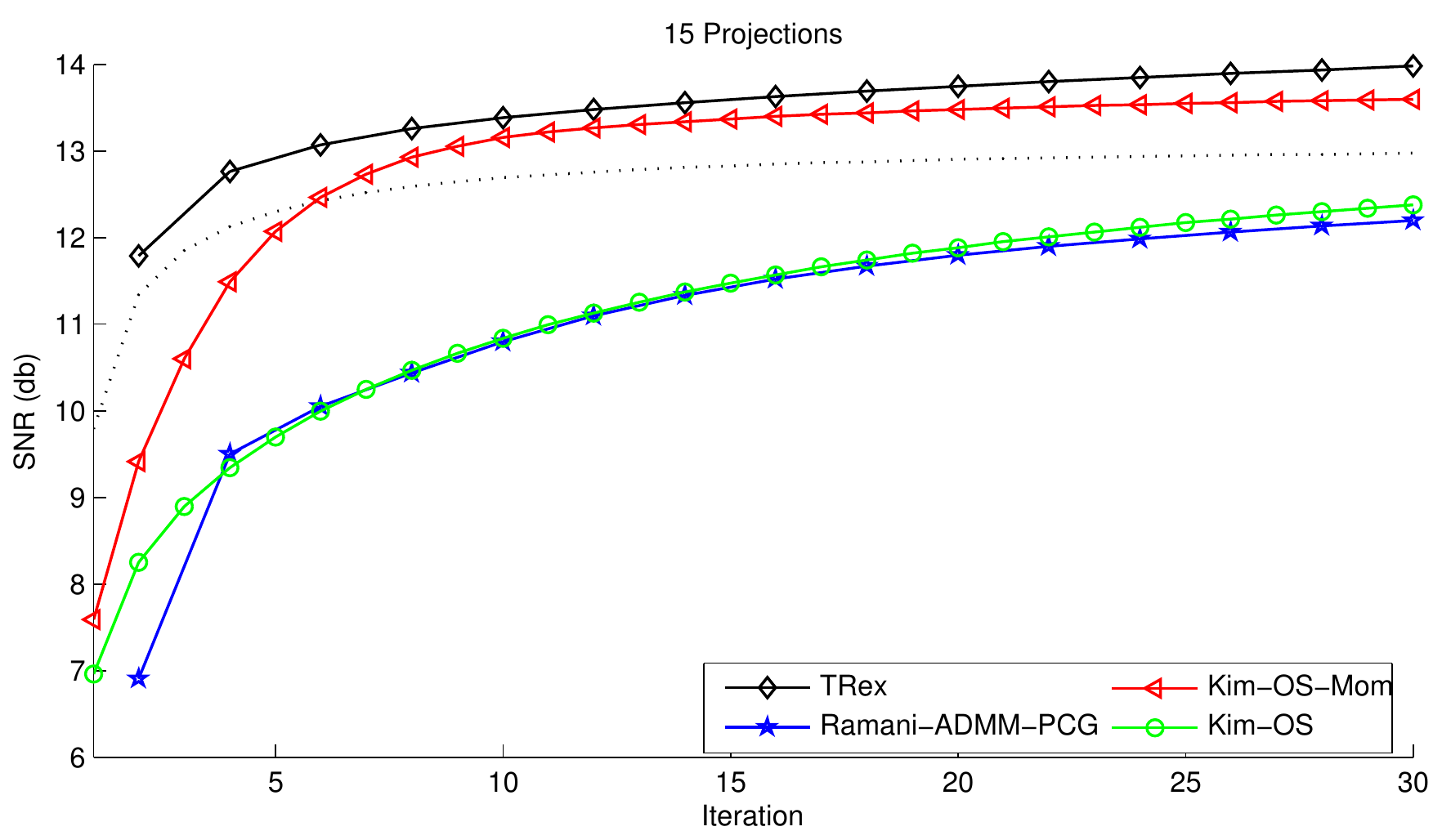}

\includegraphics[width=1\textwidth]{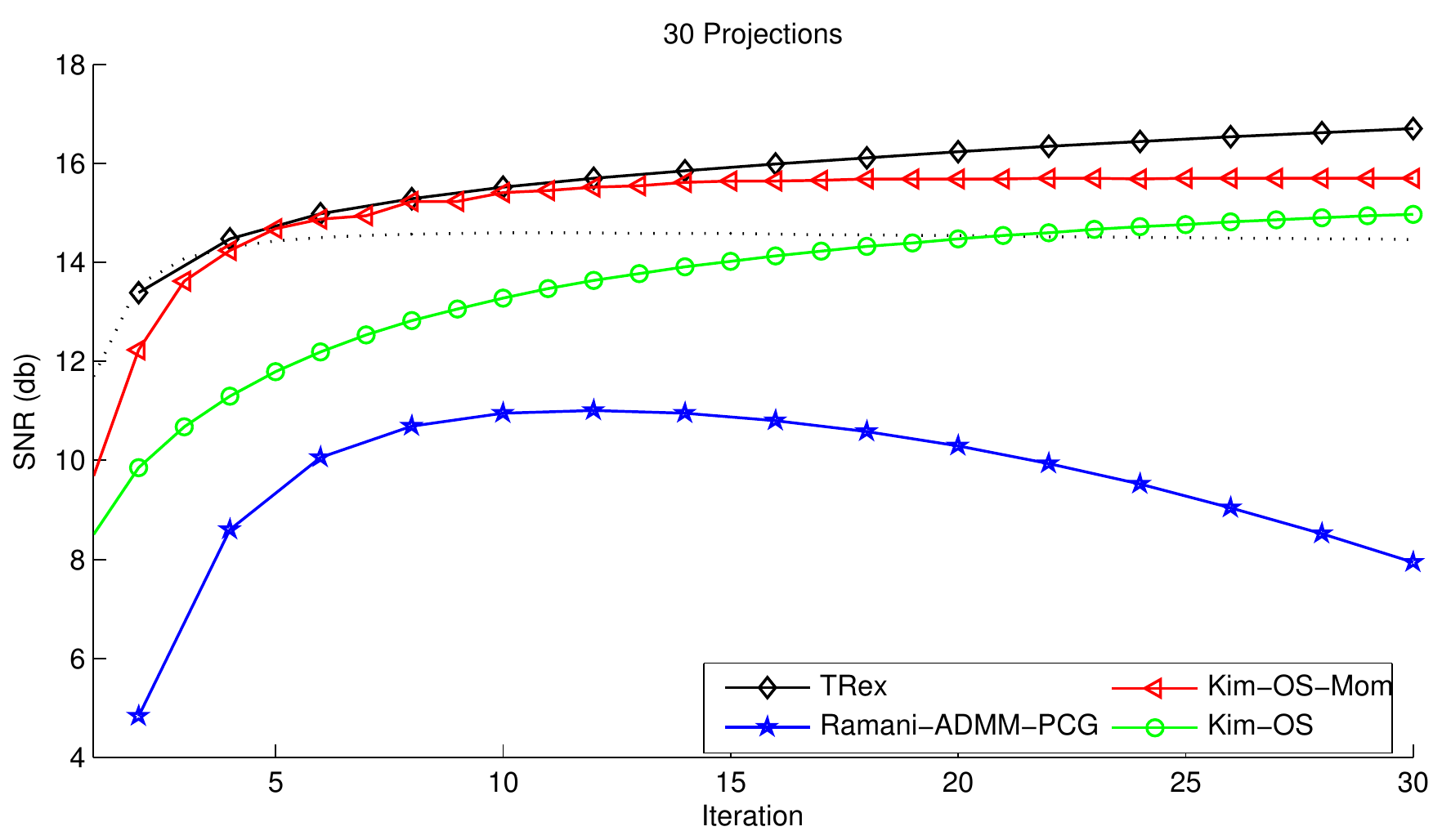}%
\end{minipage}

}\subfloat[Mouse]{%
\begin{minipage}[t]{0.33\textwidth}%
\includegraphics[width=1\textwidth]{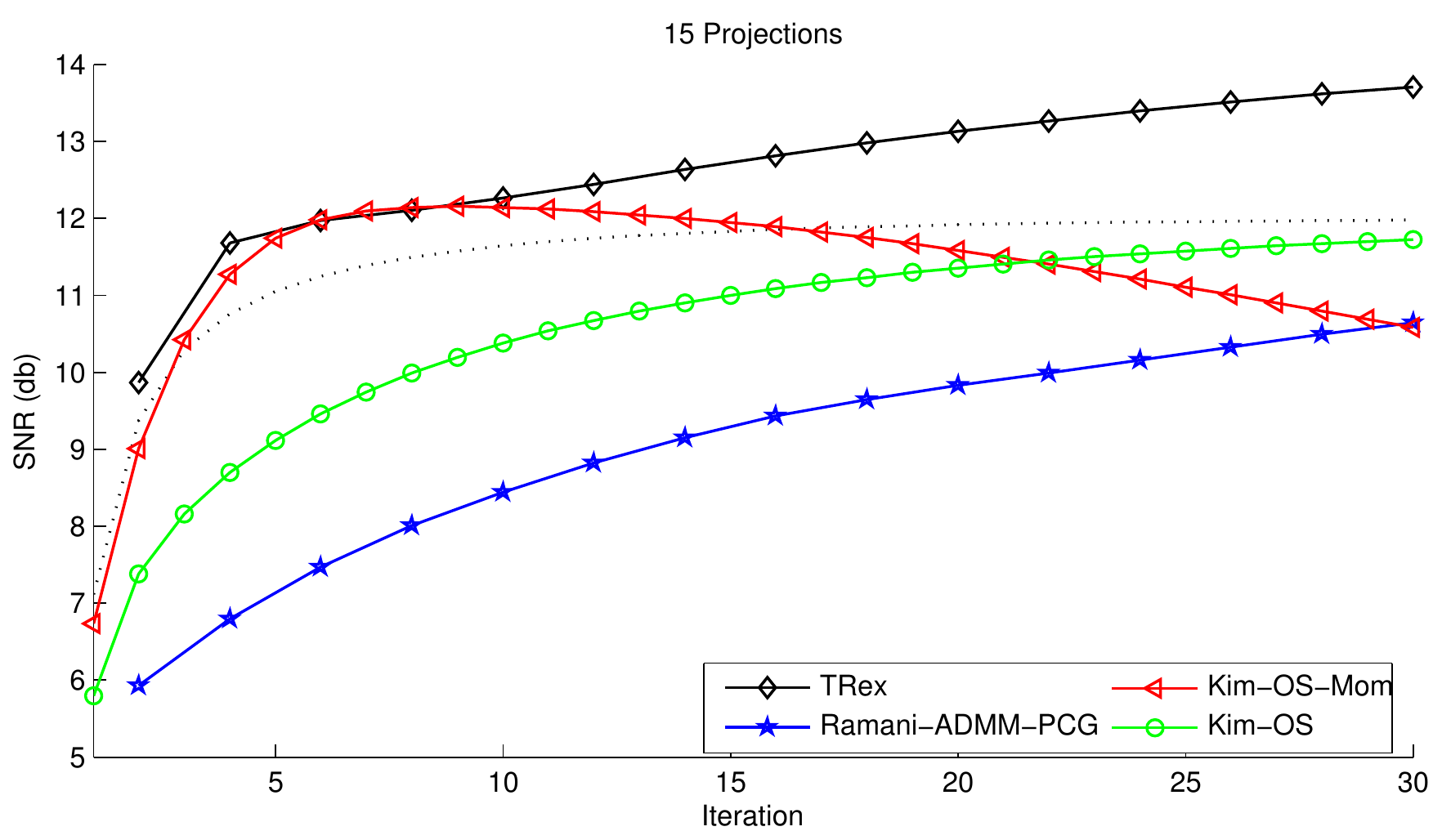}

\includegraphics[width=1\textwidth]{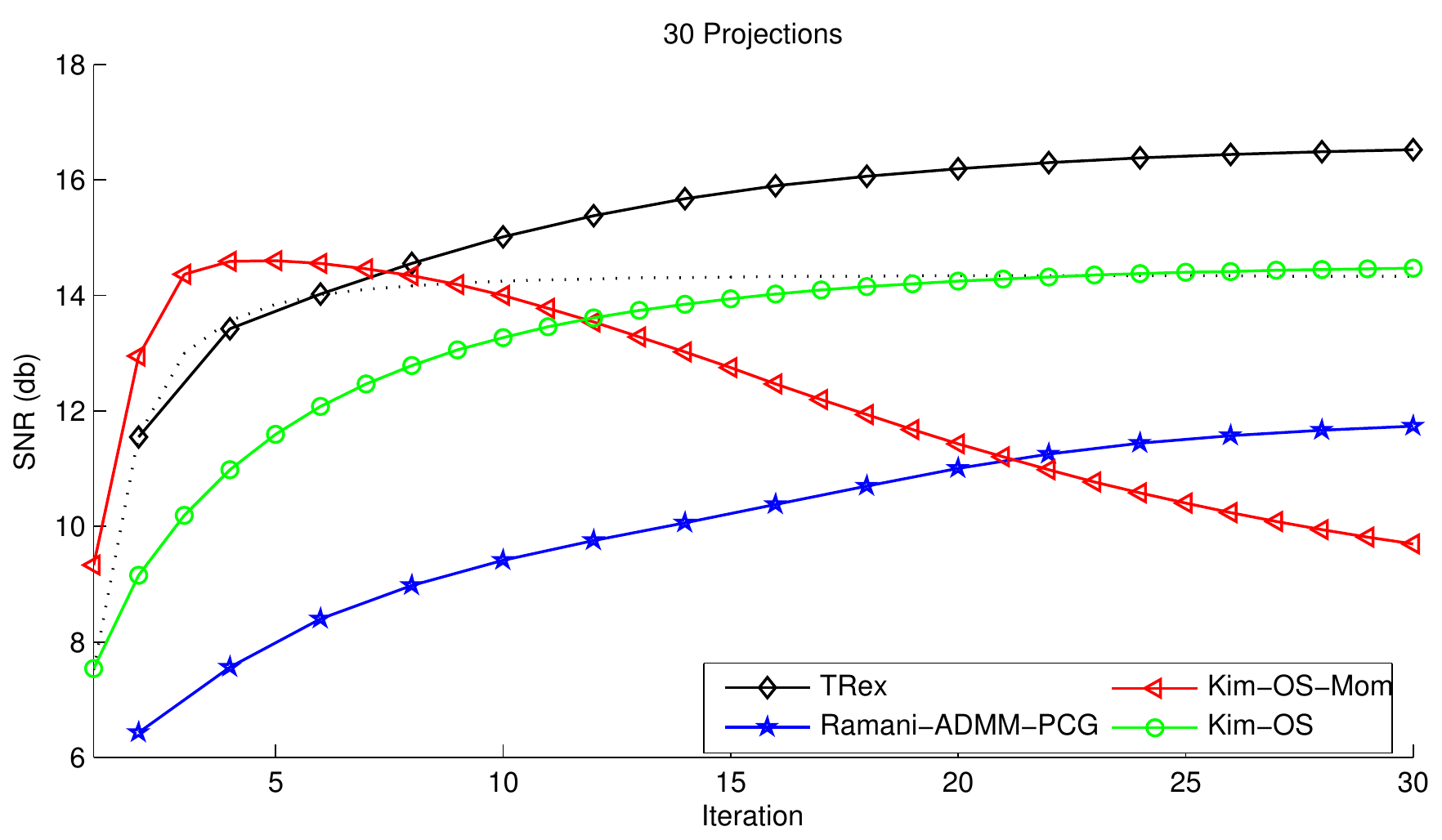}%
\end{minipage}

}

\protect\caption{\textbf{TRex Comparison to Sate of the Art}. Plots show SNR per iteration.
The TRex framework uses SART with Poisson model and SAD regularizer.
The dotted curve shows the baseline plain SART (from Sec. \ref{sub:Iterative-Algorithms-Comparison}).
See Sec. \ref{sub:Comparison-to-State-of-Art}. \label{fig:SNR-per-iter-state-of-art-comp}}
\end{figure*}

\begin{figure*}
\center\subfloat[Modified Shepp-Logan]{\includegraphics[width=1\textwidth]{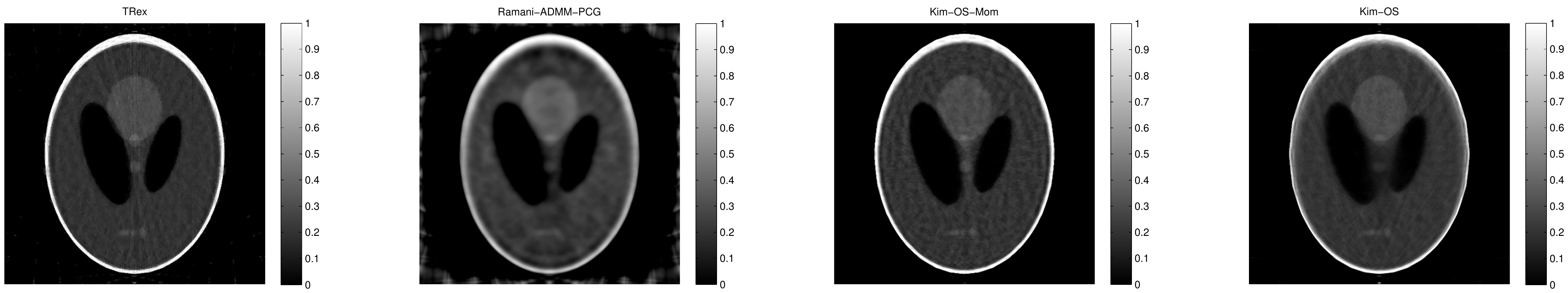}

}\\
\subfloat[NCAT]{\includegraphics[width=1\textwidth]{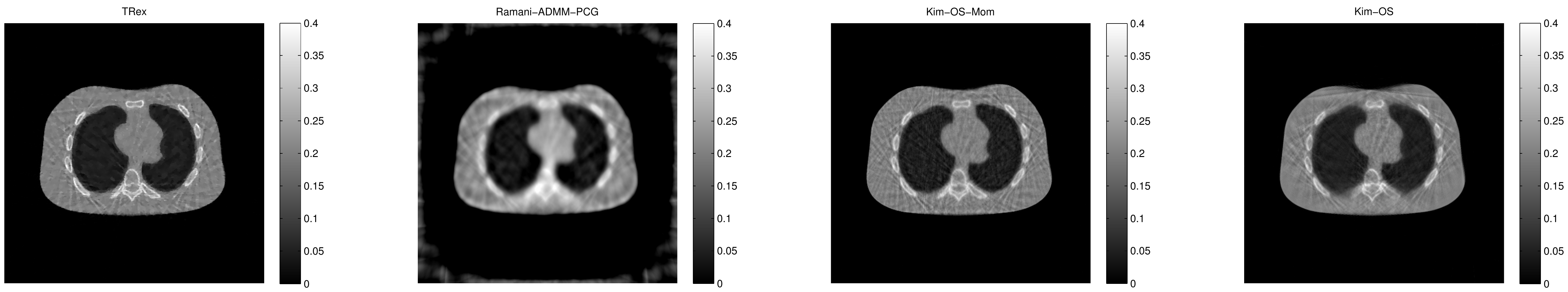}

}\\
\subfloat[Mouse]{\includegraphics[width=1\textwidth]{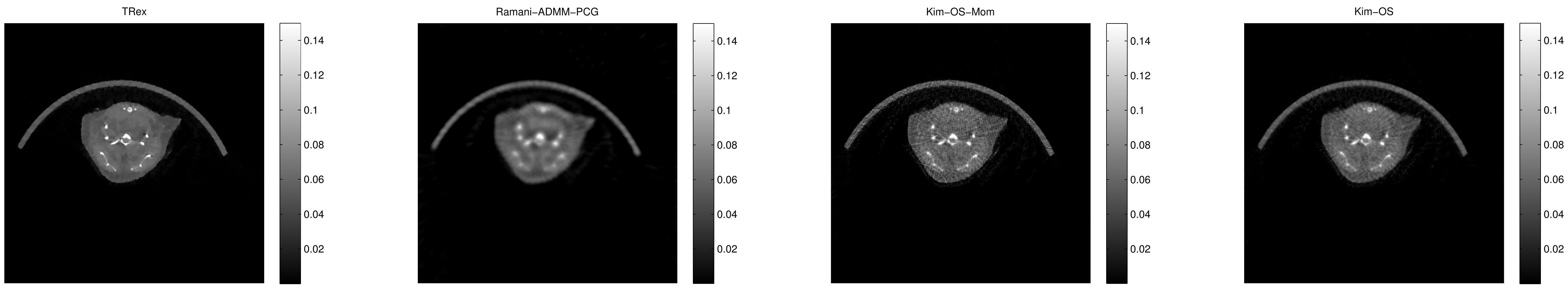}

}

\protect\caption{\textbf{TRex Comparison to Sate of the Art}. Reconstruction results
for 30 projections after 30 iterations. TRex uses SART with Poisson
model and SAD regularizer. See Sec. \ref{sub:Comparison-to-State-of-Art}.
\label{fig:reconstruction-state-of-art-comp}}
\end{figure*}

Fig. \ref{fig:SNR-per-iter-state-of-art-comp} shows a comparison
with these two algorithm. The TRex uses the SART proximal operator
with Poisson noise model and $r_{1}$ mapping and SAD regularizer.
We set $\sigma=0.05$ and $\rho=25$ for 15 projections and $\sigma=0.1$
and $\rho=50$ for 30 projections, and set $\mu=\nicefrac{1}{\rho\Vert K\Vert^{2}}$.
We initialize all methods with a uniform volume $x^{(0)}=\mathbf{0}_{n}$.
Fig. \ref{fig:reconstruction-state-of-art-comp} shows sample reconstruction
results for 30 projections after 30 iterations. We note the following:
\begin{itemize}
\item The ADMM-PCG method's performance is quite bad with these datasets.
They require a lot of tweaking to get them to work correctly since
the automated estimation methods described in \cite{ramani2012splitting}
didn't work. Moreover, the method seems very sensitive to the values
of the parameters, and thus is harder to tweak.
\item The OS-Mom method indeed accelerates the convergence of the OS method
at early iterations \cite{kim2015combining}. However, its performance
is not consistent across datasets, where some times it is good and
most of the time the SNR starts decreasing after a while, even with
relaxation.
\item TRex with SART and SAD consistently performs better and ends up with
higher SNR than ADMM-PCG or OS-Mom. Moreover, it is easier to tweak
and not very sensitive to the choice of parameters. Note that for
NCAT and Mouse, using TRex, we get SNR with 15 projections that equals
the SNR we get with 30 projections using plain SART.
\end{itemize}

\section{Conclusions\label{sec:Discussion-and-Conclusion}}

We presented TRex, a flexible proximal framework for robust tomography
reconstruction in sparse view applications. TRex relies on using iterative
methods, e.g. SART, for directly solving the tomography proximal operator.
We first compare the famous tomography iterative solvers, and then
derive proximal operators for the best four methods. We then show
how to use TRex to solve using different noise models (Gaussian and
Poisson) and using different powerful regularizers (ITV, ATV, and
SAD). We show that TRex outperforms state of the art methods, namely
ADMM-PCG \cite{ramani2012splitting} and OS-Mom \cite{kim2015combining},
and is easy to tune. We conclude that SART---even though is not guaranteed
to converge---offers the best tomography solver for sparse view applications,
followed closely by ART and BICAV. 

We plan to extend this work in several ways: (a) study how to incorporate
momentum acceleration into SART as in \cite{kim2015combining}; (b)
study how to use preconditioners with SART such as the Fourier-based
cone filter preconditioners \cite{fessler1999conjugate}; (c) study
other applications such as low-dosage X-ray tomography, which changes
the nature of the measurement noise \cite{xu2014quantifying}; and
(d) implement and apply TRex to 3D cone beam reconstruction and compare
to other famous packages such as RTK \cite{mory2012ecg}.

\section*{Acknowledgments}

This work was supported by KAUST baseline and research center funding.

\bibliographystyle{ieeetr}
\addcontentsline{toc}{section}{\refname}\bibliography{papers}

\end{document}